\documentclass[10pt]{article}
\usepackage{amssymb,amsmath,amsthm,latexsym, bbm}
\usepackage{enumerate}
\usepackage{color}
\usepackage{graphicx}
\usepackage{mathrsfs}
\usepackage[mathscr]{euscript}
\usepackage{amscd}

\newtheorem{thm}{Theorem}[section]

\newtheorem{lemma}[thm]{Lemma}
\newtheorem{prop}[thm]{Proposition}
\newtheorem{defn}[thm]{Definition}

\newtheorem{remark}[thm]{Remark}

\numberwithin{equation}{section}

\def\Ex {{\mathbb E}}

\def\R {{\mathbb R}}

\def\F {{\mathcal F}}
\def\P {{\mathbb P}}

\def\D {{\mathcal D}}

\def\vH {{\mathbf H}}

\def\v {{\mathbf v}}
\def\u {{\mathbf u}}
\def\w {{\mathbf w}}
\def\V {{\mathbf V}}
\def\d {{\mathbf d}}
\def\b{{\mathbf b}}

\def\dH {{\mathbb H}}

\newcommand{\norm}[1]{\lVert#1\rVert}
\newcommand{\qv}[1]{\langle#1\rangle}
\def\wt{\widetilde}

\def\1{{\mathbf 1}}

\def\<{{\langle}}
\def\>{{\rangle}}
\def\eps{\varepsilon}

\def\proof{{\medskip\noindent {\bf Proof. }}}

\def\qed{{\hfill $\square$ \bigskip}}

\makeatletter
\def\tank#1{\protected@xdef\@thanks{\@thanks
		\protect\footnotetext[0]{#1}}}
\def\bigfoot{
	
	\@footnotetext}
\makeatother

\begin{document}
		
\title
{\bf Global well-posedness of stochastic nematic liquid crystals with random initial and random boundary conditions driven by multiplicative noise \thanks{This work was partially
supported by NNSF of China(Grant No. 11401057, 11801283),   Natural Science Foundation Project of CQ  (Grant No. cstc2016jcyjA0326),
Fundamental Research Funds for the Central Universities(Grant No.Grant No. 2018CDXYST0024, 63181314) and China Scholarship Council (Grant No.:201506055003).} }
		\author{
Lidan Wang
\thanks{ Nankai University, Tianjin, 300071, P.R. China}
\tank{lidanw.math@gmail.com}
\qquad
Jiang-Lun Wu
\thanks{Swansea University, Bay Campus
Fabian Way, Swansea SA1 8EN
Wales/UK }
\tank{E-mail: j.l.wu@swansea.ac.uk.}
\qquad
Guoli Zhou
\thanks{Corresponding author. Chongqing University, Chongqing, 400044, P.R. China }
\tank{E-mail:zhouguoli736@126.com.}
}

\date{}
\maketitle

\begin{abstract}
	The flow of nematic liquid crystals can be described by a highly nonlinear stochastic hydrodynamical model, thus is often influenced by random fluctuations, such as uncertainty in specifying initial conditions and  boundary conditions. In this article, we consider the $2$-D stochastic nematic liquid crystals with the velocity field perturbed by affine-linear multiplicative white noise, with random initial data and random boundary conditions. Our main objective is to establish the global well-posedness of the stochastic equations under certain sufficient Malliavin regularity of the initial conditions and the boundary conditions. The Malliavin calculus techniques play important roles in proving the global existence of the solutions to the stochastic nematic liquid crystal models with random initial and random boundary conditions. It should be pointed out that the global well-posedness is also true when the stochastic system is perturbed by the noise on the boundary.
\end{abstract}	

\bigskip


\bigskip


\bigskip
	
%%%%%%%%%%%%%%%%%%%%%%
\section{Introduction}
%%%%%%%%%%%%%%%%%%%%%%	
The liquid crystal is an intermediate state of a matter, which possesses some typical properties of a liquid as well as some crystalline properties. One can observe the flow of nematic liquid crystals as slowly moving particles where the alignment of particles and the velocity of the fluid sway each other. The history of the hydrodynamic theory for liquid crystals traces back to 1960's, Ericksen \cite{E} and Leslie \cite{L} expanded the continuum theory to design the dynamics of the nematic liquid crystals. The so-called Ericksen-Leslie system is well designed for describing many special flows for the materials, especially for those with small molecules, and is widely applied in the engineering and mathematical communities for studying liquid crystals.

Later on, the most fundamental formulation of dynamical system describing the orientation as well as the macroscopic motion for the nematic liquid crystals was introduced by Lin-Liu \cite{LL}:
\begin{align*}
	&d\v+[(\v\cdot\nabla)\v-\mu\Delta\v+\nabla p]dt=-\lambda\nabla\cdot(\nabla\d\odot\nabla\d)dt, \nabla\cdot\v=0,\\
	&d\d+(\v\cdot\nabla)\d dt=\gamma(\Delta\d+|\nabla\d|^2\d)dt, |\d|^2=1.
\end{align*}
In order to avoid the difficulties arising from bounding the nonlinear gradient term in the second equation for the orientation field, as suggested in  Lin-Liu \cite{LL}, one can use the Ginzburg-Landau approximation to ease the constraint $|\d|^2=1$, and the corresponding approximation energy is
$$\int_D\left[\frac{1}{2}|\nabla\d|^2+\frac{1}{4\beta^2}(|\d|^2-1)^2\right]dx.$$
Then one arrives at the following approximating system
\begin{align*}
	&d\v+[(\v\cdot\nabla)\v-\mu\Delta\v+\nabla p]dt=-\lambda\nabla\cdot(\nabla\d\odot\nabla\d)dt, \nabla\cdot\v=0,\\
	&d\d+(\v\cdot\nabla)\d dt=\gamma\left(\Delta\d-\frac{1}{\eta^2}(|\d|^2-1)\d\right)dt.
\end{align*}
The above system can be viewed as the simplest mathematical model keeping the most important mathematical structure as well as essential difficulties of the original Ericksen-Leslie system (see \cite{LL}). This deterministic system with Dirichlet boundary conditions has been well studied in a series of work both theoretically (see \cite{LL,LL2})and numerically.

Along with the developments of deterministic system, the random case has also drawn a lot interests in recent years. In the papers \cite{BHR1,BHR2}, Brze\'zniak-Hausenblas-Razafimandimby have studied the nematic crystal flow model perturbed by multiplicative Gaussian noise and give the global well-posedness for the weak and strong solutions in $2$-D case. For the pure jump noise case in $2$-D, Brze\'zniak-Manna-Panda \cite{BMP} have obtained the global well-posedness for the martingale solution. A weak martingale solution result is also established for three dimensional stochastic nematic liquid crystals with pure jump noise in \cite{BMP}.

As far as we know, the present work is the first attempt to study stochastic nematic liquid crystal equations with random initial and random boundary conditions, especially when the orientation field is perturbed by the noise on the boundary. Our motivation firstly derives from the limitation of predicting dynamical behavior in nonlinear systems due to uncertainty in initial data, which has been widely investigated (see \cite{GHs}). The related study has drawn a lot attention in the geophysical community (see \cite{OSB, Pa, PH}). Our main result in this article implies that each stationary point of the present stochastic model generates a pathwise anticipating stationary solution of the Stratonovich stochastic equations. Another motivation of our work is, near stationary solutions, multiplicative ergodic theory techniques ensure the existence of local random invariant manifolds which necessarily anticipate the driven noise. One can refer to \cite{DLS, MZ} and related works. Hence, the study of a dynamic characterization of semiflows as well as invariant manifolds will appeal to the analysis of the stochastic nematic liquid crystal equations with anticipating initial date and corresponding random boundary conditions. This can be viewed as a necessary first step in the analysis of the regularity of invariant manifolds.

In this article, we consider in $\mathbf{D}\times\R_+$, where $\mathbf{D}\subset\R^2$ is a bounded domain with smooth boundary, the stochastic version of the nematic liquid crystals flows with random initial and boundary conditions. The model is formulated in the following maner:
\begin{align}
&\label{eqn:1.1}\v_t+(\v\cdot\nabla)\v-\mu\Delta\v+\nabla p+\lambda\nabla\cdot(\nabla\d \odot\nabla\d)=\sum_{k=1}^{\infty}\sigma_k\v\circ\dot W_k+\sigma_0\dot{\wt{W_0}},\\
&\nabla\cdot\v(t)=0,\\
&\label{eqn:1.3}\d_t+(\v\cdot\nabla)\d-\gamma\left(\Delta\d-\frac{1}{\eta^2}(|\d|^2-1)\d\right)=0.
\end{align}
The unknowns are the fluid velocity field $\v=(v^1,v^2)\in\R^2$, the averaged macroscopic/continuum molecular orientation field $\d=(d^1,d^2,d^3)\in\R^3$, and the pressure function $p(x,t)$, where $\mu,\lambda, \gamma$ are positive constants and stand for viscosity, the competition between kinetic and potential energies, and macroscopic elastic relaxation time for $\d$ respectively.  The operation $[\nabla\d \odot\nabla\d]_{ij}$ yields a $2\times2$ matrix whose  entry is given by the following
$$[\nabla\d \odot\nabla\d]_{ij}=\sum_{k=1}^{3}\partial_{x_i}d^k\partial_{x_j}d^k,\ i,j=1,2.$$

For the stochastic term, $\{W_k(t)_{t\in[0,T]}\}_{k\geq1}$ is a sequence of independent, identically distributed one dimensional Brownian motions which are also independent of a space-time noise $\wt W_0(t,x)$. The space-time noise $\tilde{W}_{0}(t,x)$ is a Brownian in the time variable $t\in \mathbb{R}^{+}$ and smooth in the space variable $x\in \mathbf{D},$ where $\dot{W_k},\dot{\wt W_0}$ are the heuristic time derivatives. The random forces are all defined on the same completely filtered Wiener space $(\Omega,\F,(\F_t)_{t\geq0},\P)$. We also assume that $\sigma_0\in\R$ and $\sum_{k=1}^{\infty}\sigma_k^2<\infty$.

We supplement the stochastic nematic liquid crystals equations with the following random initial and boundary conditions:
\begin{align}
\label{IC}&\v(x,0)=R_{\nu}(x),\ \d(x,0)=R_d(x),\ \text{ for }x\in \mathbf{D};\tag{IC}\\
\label{BC}&\v(x,t)=0,\ \d(x,t)=R_d(x),\ \text{ for }(x,t)\in\partial \mathbf{D}\times\R^+,\tag{BC}
\end{align}
where the initial conditions $R_{\nu},R_d$ are  $\F\otimes\mathcal B(\mathrm{V})\times \mathcal B(\mathbb{H}^{2})$-measurable random fields on $\mathbf{D}$ with $\mathrm{V}$ and $\mathbf{H}^{2}$ defined in Section 2.

%Our main result in this paper is the uniqueness of ergodic invariant measures for the model \eqref{eqn:1.1}-\eqref{eqn:1.3}.
%\begin{thm}
%	Under certian conditions for $(\v_0,\d_0)$, there exists a unique ergodic invariant measure for the stochastic nematic liquid crystals flows.
%\end{thm}
Compared to $2$-D Navier-Stokes equations \cite{MZ}, the stochastic nematic liquid crystal model is more complicated since it is highly nonlinear and coupled with non-homogenous boundary conditions. This causes essential difficulties in obtaining energy estimates and moment estimates, see Proposition \ref{prop:2.5}, Proposition \ref{prop:2.6}, Theorem \ref{thm:3.2} and Theorem \ref{prop:3.4}. To overcome these difficulties, we take the advantage of the special geometric structure of the nematic liquid crystal equation to obtain the adjoint estimate of $\nabla\cdot(\nabla\d\odot\nabla\d)$ and $(\v\cdot\nabla)\d$. This is vital  to establish the \textit{a priori} estimates for the solutions. In \cite{MZ}, the authors imposed sufficiently smooth initial condition to ensure the Malliavin differentiability of the weak solutions. Further, by utilising approximation argument, the authors are able to remove the smoothness assumption on the initial condition. Here the aim of our article is to obtain the global well-posedness of the strong solution to the nematic liquid crystal equation with random initial and random boundary conditions. Then it will further implies the global well-posedness of the weak  solution to the nematic liquid crystal equation with random initial and random boundary conditions. In our procedure, see Theorem \ref{thm:3.2} and Theorem \ref{prop:3.4}, we observe that the strong solution is sufficiently regular to ensure nice bounds for the nonlinear terms and the coupling terms, respectively, without using the approximation argument(see Proposition 3.4 and Theorem 3.1 in \cite{MZ}). From \cite{MZ} and the present work, it is clear that the regularity of the solution plays an important role for the Malliavin differentiability of the solution. In other words, if the solution is not regular enough then the solution can not be
Malliavin differentiable. To be more precise, if the global well-posedness of is only available for the weak solution to a nonlinear partial differential equation, one can not achieve the global well-posedness of the equation when the initial condition is randomized, which is also true for the equations with random boundary conditions. This reveals a clear difference between the partial differential equations and stochastic partial differential equations (see Theorem \ref{prop:3.4}).

As shown in \eqref{eqn:1.1}-\eqref{eqn:1.3}, our model deals with the case that only the velocity field is perturbed by the noise. This is because one needs to make use of the particular geometric structure of the nematic liquid crystal equation, \textit{the basic balance law} (see \cite{LL} for reference), to obtain the energy estimates of velocity field as well as orientation field in certain regular spaces. If both the velocity equation and the orientation equation
were forced by noises, this would then destroy the basic balance law and one could not obtain a priori estimates, which are needed for proving the global well-posedness. However, we notice that in \cite{BHR1,BMP}, a noise was added in a special and smart manner to the orientation field equation but it did not bring essential difficulty while analysing the system. Moreover, we would be like to point out that, according to our work, the structure of the stochastic equations should be regular enough in order to obtain the global well-posedness of the stochastic model with random initial condition or random boundary condition. For instance, if the boundary condition of the orientation field is of Neumann type, then the global well-posedness is only true for the weak solution and one can not obtain the strong solution (see \cite{BHR1,BMP}). In this case, if either the initial condition or the boundary condition is randomised, one can not get the global well-posedness for the system,
see Theorem \ref{prop:3.4}.

\par
In the present paper, we are concerned with the initial and boundary problems for the nematic liquid crystal equations with multiplicative noise, here both the initial and the boundary conditions are randomised, which leads to the stochastic integrals defined via Skorohod integral, instead of It\^o integral. Thus, in order to show the global well-posedness result for the random initial and the random boundary problems (see Theorem \ref{thm:main}), we must establish the regularities of the solutions with respect to the initial data as well as the sample path. Specifically saying, we need to show the solutions $\v(t,R_{\nu},\omega),\d(t,R_d,\omega)$ are differentiable with respect to the random fields $R_{\nu},R_d$ and the sample path $\omega$. We would like to mention that  the regularity results established  in Theorem \ref{thm:3.2} and Theorem \ref{prop:3.4} are new and profound which do not exist  in previous work even for the deterministic case. As shown in our proof, the main difficulties lie in the process of dominating the highly nonlinear terms and the coupling terms. So in order to conquer that,  we make full usage of the geometric structure obtaining more delicate estimates: Propostion \ref{prop:2.5}, Proposition \ref{prop:2.6}, which are key \textit{a priori} estimates to establish the regularities of the present stochastic system with random initial and random boundary conditions.
Technically, we are not able to follow the standard arguments to obtain the estimates as in the proofs of Theorem \ref{thm:3.2} and Theorem \ref{prop:3.4}.
The feature of our derivations is as follows. By noticing that the boundary condition is independent of the time variable and the special geometric structure of the equation (see the equation (\ref{eqn:3.7})), we formulate an equivalent equation (see \eqref{eqn:3.9a}). Further, by estimating the equation \eqref{eqn:3.9a}, we then obtain the energy estimates for the equation (\ref{eqn:3.7}) in $\mathbb{H}^{1}$ space. But what we need is to get energy estimates
of the orientation field in $\mathbb{H}^{2}$ space. This creates a new difficulty. If we try to follow the standard argument to establish estimates of orientation field in $\mathbb{H}^{2}$ space, that is, naturally by taking inner product between $\Delta\hat\d(t,\d_0)(\b_0)$ and $\Delta\hat\d_t(t,\d_0)(\b_0)$ in $L^2(\mathbf{D})$ space, then we will encounter a non-treatable term $\langle \Delta \hat\d,  \Delta^{2} \hat\d \rangle.$ We even do not know its sign. To get rid of this difficulty, we recall that the orientation is independent of time on the boundary, a key observation which leads to an important integration by parts assuring the derivation of the energy estimates in $\mathbb{H}^{2}$ space. Precisely, noting that $ \langle \Delta\hat\d(t,\d_0)(\b_0), \partial_{t}  \Delta\hat\d(t,\d_0)(\b_0)   \rangle=\langle \nabla\Delta\hat\d(t,\d_0)(\b_0), \nabla \partial_{t}\d\rangle$, then $\langle \Delta \hat\d,  \Delta^{2} \hat\d \rangle$ is replaced with $ \langle \nabla\Delta \hat\d,  \nabla\Delta \hat\d \rangle\ge0$ which then opens the way to obtain the energy estimates of orientation field $\hat\d$ in $\mathbb{H}^{2}$ space, see Equation (3.14). These techniques, together with Propositions \ref{prop:2.5} and \ref{prop:2.6} play essentially  important roles in obtaining Theorem \ref{thm:3.2}, Theorem \ref{prop:3.4} and Proposition \ref{prop:3.5}, which ensure the establishment of our main result Theorem \ref{thm:main} and Theorem \ref{thm:main1}.

The rest of this paper is organized as follows: in Section 2, we define some functional spaces and give the abstract model expression for the stochastic model. The main result is also given in this section. In Section 3, we derive \textit{a priori} estimates and discuss Malliavin differentiability of the stochastic model with deterministic initial conditions and deterministic boundary conditions. In Section 4, we get back to the anticipating model and prove the stochastic nematic liquid crystals flows with random initial conditions and with random boundary conditions.

As usual, the constant $C$ may change from one line to another except that we give a special declaration, we denote by $C(a)$ a constant that depends on some parameter $a$.

%%%%%%%%%%%%%%%%%%%%%%
\section{Preliminaries and the main result}

We first set a space
$$ V=\{\v\in(C_0^\infty(\mathbf{D}))^2:\nabla\cdot\v=0\}. $$
Now we define spaces $ \vH, \mathbf V$, and $\mathbf H^m$  as the closure of $V$ in $(L^2(\mathbf{D}))^2,(H^1(\mathbf{D}))^2$  and $(H^m(\mathbf{D}))^2$, respectively. Let $|\cdot|_{2}$ and $\qv{\cdot,\cdot}$ be the norm and inner product in the space $\vH$, and  let $\norm{\cdot}_{1}$ and $\langle, \rangle_{\V}$ stand for the  norm and the inner product in the space $\V$, where $\langle, \rangle_{\V}$ is defined by
$$\qv{\v,\u}_\V:=\int_D\nabla \v\cdot \nabla \u dx,\ \ \text{ for }\  \mathbf{v}, \mathbf{u}\in \V.$$
Moreover, by  Poincar\'{e}'s inequality,  there exists a constant $c$ such that for any $\mathbf{v}\in \V$ we have $\|\mathbf{v}\|_{1}\leq c|\nabla \mathbf{v}|_{2}$.

Let $\mathbb{ H}^{m}=({H}^{m}(\mathbf{D}))^{3},m=0, 1,2,.... $ When $m=0,$ set $\mathbb{H}=\mathbb{H}^{0}= (L^{2}(\mathbf{D}))^{3} $ for simplicity. Then similarly, let $|\cdot|_{2}$ and $\qv{\cdot,\cdot}$ be the norm and inner product in the space $\dH$, and  let $\norm{\cdot}_{1}$ and $\langle, \rangle_{\dH^1}$ stand for the  norm and the inner product in the space $\dH^1$, where $\langle, \rangle_{\dH^1}$ is defined by
$$\qv{\d,\b}_{\dH^1}:= \int_D \d\cdot \b dx  +\int_D\nabla \d\cdot \nabla \b dx,\ \ \text{ for }\  \mathbf{d}, \mathbf{b}\in \dH^1.$$
Denote by $\V'$ the dual space of $\V$. And define the linear operator $A_{1}:\V \mapsto \V'$ as the following:
$$\qv{A_1\v,\u}=\qv{\v,\u}_\V, \text{ for }\v,\u\in\V.$$
Since the operator $A_{1}$ is positive self-adjoint with compact resolvent,
by the classical spectral theorem, $A_1$ admits an increasing sequence of eigenvalues $\alpha_j$ diverging to infinity with the corresponding eigenvectors $e_j$. Assume
\begin{equation}
\label{eqn:2.1}
\sum_{i=1}^{\infty} \lambda_{i}\alpha_{i}^{2}<\infty.
\end{equation}
Let $D(A_{1}):=\{\v\in \V, A_{1}\v\in \vH  \}$, since $A_1^{-1}$ is a self-adjoint compact operator as well, due to the classic spectral theory,  we can define the power $A_{1}^{s}$ for any $s\in \mathbb{R}.$ Moreover, $D(A_{1})'= D(A_{1}^{-1})$ is the dual space of $D(A_{1})$. And we have the compact embedding relationship
$$
D(A_{1})\subset \V \subset \vH\cong\vH' \subset \V'\subset D(A_{1})',
\text{ and }
\langle \cdot,  \cdot \rangle_{\mathbf{V}}=\langle A_{1}\cdot, \cdot    \rangle= \langle A_{1}^{\frac{1}{2}}\cdot, A_{1}^{\frac{1}{2}}\cdot\rangle.
$$
We define another operator $A_2:\dH^2\to \dH$ by $- \Delta$ satisfying
% $$\qv{A_2\d,\b}=\qv{\d,\b}_{\dH^1}=\int_D\nabla\d\cdot\nabla\b dx,\ \d,\b\in\dH^1.$$
 $D(A_2):=\{\d\in\dH^2;\d =R_d(x)\in \dH^2\ \mathrm{on}\ \partial \mathbf{D}\}$. Obviously, we have the compact embedding relationship
$$\dH^2\subset\dH^1\subset\dH\cong\dH'\subset(\dH^1)'\subset(\dH^2)'.$$
Define the trilinear form $b_1$ by
$$b_1(\u,\v,\w)=\sum_{i,j=1}^{2}\int_Du^i\partial_{x_i}v^jw^jdx, \text{ for }\u  \in \mathbf{H}\ \text{and } \v, \w \in \mathbf{V} \text{ and integral exists.}$$
If $\u,\v,\w\in\V$, then
$$|b_1(\u,\v,\w)|\leq c|\u|_2\norm{\v}_1\norm{\w}_1.$$
Now we define a bilinear form $B_1(\u,\v):=b_1(\u,\v,\cdot)$, then $B_1(\u,\v)\in\V'$ for $\u,\v\in\V$ and enjoys the following bound:
\begin{equation}
	\norm{B_1(\u,\v)}_{\V'}\leq c|\u|_2\norm{\v}_1.
\end{equation}
\begin{lemma}
	\label{lemma:2.1}
	The mapping $B_1:\V\times\V\to\V'$ is bilinear and continuous, and $b_1,B_1$ have the following properties:
$$	\begin{array}{lll}
b_1(\u,\v,\w)=-b_1(\u,\w,\v),&\ \qv{B_1(\u,\v),\w}=-\qv{B_1(\u,\w),\v},&\text{ for }\u,\v,\w\in\V.\\
b_1(\u,\v,\v)=0,&\ \qv{B_1(\u,\v),\v}=0,&\text{ for }\u,\v\in\V.
\end{array}
$$
Moreover, if $\u,\v,\w\in\V$, we have
\begin{equation}
\label{eqn:2.3}
	|b_1(\u,\v,\w)|=\qv{B_1(\u,\v),\w}\leq 2|\u|_2^{\frac{1}{2}}\norm{\u}_1^{\frac{1}{2}}\norm{\v}_1|\w|_2^{\frac{1}{2}}\norm{\w}_1^{\frac{1}{2}}.
\end{equation}
\end{lemma}
Define another trilinear form $b_2$ by
$$b_2(\v,\d,\b)=\sum_{j=1}^{3}\sum_{i=1}^{2}\int_Dv^i\partial_{x_i}d^jb^jdx\text{ for }\v\in\vH,\d\ \mathrm{and}\ \b \in \mathbb{H}^{1}.$$
Define another bilinear map $B_2$ on $\vH\times \mathbb{H}^{1}$ taking values in $\mathbb{ H}^{-1}$ such that $\qv{B_2(\v,\d),\b}:=b_2(\v,\d,\b)$.
\begin{lemma}
	\label{lemma:2.2}
	For $\v\in\V, \b\in\dH^1, \d\in\dH^2$, there exists a constant $c$ such that
	\begin{equation}
		|b_2(\v,\d,\b)|=|\qv{B_2(\v,\d),\b}|\leq c|\v|_2\norm{\d}_1\norm{\b}_1.
	\end{equation}
	Moreover, we have
	\begin{equation}
		|B_2(\v,\d)|_{2}\leq c\norm{\v}_1\norm{\d}_2,\hspace{1cm}\qv{B_2(\v,\d),\d}=0.
	\end{equation}
\end{lemma}
Now define the trilinear form $m$ by setting
$$m(\d,\b,\v)=\sum_{i,j=1}^{2}\sum_{k=1}^{3}\int_D\partial_{x_i}d^k\partial_{x_j}b^k\partial_{x_i}v^jdx.$$
There exists a bilinear operator $M$ defined on $\dH^2\times\dH^2$ taking values in $\V'$ such that $\qv{M(\d,\b),\v}:=m(\d,\b,\v).$
By interpolation inequality,  we can easily obtain
\begin{lemma}
	\label{lemma:2.3}
For any $\d,\b\in\dH^2,\v\in\V$, there exists a constant $c$ such that
$$|m(\d,\b,\v)|\leq c\norm{\d}_1^{\frac{1}{2}}\|\d\|_2^{\frac{1}{2}}\norm{\b}_1^{\frac{1}{2}}\|\b\|_2^{\frac{1}{2}}\norm{\v}_1.$$
Thus, for any $\d,\b\in\dH^2$,
\begin{equation}
	\norm{M(\d,\b)}_{\V'}\leq c\norm{\d}_1^{\frac{1}{2}}\|\d\|_2^{\frac{1}{2}}\norm{\b}_1^{\frac{1}{2}}|\Delta\b|_2^{\frac{1}{2}}.
\end{equation}
\end{lemma}
Now we arrive at the useful basic balance law and we include the proof here for reader's convenience.
\begin{lemma}
	\label{lemma:2.4}
	For $\u\in\V, \d\in\dH^2$, we have
	$$\qv{M(\d,\d),\u}=\qv{B_2(\u,\d),\Delta\d}.$$
\end{lemma}
\proof By integration by parts and the boundary conditions \eqref{BC}, we have
\begin{align*}
	\qv{M(\d,\d),\u}=&\qv{\nabla\cdot(\nabla\d\odot\nabla\d),\u}=\int_D\partial_{x_i}(\partial_{x_i}d^k\partial_{x_j}d^k)u^jdx\\
	=&-\int_D\partial_{x_i}d^k\partial_{x_j}d^k\partial_{x_i}u^jdx,
\end{align*}
and
\begin{align*}
	\qv{B_2(\u,\d),\Delta\d}=&\qv{\u\cdot\nabla\d,\Delta\d}=\int_Du^i\partial_{x_i}d^k\partial_{x_jx_j}d^kdx\\
	=&-\int_D\partial_{x_j}u^i\partial_{x_i}d^k\partial_{x_j}d^kdx-\int_Du^i\partial_{x_ix_j}d^k\partial_{x_j}d^kdx\\
	=&-\int_D\partial_{x_j}u^i\partial_{x_i}d^k\partial_{x_j}d^kdx=\qv{M(\d,\d),\u}.
\end{align*}
\qed

In the following, we will state two important results that are used several times in the rest of the paper.
\begin{prop}
	\label{prop:2.5}
	For $\d, \b \in \mathbb{H}^{2}$ and $\u \in \V,$ we have
	\begin{equation}
		\qv{M(\d,\b),\u}+\qv{M(\b,\d),\u}=\qv{B_2(\u,\d),\Delta\b}+\qv{B_2(\u,\b),\Delta\d}.
	\end{equation}
	
\end{prop}
\proof By the bilinear property of the operator $M$, and the basic balance law in Lemma \ref{lemma:2.4},
\begin{align*}
	&\qv{M(\d,\b),\u}+\qv{M(\b,\d),\u}\\
	=&\qv{M(\d,\d),\u}-\qv{M(\d-\b,\d-\b),\u}+\qv{M(\b,\b),\u}\\
	=&\qv{B_2(\u,\d),\Delta\d}-\qv{B_2(\u,\d-\b),\Delta(\d-\b)}+\qv{B_2(\u,\b),\Delta\b}\\
	=&\qv{B_2(\u,\d),\Delta\b}+\qv{B_2(\u,\b),\Delta\d}.
\end{align*}
\qed
\begin{prop}
	\label{prop:2.6}
For $\d, \b \in \mathbb{H}^{3}$ and $\u \in \V,$ and continuous functions $\alpha(s),\beta(s), s\in [0, t] $, we get
	\begin{align}
	&\int_0^t\alpha(s)\qv{M(\d,\b),\u}ds+\int_0^t\alpha(s)\qv{M(\b,\d),\u}ds-\int_0^t\beta(s)\qv{B_2(\u,\b),\Delta\d}ds\notag\\
	\leq&2(|\alpha|_\infty+|\beta|_\infty)\int_0^t\norm{\d}_1^{1/2}\norm{\d}_2^{1/2}\norm{\b}_1^{1/2}\norm{\b}_2^{1/2}\norm{\u}_1ds+|\beta|_\infty\int_0^t|\u|_2\norm{\d}_1\norm{\b}_3ds.
	\end{align}
	Or
	\begin{align}
	&\int_0^t\alpha(s)\qv{M(\d,\b),\u}ds+\int_0^t\alpha(s)\qv{M(\b,\d),\u}ds-\int_0^t\beta(s)\qv{B_2(\u,\b),\Delta\d}ds\notag\\
\leq&|\alpha|_\infty\int_0^t|\u|_2\norm{\d}_1\norm{\b}_3ds+(|\alpha|_\infty+|\beta|_\infty)\int_0^t|\u|_2\norm{\b}_1\norm{\d}_3ds.
	\end{align}
\end{prop}
\proof
With different time function coeffecients, we apply the  identity in Proposition \ref{prop:2.5}, together with Lemma \ref{lemma:2.2}, \ref{lemma:2.3},
	\begin{align*}
&\int_0^t\alpha(s)\qv{M(\d,\b),\u}ds+\int_0^t\alpha(s)\qv{M(\b,\d),\u}ds-\int_0^t\beta(s)\qv{B_2(\u,\b),\Delta\d}ds\\
=&\int_0^t(\alpha(s)-\beta(s))\qv{M(\d,\b),\u}ds+\int_0^t(\alpha(s)-\beta(s))\qv{M(\b,\d),\u}ds\\
&+\int_0^t\beta(s)\qv{M(\d,\b),\u}ds+\int_0^t\beta(s)\qv{M(\b,\d),\u}ds-\int_0^t\beta(s)\qv{B_2(\u,\b),\Delta\d}ds\\
=&\int_0^t(\alpha(s)-\beta(s))\qv{M(\d,\b),\u}ds+\int_0^t(\alpha(s)-\beta(s))\qv{M(\b,\d),\u}ds+\int_0^t\beta(s)\qv{B_2(\u,\d),\Delta\b}ds\\
\leq&2(|\alpha|_\infty+|\beta|_\infty)\int_0^t\norm{\d}_1^{1/2}\norm{\d}_2^{1/2}\norm{\b}_1^{1/2}\norm{\b}_2^{1/2}\norm{\u}_1ds+|\beta|_\infty\int_0^t|\u|_2\norm{\d}_1\norm{\b}_3ds.
\end{align*}
Or, direclty applying Proposition \ref{prop:2.5} and Lemma \ref{lemma:2.2},
\begin{align*}
&\int_0^t\alpha(s)\qv{M(\d,\b),\u}ds+\int_0^t\alpha(s)\qv{M(\b,\d),\u}ds-\int_0^t\beta(s)\qv{B_2(\u,\b),\Delta\d}ds\\
=&\int_0^t\alpha(s)\qv{B_2(\u,\d),\Delta\b}ds+\int_0^t\alpha(s)\qv{B_2(\u,\b),\Delta\d}ds-\int_0^t\beta(s)\qv{B_2(\u,\b),\Delta\d}ds\\
=&\int_0^t\alpha(s)\qv{B_2(\u,\d),\Delta\b}ds+\int_0^t(\alpha(s)-\beta(s))\qv{B_2(\u,\b),\Delta\d}ds\\
\leq&|\alpha|_\infty\int_0^t|\u|_2\norm{\d}_1\norm{\b}_3ds+(|\alpha|_\infty+|\beta|_\infty)\int_0^t|\u|_2\norm{\b}_1\norm{\d}_3ds.
\end{align*}
\qed
\begin{remark}
	Proposition \ref{prop:2.5} and Proposition \ref{prop:2.6} are very important to bound the nonlinear terms when we try to obtain the regularities of the solutions with respect to initial data and sample path, please see  Theorem \ref{thm:3.2}, Theorem \ref{prop:3.4} and Proposition \ref{prop:3.5}. In fact, these kinds of regularities are profound results which do not exist in previous work even for the deterministic case. In the proving process of these result (see Theorem \ref{thm:3.2}, Theorem \ref{prop:3.4} and Proposition \ref{prop:3.5}), the difficulties lie in bounding the highly nonlinear term which obliges us to take full advantage of the delicate geometric structure of the stochastic nematic liquid crystals equations. Hence, Proposition \ref{prop:2.5} and Proposition \ref{prop:2.6} are just two of the key observations to study the regularities of this stochastic model with random initial conditions and random boundary conditions.
	
	%As it is shown in the equivalent model \eqref{eqn:3.3}, when dealing with \textit{a priori} estimates, higher order regularity arises from the estimate of the bilinear operator $B_2$ involved with both velocity and orientation fields. Hence, it requires to make use of \textit{the basic balance law} to reduce the higher order regularity. However, one can note that there are different time function coefficients for the operators $M, B_2$, so we need to apply Proposition \ref{prop:2.5} to have proper regularity results.
\end{remark}

Finally, $f(\mathbf{d})$ and $F(\mathbf{d})$  are given by
\begin{eqnarray}
\label{eqn:2.4}
f(\mathbf{d})=\frac{1}{\eta^{2}}(|\mathbf{d}|^{2}-1)\mathbf{d}\ \ \mathrm{and}\ \ F(\mathbf{d})= \frac{1}{4\eta^{2}}(|\mathbf{d}|^{2}-1)^{2}.
\end{eqnarray}
We define a function $ \wt{f}:[0, \infty)\rightarrow \mathbb{R}$ by
\begin{eqnarray}
\label{eqn:2.5}
\wt{f}(x)=\frac{1}{\eta^2}(x-1),  \ \ \   x\in \mathbb{R}_{+},
\end{eqnarray}
then $f(\d)=\wt f(|\d|^2)\d$ and denote by ${F}: \mathbb{R}^{3} \rightarrow \mathbb{R}$ the Fr\'echet differentiable map such that for any $\mathbf{d}\in \mathbb{R}^{3}$ and $\xi\in \mathbb{R}^{3}$
\begin{eqnarray}
\label{eqn:2.6}
{F}'(\mathbf{d})[\xi]=f(\mathbf{d})\cdot \xi.
\end{eqnarray}
Set $\wt{F}$ to be an antiderivative of $\wt{f}$ such that $\wt {F}(0)=0.$ Then
$$\wt{F}(x)= \frac{1}{2\eta^2}(x^2-2x),  \ \ \   x\in \mathbb{R}_{+}.$$
\begin{defn}\label{def-2.1}
We say a continuous $\mathbf{H}\times \mathbb{H}^{1}$ valued random field $(\mathbf{v}(.,t), \mathbf{B}(.,t) )_{t\in [0,T]}$ defined on $(\Omega, \mathcal{F}, \mathbb{P})$ is a weak solution to problem (\ref{eqn:1.1})-(\ref{eqn:1.3}) with initial and boundary conditions \eqref{IC} and \eqref{BC} if for $(\mathbf{v}_{0}, \mathbf{d}_{0})\in \mathbf{H}\times \mathbb{H}^{1}$ the following conditions hold:
\begin{align*}
& \mathbf{v}\in C([0,T];\mathbf{H})\cap L^{2}([0,T]; \mathbf{V}),\\
& \mathbf{d}\in C([0,T];\mathbb{H}^{1})\cap L^{2}([0,T]; \mathbb{H}^{2}),
\end{align*}
and the integral relation
\begin{align*}
\langle \mathbf{v}(t), v   \rangle +&\int_{0}^{t}\langle  A_{1}\mathbf{v}(s),  v  \rangle ds+\int_{0}^{t}\langle \mathbf{v}(s)\cdot\nabla \mathbf{v}(s),  v  \rangle ds\\
&+ \int_{0}^{t}\langle\nabla \cdot (\nabla \mathbf{d}(s)\odot \nabla \mathbf{d}(s)), v\rangle ds=\langle \mathbf{v}_{0}, v \rangle +\int_{0}^{t} \langle \sum_{k=1}^{\infty}\sigma_{k}\v\circ dW_{k}(s), v\rangle+\langle W_{0}(t), v   \rangle ,\\
\langle \mathbf{d}(t),  d  \rangle+&\int_{0}^{t}\langle  A_{2}\mathbf{d}(s),  d  \rangle ds+\int_{0}^{t}\langle \mathbf{v}(s)\cdot\nabla \mathbf{d}(s), d\rangle ds\\
&=\langle \mathbf{d}_{0}, d   \rangle-\int_{0}^{t}\langle f(\mathbf{d}(s)), d\rangle ds,
\end{align*}
hold $a.s.$ for all $t\in [0,T]$ and $(v,d)\in \mathbf{V}\times \mathbb{H}.$
\end{defn}

\begin{defn}\label{def-2.2}
We say a continuous $\mathbf{V}\times \mathbb{H}^{2}$ valued random field $(\mathbf{v}(.,t), \mathbf{B}(.,t) )_{t\in [0,T]}$ defined on $(\Omega, \mathcal{F}, \mathbb{P})$ is a strong solution to problem (\ref{eqn:1.1})-(\ref{eqn:1.3}) with initial and boundary conditions \eqref{IC} and \eqref{BC} if for $(\mathbf{v}_{0}, \mathbf{d}_{0})\in \mathbf{V}\times \mathbb{H}^{2}$ the following conditions hold:
\begin{align*}
& \mathbf{v}\in C([0,T];\mathbf{V})\cap L^{2}([0,T]; \mathbf{H}^{2}),\\
& \mathbf{d}\in C([0,T];\mathbb{H}^{2})\cap L^{2}([0,T]; \mathbb{H}^{3}),
\end{align*}
and the integral relation
\begin{align*}
 \mathbf{v}(t) +&\int_{0}^{t}  A_{1}\mathbf{v}(s) ds+\int_{0}^{t}\mathbf{v}(s)\cdot\nabla \mathbf{v}(s)ds\\
&+ \int_{0}^{t}\nabla \cdot (\nabla \mathbf{d}(s)\odot \nabla \mathbf{d}(s))ds=\mathbf{v}_{0}+\int_{0}^{t} \sum_{k=1}^{\infty}\sigma_{k}\v\circ dW_{k}(s)+ W_{0}(t) ,\\
\mathbf{d}(t)+&\int_{0}^{t}  A_{2}\mathbf{d}(s) ds+\int_{0}^{t} \mathbf{v}(s)\cdot\nabla \mathbf{d}(s) ds=\mathbf{d}_{0}-\int_{0}^{t} f(\mathbf{d}(s)) ds,
\end{align*}
hold $a.s.$ for all $t\in [0,T]$.
\end{defn}

Now the equations \eqref{eqn:1.1}-\eqref{eqn:1.3} can be written as
\begin{align}
\label{eqn:2.13}
	&d\v(t)+[A_1\v(t)+B_1(\v(t))+M(\d(t))]dt=\sum_{k=1}^{\infty}\sigma_k\v(t)\circ W_k(t)+\sigma_0dW_0(t),\\
	\label{eqn:2.14}
	&d\d(t)+[A_2\d(t)+B_2(\v(t),\d(t))+f(\d(t))]dt=0,
\end{align}
with the initial conditions $\v(0)=R_{\nu},\d(0)=R_d$.

Throughout the paper, we denote by $\D$ the Malliavin differentiation of random variables on the Wiener space $(\Omega,\F,\P)$. And we denote by $\D^{1,2}(\vH) $ the Malliavin Sobolev space of all $\F$-measurable and Malliavin differentiable random variables $\Omega\to\vH $ with Malliavin derivatives owing second order moments. Correspondingly, $\D_{\text{loc}}^{1,2}(\vH)$ represents the space of random variables $\xi:\Omega\to\vH$ that are locally in $\D^{1,2}(\vH)$.

We end up this section by presenting our main theorems, which give the existence and uniqueness of solutions to the stochastic model \eqref{eqn:1.1}-\eqref{eqn:1.3}, or \eqref{eqn:2.13}-\eqref{eqn:2.14}, with  random boundary conditions \eqref{BC} and  random initial conditions \eqref{IC}. The proof is given in Section 4.
\begin{thm}
	\label{thm:main}
	Assume the initial random field $R_{\nu}\in\D^{1,2}_{\text{loc}}(\vH)\cap\V,R_d\in\D^{1,2}_{\text{loc}}(\dH^1)\cap\dH^2$, then the stochastic nematic liquid crystal flows (\ref{eqn:1.1})-(\ref{eqn:1.3}) together with equation (\ref{IC}) and equation (\ref{BC}) have a unique strong solution $(\v(t,R_{\nu}),\d(t,R_d))$ for $t\in[0,T]$. Moreover, $\v(t,R_{\nu})\in\D_{\text{loc}}^{1,2}(\vH), \d(t,R_d)\in\D_{\text{loc}}^{1,2}(\dH^1)$.
\end{thm}

In fact, following the argument of Theorem \ref{thm:main} with minor modification, we can extend our result to the random boundary conditions with stochastic force. That is our second main result which we state as follows without proof.
\begin{thm}
	\label{thm:main1}
In addition to the assumption the initial random field $R_{\nu}\in\D^{1,2}_{\text{loc}}(\vH)\cap\V,R_d\in\D^{1,2}_{\text{loc}}(\dH^1)\cap\dH^2$, let
\begin{align}
\label{BC'}&\v(x,t)=0,\ \d(x,t)=R_{d}(x)+ N(t),\ \text{ for }(x,t)\in\partial \mathbf{D}\times\R^+,\tag{BC\,$'$}
\end{align}
 with $N(t)$ being a compound Poisson process taking values in $\mathbb{H}^2$ space, then the stochastic nematic liquid crystals flows (\ref{eqn:1.1})-(\ref{eqn:1.3}) together with equation (\ref{IC}) and equation (\ref{BC'}) have a unique strong solution $(\v(t,R_{\nu}),\d(t,R_d))$ for $t\in[0,T]$. Moreover, $\v(t,R_{\nu})\in\D_{\text{loc}}^{1,2}(\vH), \d(t,R_d)\in\D_{\text{loc}}^{1,2}(\dH^1)$.
\end{thm}

We would like to point out that the external force on the boundary should be specified, based on the following two reasons. First, the sample paths
of the random force should be piece wise differentiable with respect to time $t$. That excludes the boundary condition involving a Brownian motion as it is well known that the Brownian paths are nowhere differentiable. However, the sample paths of a compound Poisson process, as random step functions, are indeed piece wise differentiable with respect to time $t$. Second, the time non-homogenous boundary condition breaks the integration by parts formula (e.g., see (\ref{eqn:3.14a})) which is vital for the
derivation of the energy estimates for the orientation fields in $\mathbb{H}^2$ space, see Theorem \ref{thm:3.2} and Theorem \ref{prop:3.4}.

 %%%%%%%%%%%% %%%%%%%%%%%
\section{A priori estimates}
%%%%%%%%%%%%%%%%%%%%%%%
Consider the stochastic model with a deterministic initial condition $(\v_0, \d_{0})\in\V\times \mathbb{H}^{2}$,
\begin{align}
\label{eqn:3.1}
	&d\v(t,\v_0)+[A_1\v(t,\v_0)+B_1(\v(t,\v_0))+M(\d(t,\d_0))]dt=\sum_{k=1}^{\infty}\sigma_k\v(t)\circ W_k(t)+\sigma_0dW_0(t),\notag\\
	&d\d(t,\d_0)+[A_2\d(t,\d_0)+B_2(\v(t,\v_0),\d(t,\d_0))+f(\d(t,\d_0))]dt=0,\notag\\
	&\v(0,\v_0)=\v_0\in\V,\ \d(0,\d_0)=\d_0\in\dH^2.
\end{align}
The global well-posedness for the strong solution  of \eqref{eqn:3.1} has been studied in \cite{BHR2} and \cite{GHZ}, and it is known that under the condition \eqref{eqn:2.1}, for any $T>0$, $\v(\cdot,\v_0)\in C([0,T];\V)\cap L^2([0,T];\mathbf H^2),\ \d(\cdot,\d_0)\in C([0,T];\mathbb H^2)\cap L^2([0,T];\mathbb H^3)$.

Define $$Q(t):=\exp\{\sum_{k=1}^{\infty}\sigma_kW_k(t)\},$$ then $Q(0)=1$, by Novikov condition and Doob's maximal inequality, we have $\Ex\sup_{0\leq t\leq T}|Q(t)|<\infty$ for arbitrary $T>0$. For simplicity of notations, we use $|Q|_\infty$ represent $\sup_{0\leq s\leq t}|Q(s)|$.

 Let $Z(t)$ be the unique solution of the stochastic equation:
\begin{align}
	&dZ(t)=-A_1Z(t)dt+\sigma_0Q(t)^{-1}dW_0(t);\notag\\
	&Z(0)=0;\notag\\
	&Z(t,x)=0,\ x\in\partial \mathbf{D}, t\geq0.
\end{align}
Now define $\u(t,\v_0):=\v(t,\v_0)Q(t)^{-1}-Z(t)$, $t\geq0$, then by It\^o's formula, $\u,\d$ satisfy the following equations:
\begin{align}
\label{eqn:3.3}
&d\u(t)+[A_1\u(t)+Q(t)B_1\big(\u(t)+Z(t)\big)+Q(t)^{-1}M(\d(t))]dt=0,\notag\\
&d\d(t)+\Big[A_2\d(t)+Q(t)B_2\big(\u(t)+Z(t),\d(t)\big)+f(\d(t))\Big]dt=0,\notag\\
&\u(0)=\v_0,\ \d(0)=\d_0.
\end{align}
Using the estimates in  \cite{BHR1} or \cite{GHZ}, we obtain the following proposition,
\begin{prop}
	\label{prop:3.1}
	For $\v_0\in\V,\d_0\in\dH^2$ and $\omega\in\Omega$. Denote by $(\u(t,\v_0,\omega),\d(t,\d_0,\omega))$ the unique solution to \eqref{eqn:3.3} on $[0,T]$. Then the following estimates hold:
	\begin{align}
	&\sup_{0\leq t\leq T}[|\u(t,\v_0,\omega)|_2^2+\norm{\d(t,\d_0,\omega)}_1^2]+\int_0^T\norm{\u(t,\v_0,\omega)}_1^2dt+\int_0^T\norm{\d(t,\d_0,\omega)}_2^2dt\notag\\
	\leq&c(|\v_0|_2,\norm{\d_0}_1,|Q|_\infty,\sup_{0\leq t\leq T}\norm{Z}_2,T),
	\end{align}
	and
	\begin{align}
		&\sup_{0\leq t\leq T}[\norm{\u(t,\v_0,\omega)}_1^2+\norm{\d(t,\d_0,\omega)}_2^2]+\int_0^T\norm{\u(t,\v_0,\omega)}_2^2dt+\int_0^T\norm{\d(t,\d_0,\omega)}_3^2dt\notag\\
		\leq&c(\norm{\v_0}_1,\norm{\d_0}_2,|Q|_\infty,\sup_{0\leq t\leq T}\norm{Z}_2,\int_0^T\norm{Z}_3^2dt,T).
	\end{align}
\end{prop}
The following theorem states the regularity of the solutions  to \eqref{eqn:3.3}, which is differentiable with respect to initial data.
\begin{thm}
	\label{thm:3.2}
	For $(\v_0,\d_0)\in\V\times\dH^2$, $(\v(t,\v_0,\omega),\d(t,\d_0,\omega))\in\V\times\dH^2$, and the solution map $(\v_0,\d_0)\mapsto(\v(t,\v_0,\omega),\d(t,\d_0,\omega))$ is $C^{1,1}$ for all $\omega\in\Omega, t\geq0$, and has bounded Fr\'echet derivatives  on bounded sets in $\V\times\dH^2$.
	
	Moreover, the Fr\'echet derivative $t\mapsto (D\v(t,\v_0,\omega),D\d(t,\d_0,\omega))\in L(\V\times\dH^2)$ is continuous in $t$, and the Fr\'echet derivative is compact for any $t>0,\omega\in\Omega$, where $ L(\V\times\dH^2)$ represents the space of bounded linear operators from $\V\times\dH^2$ to $\V\times\dH^2$.
\end{thm}
\begin{remark}
Since the stochastic equations (\ref{eqn:1.1})-(\ref{eqn:1.3}) is a coupled system of the velocity field equation (\ref{eqn:1.1})and and the orientation field equation (\ref{eqn:1.3}), then when we consider the differentiability of the solutions to (\ref{eqn:1.1})-(\ref{eqn:1.3}) with respect to the initial data, we need to calculate the derivative of equation (\ref{eqn:1.1}) and (\ref{eqn:1.3}) with respect to the initial data $(\v_{0}, \d_{0})$ at the same time. In other words, if we only consider the derivative of the equation of $\u$ with respect to $\v_{0}$ or the derivative of the equation of $\d$ with respect to $\d_{0}$ is not true. See, for example, (\ref{eqn:3.6}) and (\ref{eqn:3.7}). Therefore, the coupled system (\ref{eqn:1.1})-(\ref{eqn:1.3}) is more difficult than 2D stochastic Navier-Stokes equations and other stochastic hydrodynamic systems. We should carefully deal with this kind of stochastic coupled system.
\end{remark}
\proof Let $(\u(t,\v_0),\d(t,\d_0))_{t\in [0,T]}$ represent the unique strong solution to \eqref{eqn:3.3}, see \cite{GHZ}, where $(\u(t,\v_0),\d(t,\d_0))$ is shown to be Lipschitz continuous with respect to $(\v_0,\d_0)$ in $\V\times\dH^2$.

Given initial conditions $\u_0\in\V,\b_0\in\dH^2$, we consider the following random equations with boundary conditions (\ref{BC}):
\begin{align}
\label{eqn:3.6}
	\hat\u(t,\v_0)(\u_0)=&\u_0-\int_0^tA_1\hat\u(s,\v_0)(\u_0)ds-\int_0^tQ(s)B_1(\hat\u(s,\v_0)(\u_0),\u(s,\v_0)+Z(s))ds\notag\\
	&-\int_0^tQ(s)B_1(\u(s,\v_0)+Z(s),\hat\u(s,\v_0)(\u_0))ds-\int_0^tQ(s)^{-1}M(\hat\d(s,\d_0)(\b_0),\d(s,\d_0))ds\notag\\
	&-\int_0^tQ(s)^{-1}M(\d(s,\d_0),\hat\d(s,\d_0)(\b_0))ds;\\
	\label{eqn:3.7}
	\hat\d(t,\d_0)(\b_0)=&\b_0-\int_0^tA_2\hat\d(s,\d_0)(\b_0)ds-\int_0^tQ(s)B_2(\hat\u(s,\v_0)(\u_0),\d(s,\d_0))ds\notag\\
	&-\int_0^tQ(s)B_2(\u(s,\v_0)+Z(s),\hat\d(s,\d_0)(\b_0))ds-\int_0^tf'(\d(s,\d_0))\hat\d(s,\d_0)(\b_0)ds,
\end{align}
Obviously, the equations (\ref{eqn:3.6})-(\ref{eqn:3.7}) is linear, the global well-posedness of the strong solutions is easy to show. We omit it here. One can see \cite{BHR1}, \cite{Z} and other references.

Multiplying \eqref{eqn:3.6} with $\hat\u(t,\v_0)(\u_0)$, then integrating over $\mathbf{D}$ yields that
\begin{align}
\label{eqn:3.8}
	|\hat\u(t,\v_0)(\u_0)|_2^2=&|\u_0|_2^2-2\int_0^t|\nabla\hat\u(s,\v_0)(\u_0)|_2^2ds\notag\\
	&-2\int_0^tQ(s)\qv{B_1(\hat\u(s,\v_0)(\u_0),\u(s,\v_0)+Z(s)),\hat\u(s,\v_0)(\u_0)}ds\notag\\
	&-2\int_0^tQ(s)\qv{B_1(\u(s,\v_0)+Z(s),\hat\u(s,\v_0)(\u_0)),\hat\u(s,\v_0)(\u_0)}ds\notag\\
	&-2\int_0^tQ^{-1}(s)\qv{M(\hat\d(s,\d_0)(\b_0),\d(s,\d_0)),\hat\u(s,\v_0)(\u_0)}ds\notag\\
	&-2\int_0^tQ^{-1}(s)\qv{M(\d(s,\d_0),\hat\d(s,\d_0)(\b_0)),\hat\u(s,\v_0)(\u_0)}ds\notag\\
	=:&i_1+\cdots+i_6.
\end{align}
First by H\"older's inequality and Young's inequality,
\begin{align*}
	i_3\leq&\sup_{0\leq t\leq T}|Q(t)|\int_0^t|\nabla(\u(s,\v_0)+Z(s))|_\infty|\hat\u(s,\v_0)(\u_0)|_2^2ds\\
	\leq& c|Q|_\infty\int_0^t\norm{\u(s,\v_0)+Z(s)}_2|\hat\u(s,\v_0)(\u_0)|_2^2ds.
\end{align*}
By Lemma \ref{lemma:2.1}, $i_4=0$. For $i_5,i_6$, we get
\begin{align*}
i_5+i_6\leq& \sup_{0\leq t\leq T}|Q^{-1}(t)|\int_0^t|\nabla\d(s,\d_0)|_\infty|\nabla\hat\d(s,\d_0)(\b_0)|_2|\nabla\hat\u(s,\v_0)(\u_0)|_2ds\\
\leq&\eps\int_0^t|\nabla\hat\u(s,\v_0)(\u_0)|_2^2ds+c|Q^{-1}|_\infty^2\int_0^t\norm{\d(s,\d_0)}_2^2|\nabla\hat\d(s,\d_0)(\b_0)|_2^2ds.
\end{align*}
%For $i_5,i_6$, we apply Proposition \ref{prop:2.5}, H\"older's inequality and Young's inequality
%\begin{align*}
%	i_5+i_6=&-2\int_0^tQ^{-1}(t)\qv{B_2(\hat\u(s,\v_0)(\u_0),\hat\d(s,\d_0)(\b_0)),\Delta\d(s,\d_0)}ds\\
%	&-2\int_0^tQ^{-1}(t)\qv{B_2(\hat\u(s,\v_0)(\u_0),\d(s,\d_0)),\Delta\hat\d(s,\d_0)(\b_0)}ds\\
%	\leq&c|Q^{-1}|_\infty\int_0^t|\Delta\d(s,\d_0)|_\infty|\hat\u(s,\v_0)(\u_0)|_2|\nabla\hat\d(s,\d_0)(\b_0)|_2ds\\
%	&+c|Q^{-1}|_\infty\int_0^t|\nabla\d(s,\d_0)|_\infty|\hat\u(s,\v_0)(\u_0)|_2|\Delta\hat\d(s,\d_0)(\b_0)|_2ds
%\end{align*}
Hence, from \eqref{eqn:3.8}, we derive that
\begin{align}
\label{eqn:3.9}
	&|\hat\u(t,\v_0)(\u_0)|_2^2-|\u_0|_2^2+\int_0^t|\nabla\hat\u(s,\v_0)(\u_0)|_2^2ds\notag\\
	\leq&c|Q|_\infty\int_0^t\norm{\u(s,\v_0)+Z(s)}_2|\hat\u(s,\v_0)(\u_0)|_2^2ds+c|Q^{-1}|_\infty^2\int_0^t\norm{\d(s,\d_0)}_2^2|\nabla\hat\d(s,\d_0)(\b_0)|_2^2ds.
\end{align}
Under the boundary condition (BC), we even can not obtain the $L^2(\mathbf{D})$ estimates for $\hat\d$ by taking inner product of (\ref{eqn:3.7}) with $\hat\d$ in the corresponding space. As the integration by parts $ \langle A_{1}\hat\d,    \hat\d \rangle= | \nabla \hat\d|_{2}^{2}$ is not true. We can not also follow the argument as in $\cite{LL}$ to obtain the estimates needed. Hence, to obtain the energy estimates of $\hat\d$, we need to construct an equivalent equation in the following
\begin{align}\label{eqn:3.9a}
\hat\d(t,\d_0)(\b_0)&=\b_0-\int_0^tA_2\hat\d(s,\d_0)(\b_0)ds-\int_0^t\hat\d(s,\d_0)(\b_0)ds\notag\\
& -\int_0^tQ(s)B_2(\hat\u(s,\v_0)(\u_0),\d(s,\d_0))ds+ \int_0^t\hat\d(s,\d_0)(\b_0)ds\notag\\
	&-\int_0^tQ(s)B_2(\u(s,\v_0)+Z(s),\hat\d(s,\d_0)(\b_0))ds\notag\\
&-\int_0^tf'(\d(s,\d_0))\hat\d(s,\d_0)(\b_0)ds.
\end{align}
Then multiplying \eqref{eqn:3.9a} with $-\Delta\hat\d(t,\d_0)(\b_0)+\hat\d(t,\d_0)(\b_0)$, and integrating over $\mathbf{D}$ give that
\begin{align}
\label{eqn:3.10}
	&|\nabla\hat\d(t,\d_0)(\b_0)|_2^2+|\hat\d(t,\d_0)(\b_0)|_2^2\notag\\
	=&\norm{\b_0}_1^2-2\int_0^t|-\Delta\hat\d(s,\d_0)(\b_0)+\hat\d(s,\d_0)(\b_0)|_2^2ds\notag\\
	&-2\int_0^tQ(s)\qv{B_2(\hat\u(s,\v_0)(\u_0),\d(s,\d_0)),-\Delta\hat\d(t,\d_0)(\b_0)+\hat\d(t,\d_0)(\b_0)}ds\notag\\
	&-2\int_0^tQ(s)\qv{B_2(\u(s,\v_0)+Z(s),\hat\d(s,\d_0)(\b_0)),-\Delta\hat\d(t,\d_0)(\b_0)+\hat\d(t,\d_0)(\b_0)}ds\notag\\
	&-2\int_0^t\qv{(f'(\d(s,\d_0))-1)\hat\d(s,\d_0)(\b_0),-\Delta\hat\d(t,\d_0)(\b_0)+\hat\d(t,\d_0)(\b_0)}ds.
\end{align}
With similar discussion as in the estimates of $\hat\u$, we arrive at
\begin{align}
\label{eqn:3.11}
	&|\nabla\hat\d(t,\d_0)(\b_0)|_2^2+|\hat\d(t,\d_0)(\b_0)|_2^2-\norm{\b_0}_1^2+\int_0^t|-\Delta\hat\d(s,\d_0)(\b_0)+\hat\d(s,\d_0)(\b_0)|_2^2ds\notag\\
	\leq&c|Q|_\infty^2\int_0^t\norm{\d(s,\d_0)}_2^2|\hat\u(s,\v_0)(\u_0)|_2^2ds+c|Q|_\infty^2\int_0^t\norm{\u(s,\v_0)+Z(s)}_1^2|\hat\d(t,\d_0)(\b_0)|_2^2ds\notag\\
	&+c\int_0^t(\norm{\d(s,\d_0)}_1^2+1)^2|\hat\d(s,\d_0)(\b_0)|_2^2.
\end{align}
Combining the above estimates \eqref{eqn:3.9}, \eqref{eqn:3.11}, applying Gronwall's inequality, we get
\begin{align}\label{eqn:3.12}
	&\sup_{0\leq t\leq T}\left[|\hat\u(t,\v_0)(\u_0)|_2^2+|\nabla\hat\d(t,\d_0)(\b_0)|_2^2+|\hat\d(t,\d_0)(\b_0)|_2^2\right]\notag\\
	&+\int_0^T|\nabla\hat\u(t,\v_0)(\u_0)|_2^2dt+\int_0^T|-\Delta\hat\d(t,\d_0)(\b_0)+\hat\d(t,\d_0)(\b_0)|_2^2dt\notag\\
	\leq&(|\u_0|_2^2+\norm{\b_0}_1^2)\exp c\bigg\{|Q|_\infty\int_0^T\norm{\u(t,\v_0)+Z(t)}_2dt+|Q^{-1}|_\infty^2\int_0^T\norm{\d(t,\d_0)}_2^2dt\notag\\
	&\hspace{3.5cm}+|Q|_\infty^2\int_0^T(\norm{\d(t,\d_0)}_2^2+\norm{\u(t,\v_0)+Z(t)}_1^2)dt\notag\\
	&\hspace{3.5cm}+\int_0^T(\norm{\d(t,\d_0)}_1^2+1)^2dt\bigg\}<\infty.
\end{align}
Now we try to follow the standard argument to establish estimates of orientation field in $\mathbb{H}^{2}$ space.  By taking inner product between $\Delta\hat\d(t,\d_0)(\b_0)$ and $\Delta\hat\d_t(t,\d_0)(\b_0)$. Then we will encounter a non treatable term $\langle \Delta \hat\d,  \Delta^{2} \hat\d \rangle.$ We even do not know its sign. To get rid of this difficulty, we recall that the orientation is independent of time on the boundary, a key observation which leads to an important integration by parts assuring the derivation of the energy estimates in $\mathbb{H}^{2}$ space. Precisely, noting that $ \langle \Delta\hat\d(t,\d_0)(\b_0), \partial_{t}  \Delta\hat\d(t,\d_0)(\b_0)   \rangle=\langle \nabla\Delta\hat\d(t,\d_0)(\b_0), \nabla \partial_{t}\d\rangle$, then $\langle \Delta \hat\d,  \Delta^{2} \hat\d \rangle$ is replaced with $ \langle \nabla\Delta \hat\d,  \nabla\Delta \hat\d \rangle\ge0$ which then helps us to obtain the energy estimates of orientation $\hat\d$ in $\mathbb{H}^{2}$ space. Namely,
\begin{align}\label{eqn:3.14a}
\partial_{t}|\Delta\hat\d(t,\d_0)(\b_0)|_2^2=&2\langle \Delta\hat\d(t,\d_0)(\b_0), \partial_{t}  \Delta\hat\d(t,\d_0)(\b_0)   \rangle\notag\\
=&-2\langle \nabla\Delta\hat\d(t,\d_0)(\b_0), \nabla \partial_{t}\d\rangle\notag\\
=&-2|\nabla\Delta\hat\d(t,\d_0)(\b_0)|_2^2\notag\\
	&+2Q(t)\qv{\nabla\Delta\hat\d(t,\d_0)(\b_0),B_2(\nabla\hat\u(t,\v_0)(\u_0),\d(t,\d_0))}\notag\\
	&+2Q(t)\qv{\nabla\Delta\hat\d(t,\d_0)(\b_0),B_2(\hat\u(t,\v_0)(\u_0),\nabla\d(t,\d_0))}\notag\\
	&+2Q(t)\qv{\nabla\Delta\hat\d(t,\d_0)(\b_0),B_2(\nabla(\u(t,\v_0)+Z(t)),\hat\d(t,\d_0)(\b_0))}\notag\\
	&+2Q(t)\qv{\nabla\Delta\hat\d(t,\d_0)(\b_0),B_2(\u(t,\v_0)+Z(t),\nabla\hat\d(t,\d_0)(\b_0))}\notag\\
&+2Q(t)\qv{\nabla\Delta\hat\d(t,\d_0)(\b_0),\nabla(f'(\d)\hat{\d}(t,\d_0)(\b_0) )}\notag\\
\end{align}
By H\"older's inequality and Young's inequality, we have
\begin{align}
	&\frac{|\Delta\hat\d(t,\d_0)(\b_0)|_2^2}{dt}+|\nabla\Delta\hat\d(t,\d_0)(\b_0)|_2^2\notag\\
	\leq&c|Q|_\infty^2\norm{\d(t,\d_0)}_2^2|\nabla\hat\u(t,\v_0)(\u_0)|_2^2+c|Q|_\infty^2\norm{\d(t,\d_0)}_3^2|\hat\u(t,\v_0)(\u_0)|_2^2\notag\\
	&+c|Q|_\infty^2\norm{\u(t,\v_0)+Z(t)}_2^2|\nabla\hat\d(t,\d_0)(\b_0)|_2^2+c|Q|_\infty^2\norm{\u(t,\v_0)+Z(t)}_1^2|\Delta\hat\d(t,\d_0)(\b_0)|_2^2\notag\\
	&+c(1+\norm{\d(t,\d_0)}_1)^2\norm{\d(t,\d_0)}_2^2|\hat\d(t,\d_0)(\b_0)|_2^2+c(1+\norm{\d(t,\d_0)}_1^2)^2|\nabla\hat\d(t,\d_0)(\b_0)|_2^2
\end{align}
Applying Gronwall's inequality, and with the estimate in \eqref{eqn:3.12}, we get
\begin{align}\label{eqn:3.15}
	&\sup_{0\leq t\leq T}|\Delta\hat\d(t,\d_0)(\b_0)|_2^2+\int_0^T|\nabla\Delta\hat\d(t,\d_0)(\b_0)|_2^2dt\notag\\
	\leq&\bigg(|\Delta\hat\d(0,\d_0)(\b_0)|_2^2+c|Q|_\infty^2\sup_{0\leq t\leq T}\norm{\d(t,\d_0)}_2^2\int_0^T|\nabla\hat\u(t,\v_0)(\u_0)|_2^2dt\notag\\
	&\hspace{5mm}+c|Q|_\infty^2\sup_{0\leq t\leq T}|\hat\u(t,\v_0)(\u_0)|_2^2\int_0^T\norm{\d(t,\d_0)}_3^2dt\notag\\
	&\hspace{5mm}+c|Q|_\infty^2\sup_{0\leq t\leq T}|\nabla\hat\d(t,\d_0)(\b_0)|_2^2\int_0^T\norm{\u(t,\v_0)+Z(t)}_2^2dt\notag\\
	&\hspace{5mm}+c\sup_{0\leq t\leq T}(1+\norm{\d(t,\d_0)}_1)^2\norm{\d(t,\d_0)}_2^2|\hat\d(t,\d_0)(\b_0)|_2^2T\notag\\
	&\hspace{5mm}+c\sup_{0\leq t\leq T}(1+\norm{\d(t,\d_0)}_1^2)^2|\nabla\hat\d(t,\d_0)(\b_0)|_2^2T\bigg)\notag\\
	&\times\exp\{c|Q|_\infty^2\sup_{0\leq t\leq T}\norm{\u(t,\v_0)+Z(t)}_1^2T\}<\infty.
\end{align}
Multiplying \eqref{eqn:3.6} with $-\Delta\hat\u(t,\v_0)(\u_0)$ and integrating over $\mathbf{D}$ yields that
\begin{align}
	\frac{d|\nabla\hat\u(t,\v_0)(\u_0)|_2^2}{dt}=&-2|\Delta\hat\u(t,\v_0)(\u_0)|_2^2\notag\\
		&+2Q(t)\qv{B_1(\hat\u(t,\v_0)(\u_0),\u(t,\v_0)+Z(t)),\Delta\hat\u(t,\v_0)(\u_0)}\notag\\
	&+2Q(t)\qv{B_1(\u(t,\v_0)+Z(t),\hat\u(t,\v_0)(\u_0)),\Delta\hat\u(t,\v_0)(\u_0)}\notag\\
	&+2Q(t)^{-1}\qv{M(\hat\d(t,\d_0)(\b_0),\d(t,\d_0)),\Delta\hat\u(t,\v_0)(\u_0)}\notag\\
	&+2Q(t)^{-1}\qv{M(\d(t,\d_0),\hat\d(t,\d_0)(\b_0)),\Delta\hat\u(t,\v_0)(\u_0)}\notag\\
	=:&j_1+\cdots+j_5.
\end{align}
Applying H\"older's inequality and Young's inequality gives that
\begin{align*}
	j_2+j_3\leq&\eps|\Delta\hat\u(t,\v_0)(\u_0)|_2^2+c|Q|_\infty^2\norm{u(t,\v_0)+Z(t)}_2^2|\hat\u(t,\v_0)(\u_0)|_2^2\\
	&+c|Q|_\infty^2\norm{u(t,\v_0)+Z(t)}_1^2|\nabla\hat\u(t,\v_0)(\u_0)|_2^2.
\end{align*}
By Proposition \ref{prop:2.5}, H\"older's inequality and Young's inequality,
\begin{align*}
	j_4+j_5=&2Q(t)^{-1}\qv{B_2(\Delta\hat\u(t,\v_0)(\u_0),\hat\d(t,\d_0)(\b_0)),\Delta\d(t,\d_0)}\\
	&+2Q(t)^{-1}\qv{B_2(\Delta\hat\u(t,\v_0)(\u_0),\d(t,\d_0)),\Delta\hat\d(t,\d_0)(\b_0)}\\
	\leq&2|Q^{-1}|_\infty|\Delta\d(t,\d_0)|_\infty|\Delta\hat\u(t,\v_0)(\u_0)|_2|\hat\d(t,\d_0)(\b_0)|_2\\
	&+2|Q^{-1}|_\infty|\d(t,\d_0)|_\infty|\Delta\hat\u(t,\v_0)(\u_0)|_2|\Delta\hat\d(t,\d_0)(\b_0)|_2\\
	\leq&\eps|\Delta\hat\u(t,\v_0)(\u_0)|_2^2+c|Q^{-1}|_\infty^2\norm{\d(t,\d_0)}_3^2|\hat\d(t,\d_0)(\b_0)|_2^2\\
	&+c|Q^{-1}|_\infty^2\norm{\d(t,\d_0)}_1^2|\Delta\hat\d(t,\d_0)(\b_0)|_2^2.
\end{align*}
Altogether, we have
\begin{align}
	&\frac{d|\nabla\hat\u(t,\v_0)(\u_0)|_2^2}{dt}+|\Delta\hat\u(t,\v_0)(\u_0)|_2^2\notag\\
	\leq&c|Q|_\infty^2\norm{u(t,\v_0)+Z(t)}_2^2|\hat\u(t,\v_0)(\u_0)|_2^2+c|Q|_\infty^2\norm{u(t,\v_0)+Z(t)}_1^2|\nabla\hat\u(t,\v_0)(\u_0)|_2^2\notag\\
	&+c|Q^{-1}|_\infty^2\norm{\d(t,\d_0)}_3^2|\hat\d(t,\d_0)(\b_0)|_2^2+c|Q^{-1}|_\infty^2\norm{\d(t,\d_0)}_1^2|\Delta\hat\d(t,\d_0)(\b_0)|_2^2.
\end{align}
According to \eqref{eqn:3.12} and $\eqref{eqn:3.15}$, $\sup\limits_{0\leq t\leq T}\norm{\hat\d(t,\d_0)(\b_0)}_2^2<\infty$, now applying Gronwall inequality to the above estimate yields that
\begin{align}\label{eqn:3.18}
	&\sup_{0\leq t\leq T}|\nabla\hat\u(t,\v_0)(\u_0)|_2^2+\int_0^T|\Delta\hat\u(t,\v_0)(\u_0)|_2^2dt\notag\\
	\leq&\bigg(\norm{\u_0}_1^2+c|Q|_\infty^2\sup_{0\leq t\leq T}|\hat\u(t,\v_0)(\u_0)|_2^2\int_0^T\norm{u(t,\v_0)+Z(t)}_2^2dt\notag\\
	&\hspace{5mm}+c|Q^{-1}|_\infty^2\sup_{0\leq t\leq T}|\hat\d(t,\d_0)(\b_0)|_2^2\int_0^T\norm{\d(t,\d_0)}_3^2dt+c|Q^{-1}|_\infty^2\sup_{0\leq t\leq T}\norm{\d(t,\d_0)}_1^2|\Delta\hat\d(t,\d_0)(\b_0)|_2^2T\bigg)\notag\\
	&\times\exp\{c|Q|_\infty^2\sup_{0\leq t\leq T}\norm{\u(t,\v_0)+Z(t)}_1^2\}<\infty.
\end{align}
Since $\hat\u(t,\v_0)(\u_0),\hat\d(t,\d_0)(\b_0)$ are linear with respect to $\u_0,\b_0$, respectively. The above estimates \eqref{eqn:3.12}, \eqref{eqn:3.15} and \eqref{eqn:3.18} imply that $(\hat\u(t,\v_0)(\u_0),\hat\d(t,\d_0)(\b_0))\in L(\V\times\dH^2)$ for any $t\in[0,T]$, and $\hat\u(t,\v_0)(\cdot)\in L(\V,L^2([0,T];\vH^2))$, $\hat\d(t,\d_0)(\cdot)\in L(\dH^2,L^2([0,T];\vH^3))$, that is,
\begin{align}\label{eqn:3.19}
	&\sup_{0\leq t\leq T}[\norm{\hat\u(t,\v_0)(\u_0)}_{L(\V)}^2+\norm{\hat\d(t,\d_0)(\b_0)}_{L(\dH^2)}^2]\notag\\
	&+\int_0^T\norm{\hat\u(t,\v_0)(\u_0)}_{L(\V,L^2([0,T];\vH^2))}^2dt++\int_0^T\norm{\hat\d(t,\d_0)(\b_0)}_{L(\dH^2,L^2([0,T];\dH^3))}^2dt\leq c(\u_0,\d_0,Q,Z,T).
\end{align}
Now we will show $(\v_0,\d_0)\mapsto(\u(t,\v_0,\omega),\d(t,\d_0,\omega))$ has continuous Fr\'echet derivatives given by
$$D\u(t,\v_0,\omega)=\hat\u(t,\v_0,\omega)(\cdot),\ D\d(t,\d_0,\omega)=\hat\d(t,\d_0,\omega)(\cdot).$$
Thus, it suffices to show
\begin{align}
\label{eqn:3.20}
	\lim_{h\to0}\sup_{\norm{\u_0}_1+\norm{\b_0}_2\leq1}&\bigg\{\bigg\lVert\frac{\u(t,\v_0+h\u_0,\omega)-\u(t,\v_0,\omega)}{h}-\hat\u(t,\v_0)(\u_0)\bigg\rVert_1\notag\\
	&+\bigg\lVert\frac{\d(t,\d_0+h\b_0,\omega)-\d(t,\d_0,\omega)}{h}-\hat\d(t,\d_0)(\b_0)\bigg\rVert_2\bigg\}=0,
\end{align}
and the map $(\v_0,\d_0)\mapsto(\hat\u(t,\v_0,\omega),\hat\d(t,\d_0,\omega))$ is continuous.

We use the equations satisfied by $(\u(t,\v_0+h\u_0,\omega),\d(t,\d_0+h\b_0,\omega))$ and $(\u(t,\v_0,\omega),\d(t,\d_0,\omega))$, and denote by $\bar\u(t,\omega)=\u(t,\v_0+h\u_0)-\u(t,\v_0)$, $\bar\d(t,\omega)=\d(t,\d_0+h\b_0)-\d(t,\d_0)$. Then $\bar\u(t,\omega),\bar\d(t,\omega)$ satisfy the following equations:
\begin{align}
\label{eqn:3.21}
	\bar\u(t,\omega)=&h\u_0-\int_0^tA_1\bar\u(s,\omega)ds-\int_0^tQ(s)B_1(\u(s,\v_0+h\u_0)+Z(s),\bar\u(s,\omega))ds\notag\\
	&-\int_0^tQ(s)B_1(\bar\u(s,\omega),\u(s,\v_0))+Z(s))ds-\int_0^tQ^{-1}(s)M(\d(s,\d_0+h\b_0),\bar\d(s,\omega))ds\notag\\
	&-\int_0^tQ^{-1}(s)M(\bar\d(s,\omega),\d(s,\d_0))ds.\\
	\label{eqn:3.22}
	\bar\d(t,\omega)=&h\b_0-\int_0^tA_2\bar\d(s,\omega)ds-\int_0^tQ(s)B_2(\u(s,\v_0+h\u_0)+Z(s),\bar\d(s,\omega))ds\notag\\
	&-\int_0^tQ(s)B_2(\bar\u(s,\omega),\d(s,\d_0))ds-\int_0^tf(\d(s,\d_0+h\b_0))ds+\int_0^tf(\d(s,\d_0))ds.
\end{align}
We first take inner product of \eqref{eqn:3.21} with $\bar\u(t,\omega)$ in $L^{2}(\mathbf{D} )$,
\begin{align}
	|\bar\u(t,\omega)|_2^2=&h^2|\u_0|_2^2-2\int_0^t|\nabla\bar\u(s,\omega)|_2^2ds\notag\\
	&-2\int_0^tQ(s)\qv{B_1(\u(s,\v_0+h\u_0)+Z(s),\bar\u(s,\omega)),\bar\u(s,\omega)}ds\notag\\
	&-2\int_0^tQ(s)\qv{B_1(\bar\u(s,\omega),\u(s,\v_0))+Z(s)),\bar\u(s,\omega)}ds\notag\\
	&-2\int_0^tQ^{-1}(s)\qv{M(\d(s,\d_0+h\b_0),\bar\d(s,\omega)),\bar\u(s,\omega)}ds\notag\\
	&-2\int_0^tQ^{-1}(s)\qv{M(\bar\d(s,\omega),\d(s,\d_0)),\bar\u(s,\omega)}ds\notag\\
	=:&k_1+\cdots+k_6.
\end{align}
By Lemma \ref{lemma:2.1}, H\"older's inequality and Young's inequality, $k_3=0$ and
\begin{align*}
	k_4\leq&c|Q|_\infty\int_0^t|\nabla(\u(s,\v_0)+Z(s))|_\infty|\bar\u(s,\omega)|_2^2ds\\
	\leq&c|Q|_\infty\int_0^t\norm{\u(s,\v_0)+Z(s)}_2\bar\u(s,\omega)|_2^2ds
\end{align*}
By H\"older's inequality and Young's inequality, we have
\begin{align*}
	k_5+k_6
	\leq&|Q^{-1}|_\infty\int_0^t|\nabla\d(s,\d_0+h\b_0)|_\infty|\nabla\bar\u(s,\omega)|_2|\nabla\bar\d(s,\omega)|_2ds\\
	&+|Q^{-1}|_\infty\int_0^t|\nabla\d(s,\d_0)|_\infty|\nabla\bar\u(s,\omega)|_2|\nabla\bar\d(s,\omega)|_2ds\\
	\leq&\eps\int_0^t|\nabla\bar\u(s,\omega)|_2^2ds+c|Q^{-1}|^2_\infty\int_0^t\norm{\d(s,\d_0+h\b_0)}|_2^2|\nabla\bar\d(s,\omega)|_2^2ds\\
	&+c|Q^{-1}|^2_\infty\int_0^t\norm{\d(s,\d_0)}_2^2|\nabla\bar\d(s,\omega)|_2^2ds.
\end{align*}
Altogether, we get
\begin{align}\label{eqn:3.24}
	&|\bar\u(t,\omega)|_2^2+\int_0^t|\nabla\bar\u(s,\omega)|_2^2ds-h^2|\u_0|_2^2\notag\\
	\leq&c|Q|_\infty\int_0^t\norm{\u(s,\v_0)+Z(s)}_2|\bar\u(s,\omega)|_2^2ds+c|Q^{-1}|^2_\infty\int_0^t\norm{\d(s,\d_0+h\b_0)}_2^2|\nabla\bar\d(s,\omega)|_2^2ds\notag\\
	&+c|Q^{-1}|^2_\infty\int_0^t\norm{\d(s,\d_0)}_2^2|\nabla\bar\d(s,\omega)|_2^2ds.
\end{align}
Similarly to derive (\ref{eqn:3.10}), we taking inner product between $\partial_t\bar\d(t,\omega)$ and $-\Delta\bar\d(t,\omega)+\bar\d(t,\omega)$ in $L^{2}(\mathbf{D})$ obtaining that
\begin{align}
		&\qv{\partial_t\bar\d(t,\omega),-\Delta\bar\d(t,\omega)+\bar\d(t,\omega)}\notag\\
		=&\frac{1}{2}\frac{d|\nabla\bar\d(t,\omega)|_2^2}{dt}+\frac{1}{2}\frac{d|\bar\d(t,\omega)|_2^2}{dt}\notag\\
		=&-\qv{-\Delta\bar\d(t,\omega)+\bar\d(t,\omega),-\Delta\bar\d(t,\omega)+\bar\d(t,\omega)}\notag\\
		&-\qv{-\Delta\bar\d(t,\omega)+\bar\d(t,\omega),Q(t)B_2(\u(t,\v_0+h\u_0)+Z(t),\bar\d(t,\omega))}\notag\\
		&-\qv{-\Delta\bar\d(t,\omega)+\bar\d(t,\omega),Q(t)B_2(\bar\u(t,\omega),\d(t,\d_0))}\notag\\
		&-\qv{-\Delta\bar\d(t,\omega)+\bar\d(t,\omega),f(\d(t,\d_0+h\b_0))-f(\d(t,\d_0))-\bar\d(t,\omega)}.
\end{align}
By H\"older's inequality and Young's inequality, we obtain that
\begin{align}\label{eqn:3.26}
	&\frac{d|\nabla\bar\d(t,\omega)|_2^2}{dt}+\frac{d|\bar\d(t,\omega)|_2^2}{dt}+|-\Delta\bar\d(t,\omega)+\bar\d(t,\omega)|_2^2\notag\\
	\leq&c|Q|_\infty^2\norm{\u(t,\v_0+h\u_0)+Z(t)}_1^2|\bar\d(t,\omega)|_2^2+c|Q|^2_\infty\norm{\d(t,\d_0)}_1^2|\bar\u(t,\omega)|_2^2\notag\\
	&+c|\bar\d(t,\omega)|_2^2+c\norm{\d(t,\d_0+h\b_0)}_1^4|\bar\d(t,\omega)|_2^2\notag\\
	&+c(\norm{\d(t,\d_0+h\b_0)}_1^2+\norm{\d(t,\d_0)}_1^2)\norm{\d(t,\d_0)}_1^2|\bar\d(t,\omega)|_2^2.
\end{align}
Combining \eqref{eqn:3.24} and \eqref{eqn:3.26}, and applying Gronwall's inequality gives that
\begin{align}\label{eqn:3.27}
	&\sup_{0\leq t\leq T}[|\bar\u(t,\omega)|_2^2+|\nabla\bar\d(t,\omega)|_2^2+|\bar\d(t,\omega)|_2^2]+\int_0^T|\nabla\bar\u(t,\omega)|_2^2dt+\int_0^T|-\Delta\bar\d(t,\omega)+\bar\d(t,\omega)|_2^2dt\notag\\
	\leq&h^2(|\u_0|_2^2+\norm{\b_0}_1^2)\exp c\bigg\{T+|Q|_\infty\int_0^T\norm{\u(t,\v_0)+Z(t)}_2dt+|Q^{-1}|^2_\infty\sup_{0\leq t\leq T}\norm{\d(t,\d_0+h\b_0)}_2^2T\notag\\
	&\hspace{37mm}+|Q^{-1}|^2_\infty\sup_{0\leq t\leq T}\norm{\d(t,\d_0)}_2^2T+|Q|_\infty^2\sup_{0\leq t\leq T}\norm{\u(t,\v_0+h\u_0)+Z(t)}_1^2T\notag\\
	&\hspace{37mm}+|Q|^2_\infty\sup_{0\leq t\leq T}\norm{\d(t,\d_0)}_1^2T+\sup_{0\leq t\leq T}\norm{\d(t,\d_0+h\b_0)}_1^4T\notag\\
	&\hspace{37mm}+\sup_{0\leq t\leq T}(\norm{\d(t,\d_0+h\b_0)}_1^2+\norm{\d(t,\d_0)}_1^2)\norm{\d(t,\d_0)}_1^2T\bigg\}\notag\\
	=:&h^2(|\u_0|_2^2+\norm{\b_0}_1^2)g_1(T).
\end{align}
Similarly to derive (\ref{eqn:3.14a}), we take inner product between $\Delta\bar\d(t,\omega)$ and $\Delta\partial_t\bar\d(t,\omega)$ in $L^2(\mathbf{D})$ obtaining that
\begin{align}
	 \frac{1}{2}\partial_t|\Delta\bar\d(t,\omega)|_2^2  =&\qv{\Delta\bar\d(t,\omega),\Delta\partial_t\bar\d(t,\omega)}=-\qv{\nabla\Delta\bar\d(t,\omega),\nabla\partial_t\bar\d(s,\omega)}\notag\\
	=&-|\nabla\Delta\bar\d(t,\omega)|_2^2+Q(t)\qv{\nabla\Delta\bar\d(t,\omega),B_2(\nabla(\u(s,\v_0+h\u_0)+Z(t)),\bar\d(t,\omega))}\notag\\
	&+Q(t)\qv{\nabla\Delta\bar\d(t,\omega),B_2(\u(s,\v_0+h\u_0)+Z(t),\nabla\bar\d(t,\omega))}\notag\\
	&+Q(t)\qv{\nabla\Delta\bar\d(t,\omega),B_2(\nabla\bar\u(t,\omega),\d(t,\d_0))}\notag\\
	&+Q(t)\qv{\nabla\Delta\bar\d(t,\omega),B_2(\bar\u(t,\omega),\nabla\d(t,\d_0))}\notag\\
	&+\qv{\nabla\Delta\bar\d(t,\omega),\nabla f(\d(t,\d_0+h\b_0)-\nabla f(\d(t,\d_0))}.
\end{align}
Applying H\"older's inequality and Young's inequality gives that
\begin{align}
	&\frac{d|\Delta\bar\d(t,\omega)|_2^2}{dt}+|\nabla\Delta\bar\d(t,\omega)|_2^2\notag\\
	\leq&c|Q|_\infty^2\norm{\u(s,\v_0+h\u_0)+Z(t)}_2|\nabla\bar\d(t,\omega)|_2^2\notag\\
	&+c|Q|_\infty^2\norm{\u(s,\v_0+h\u_0)+Z(t)}_1|\Delta\bar\d(t,\omega)|_2^2\notag\\
	&+c|Q|_\infty^2\norm{\d(t,\d_0)}_2^2|\nabla\bar\u(t,\omega)|_2^2+c|Q|_\infty^2\norm{\d(t,\d_0)}_3^2|\bar\u(t,\omega)|_2^2\notag\\
	&+c(\norm{\d(t,\d_0+h\b_0)}_1^4+\norm{\d(t,\d_0)}_1^4)|\nabla\bar\d(t,\omega)|_2^2\notag\\
	&+c(\norm{\d(t,\d_0+h\b_0)}_1^2+\norm{\d(t,\d_0)}_1^2)(\norm{\d(t,\d_0+h\b_0)}_2^2+\norm{\d(t,\d_0)}_2^2)|\bar\d(t,\omega)|_2^2.
\end{align}
Applying Gronwall's inequality, and together with \eqref{eqn:3.27}, yields that
\begin{align}
	&\sup_{0\leq t\leq T}|\Delta\bar\d(t,\omega)|_2^2+\int_0^T|\nabla\Delta\bar\d(t,\omega)|_2^2dt\notag\\
	\leq& ch^2(|\u_0|_2^2+\norm{\b_0}_2^2)g_1(T)\notag\\
	&\times\bigg[|Q|_\infty^2\int_0^T\norm{\u(s,\v_0+h\u_0)+Z(t)}_2dt+|Q|_\infty^2\sup_{0\leq t\leq T}\norm{\u(s,\v_0+h\u_0)+Z(t)}_1^2T\notag\\
	&\hspace{5mm}+|Q|_\infty^2\int_0^T\norm{\d(t,\d_0)}_3^2dt+\sup_{0\leq t\leq T}(\norm{\d(t,\d_0+h\b_0)}_1^4+\norm{\d(t,\d_0)}_1^4)T\notag\\
	&\hspace{5mm}+\sup_{0\leq t\leq T}(\norm{\d(t,\d_0+h\b_0)}_1^2+\norm{\d(t,\d_0)}_1^2)(\norm{\d(t,\d_0+h\b_0)}_2^2+\norm{\d(t,\d_0)}_2^2)T\bigg]\notag\\
	&\times\exp\{c|Q|_\infty^2\sup_{0\leq t\leq T}\norm{\u(s,\v_0+h\u_0)+Z(t)}_1^2T\}.
\end{align}
Now we multiply \eqref{eqn:3.21} with $-\Delta\bar\u(t,\omega)$ and integrating over $D$ gives that
\begin{align}\label{eqn:3.30}
	\frac{d|\nabla\bar\u(t,\omega)|_2^2}{dt}=&-2|\Delta\bar\u(t,\omega)|_2^2\notag\\
	&+2Q(t)\qv{B_1(\u(t,\v_0+h\u_0)+Z(t),\bar\u(t,\omega)),\Delta\bar\u(t,\omega)}\notag\\
	&+2Q(t)\qv{B_1(\bar\u(t,\omega),\u(t,\v_0)+Z(t)),\Delta\bar\u(t,\omega)}\notag\\
	&+2Q^{-1}(t)\qv{M(\d(t,\d_0+h\b_0),\bar\d(t,\omega)),\Delta\bar\u(t,\omega)}\notag\\
	&+2Q^{-1}(t)\qv{M(\bar\d(t,\omega),\d(t,\d_0)),\Delta\bar\u(t,\omega)}\notag\\
	=:l_1+\cdots+l_5.
\end{align}
By H\"older's inequality and Young's inequality, we get
\begin{align*}
	l_2+l_3\leq&|Q|_\infty|\u(t,\v_0+h\u_0)+Z(t)|_\infty|\nabla\bar\u(t,\omega)|_2|\Delta\bar\u(t,\omega)|_2\\
	&+|Q|_\infty|\nabla(\u(t,\v_0)+Z(t))|_\infty|\bar\u(t,\omega)|_2|\Delta\bar\u(t,\omega)|_2\\
	\leq&\eps|\Delta\bar\u(t,\omega)|_2^2+c|Q|^2_\infty\norm{\u(t,\v_0+h\u_0)+Z(t)}_1^2|\nabla\bar\u(t,\omega)|_2^2\\
	&+c|Q|^2_\infty\norm{\u(t,\v_0)+Z(t)}_2^2|\bar\u(t,\omega)|_2^2
\end{align*}
By Lemma \ref{lemma:2.4} and Proposition \ref{prop:2.5}, then applying Lemma \ref{lemma:2.2}, H\"older's inequality and Young's inequality, we have
\begin{align*}
l_4+l_5=&2Q^{-1}(t)\qv{B_2(\Delta\bar\u(t,\omega),\bar\d(t,\omega)),\Delta\d(t,\d_0+h\b_0)}\\
&+2Q^{-1}(t)\qv{B_2(\Delta\bar\u(t,\omega),\d(t,\d_0)),\Delta\bar\d(t,\omega)}\\
\leq&|Q^{-1}|_\infty|\Delta\d(t,\d_0+h\b_0)|_\infty|\Delta\bar\u(t,\omega)|_2|\nabla\bar\d(t,\omega)|_2\\
&+|Q^{-1}|_\infty|\nabla\d(t,\d_0)|_\infty|\Delta\bar\u(t,\omega)|_2|\Delta\bar\d(t,\omega)|_2\\
\leq&\eps|\Delta\bar\u(t,\omega)|_2^2+c|Q^{-1}|^2_\infty\norm{\d(t,\d_0+h\b_0)}_3^2|\nabla\bar\d(t,\omega)|_2^2\\
&+c|Q^{-1}|^2_\infty\norm{\d(t,\d_0)}_2^2|\Delta\bar\d(t,\omega)|_2^2.
\end{align*}
Hence, one can get
\begin{align}\label{eqn:3.32}
	&\frac{d|\nabla\bar\u(t,\omega)|_2^2}{dt}+|\Delta\bar\u(t,\omega)|_2^2\notag\\
	\leq&c|Q|^2_\infty\norm{\u(t,\v_0+h\u_0)+Z(t)}_1^2|\nabla\bar\u(t,\omega)|_2^2+c|Q|^2_\infty\norm{\u(t,\v_0)+Z(t)}_2^2|\bar\u(t,\omega)|_2^2\notag\\
	&+c|Q^{-1}|^2_\infty\norm{\d(t,\d_0+h\b_0)}_3^2|\nabla\bar\d(t,\omega)|_2^2+c|Q^{-1}|^2_\infty\norm{\d(t,\d_0)}_2^2|\Delta\bar\d(t,\omega)|_2^2.
\end{align}
According to \eqref{eqn:3.27} and \eqref{eqn:3.30},
\begin{equation}
	\sup_{0\leq t\leq T}[|\bar\u(t,\omega)|_2^2+|\nabla\bar\d(t,\omega)|_2^2+|\Delta\bar\d(t,\omega)|_2^2]\leq h^2(|\u_0|_2^2+\norm{\b_0}_2^2)g_2(T).
\end{equation}
Applying Gronwall's inequality to \eqref{eqn:3.32}, we obtain that
\begin{align}\label{eqn:3.34}
	&\sup_{0\leq t\leq T}|\nabla\bar\u(t,\omega)|_2^2+\int_0^T|\Delta\bar\u(t,\omega)|_2^2dt\notag\\
	\leq&ch^2(|\u_0|_2^2+\norm{\b_0}_2^2)g_2(T)\notag\\
	&\times\bigg[|Q|^2_\infty\int_0^T\norm{\u(t,\v_0)+Z(t)}_2^2dt+|Q^{-1}|^2_\infty\int_0^T\norm{\d(t,\d_0+h\b_0)}_3^2dt+|Q^{-1}|^2_\infty\sup_{0\leq t\leq T}\norm{\d(t,\d_0)}_2^2T\bigg]\notag\\
	&\times\exp\{c|Q|^2_\infty\sup_{0\leq t\leq T}\norm{\u(t,\v_0+h\u_0)+Z(t)}_1^2T\}.
\end{align}
With the estimates \eqref{eqn:3.27}, \eqref{eqn:3.30} and \eqref{eqn:3.34}, we arrive at the conclusion that for all $\v_0,\u_0\in\V, \d_0,\b_0\in\dH^2$, and any $h\in\R$,
\begin{align}
\label{eqn:3.35}
	\lim_{h\to0}\sup_{\norm{\u_0}_1+\norm{\b_0}_2\leq1}\bigg\{&\sup_{0\leq t\leq T}[\norm{\u(t,\v_0+h\u_0)-\u(t,\v_0)}_1^2+\norm{\d(t,\d_0+h\b_0)-\d(t,\d_0)}_2^2]\notag\\
	&+\int_0^T\norm{\u(t,\v_0+h\u_0)-\u(t,\v_0)}_2^2dt+\int_0^T\norm{\d(t,\d_0+h\b_0)-\d(t,\d_0)}_3^2dt
	\bigg\}=0.
\end{align}
Now we set for $t\in[0,T], h\in\R\setminus\{0\}$,
$$U(t,\v_0,\u_0,h)=\frac{\u(t,\v_0+h\u_0,\omega)-\u(t,\v_0,\omega)}{h},\ X(t,\v_0,\u_0,h)=U(t,\v_0,\u_0,h)-\hat\u(t,\v_0)(\u_0);$$
$$D(t,\d_0,\b_0,h)=\frac{\d(t,\d_0+h\b_0,\omega)-\d(t,\d_0,\omega)}{h},\ Y(t,\d_0,\b_0,h)=D(t,\d_0,\b_0,h)-\hat\d(t,\d_0)(\b_0).$$
Then
\begin{align}\label{eqn:3.36}
	X(t,\v_0,\u_0,h)=&-\int_0^tA_1X(s,\v_0,\u_0,h)ds-\int_0^tQ(s)B_1(\u(s,\v_0)+Z(s),X(s,\v_0,\u_0,h))ds\notag\\
	&-\int_0^tQ(s)B_1(X(s,\v_0,\u_0,h),\u(s,\v_0+h\u_0)+Z(s))ds\notag\\
	&-\int_0^tQ(s)B_1(\hat\u(s,\v_0)(\u_0),\u(s,\v_0+h\u_0)-\u(s,\v_0))ds\notag\\
	&-\int_0^tQ(s)^{-1}M(\d(s,\d_0),Y(s,\d_0,\b_0,h))ds-\int_0^tQ(s)^{-1}M(Y(s,\d_0,\b_0,h),\d(s,\d_0+h\b_0))ds\notag\\
	&-\int_0^tQ(s)^{-1}M(\hat\d(s,\d_0)(\b_0),\d(s,\d_0+h\b_0)-\d(s,\d_0))ds.
\end{align}
\begin{align}\label{eqn:3.37}
Y(t,\d_0,\b_0,h)=&-\int_0^tA_2Y(s,\d_0,\b_0,h)ds-\int_0^tQ(s)B_2(X(s,\v_0,\u_0,h),\d(s,\d_0))ds\notag\\
&-\int_0^tQ(s)B_2(\u(s,\v_0+h\u_0)+Z(s),Y(s,\d_0,\b_0,h))ds\notag\\
&-\int_0^tQ(s)B_2(\u(s,\v_0+h\u_0)-\u(s,\u_0),\hat\d(s,\d_0)(\b_0))ds\notag\\
&-\int_0^t\frac{1}{h}[f(\d(s,\d_0+h\b_0))-f(\d(s,\d_0))]ds+\int_0^t f'(\d(s,\d_0))\hat\d(s,\d_0)(\b_0)ds.
\end{align}
We first multiply \eqref{eqn:3.36} with $X(t,\v_0,\u_0,h)$, and integrate over $\mathbf{D}$,
\begin{align}
	|X(t)|_2^2=&-2\int_0^t|\nabla X(s)|_2^2ds-2\int_0^tQ(s)\qv{B_1(\u(s,\v_0)+Z(s),X(s)),X(s)}ds\notag\\
	&-2\int_0^tQ(s)\qv{B_1(X(s),\u(s,\v_0+h\u_0)+Z(s)),X(s)}ds\notag\\
	&-2\int_0^tQ(s)\qv{B_1(\hat\u(s,\v_0)(\u_0),\u(s,\v_0+h\u_0)-\u(s,\v_0)),X(s)}ds\notag\\
	&-2\int_0^tQ(s)^{-1}\qv{M(\d(s,\d_0),Y(s)),X(s)}ds\notag\\
	&-2\int_0^tQ(s)^{-1}\qv{M(Y(s),\d(s,\d_0+h\b_0)),X(s)}ds\notag\\
	&-2\int_0^tQ(s)^{-1}\qv{M(\hat\d(s,\d_0)(\b_0),\d(s,\d_0+h\b_0)-\d(s,\d_0)),X(s)}ds.
\end{align}
By Lemma \ref{lemma:2.1}, H\"older's inequality and Young's inequality, we get
\begin{align}
\label{eqn:3.39}
	&|X(t)|_2^2+\int_0^t|\nabla X(s)|_2^2ds\notag\\
	\leq&c|Q|_\infty\int_0^t\norm{\u(s,\v_0+h\u_0)+Z(s)}_2^{2}|X(s)|_2^2ds+c\int_0^t|\nabla\bar\u(s,\omega)|_2^2ds\notag\\
	&+c|Q|_\infty^2\int_0^t\norm{\hat\u(s,\v_0)(\u_0)}_1^2|X(s)|_2^2ds+c|Q^{-1}|_\infty^2\int_0^t\norm{\d(s,\d_0)}_2^2|\nabla Y(s)|_2^2ds\notag\\
	&+c|Q^{-1}|_\infty^2\int_0^t\norm{\d(s,\d_0+h\b_0)}_2^2|\nabla Y(s)|_2^2ds+c|Q^{-1}|_\infty^2\int_0^t\norm{\hat\d(s,\d_0)(\b_0)}_2^2|\nabla\bar\d(s,\omega)|_2^2ds.
\end{align}
Similarly to derive (\ref{eqn:3.10}), we take inner product between $\partial_tY(t)$ and $-\Delta Y(t)+Y(t)$ in $L^2(\mathbf{D})$ obtaining that
\begin{align}
	&\qv{\partial_tY(t),-\Delta Y(t)+Y(t)}\notag\\
	=&\frac{1}{2}\frac{d|\nabla Y(t)|_2^2}{dt}+\frac{1}{2}\frac{d|Y(t)|_2^2}{dt}\notag\\
	=&-|-\Delta Y(t)+Y(t)|_2^2-\qv{-\Delta Y(t)+Y(t),Q(t)B_2(X(t),\d(t,\d_0))}\notag\\
	&-\qv{-\Delta Y(t)+Y(t),Q(t)B_2(\u(t,\v_0+h\u_0)+Z(t),Y(t))}\notag\\
	&-\qv{-\Delta Y(t)+Y(t),Q(t)B_2(\u(t,\v_0+h\u_0)-\u(t,\v_0),\hat\d(t,\d_0)(\b_0))}\notag\\
	&-\qv{-\Delta Y(t)+Y(t),\frac{f(\d(t,\d_0+h\b_0))-f(\d(t,\d_0))}{h}-f'(\d(s,\d_0))\hat\d(t,\d_0)(\b_0)}\notag\\
	&+\qv{-\Delta Y(t)+Y(t),Y(t)}.
\end{align}
With similar discussion, we get
\begin{align}
\label{eqn:3.41}
	&\frac{d|\nabla Y(t)|_2^2}{dt}+\frac{d|Y(t)|_2^2}{dt}+|-\Delta Y(t)+Y(t)|_2^2\notag\\
	\leq&c|Q|_\infty^2\norm{\d(t,\d_0)}_2^2|X(t)|_2^2+c|Q|_\infty^2\norm{\u(t,\v_0+h\u_0)}_1^2|\nabla Y(t)|_2^2+c|Q|_\infty^2\norm{\bar\u(t,\omega)}_1^2|\nabla\hat\d(t,\d_0)(\b_0)|_2^2dt\notag\\
	&+c\sup_{h}\norm{\d(t,\d_0+h\b_0)}_1^4|Y(t)|_2^2+c|Y(t)|_2^2.
\end{align}
Combining \eqref{eqn:3.39} and \eqref{eqn:3.41}, and applying Gronwall's inequality,
\begin{align}
	&\sup_{0\leq t\leq T}[|X(t)|_2^2+|Y(t)|_2^2+|\nabla Y(t)|_2^2]+\int_0^T|\nabla X(t)|_2^2dt+\int_0^T|-\Delta Y(t)+Y(t)|_2^2dt\notag\\
	\leq&cT[\sup_{0\leq t\leq T}|\nabla\bar\u(t,\omega)|_2^2+|Q^{-1}|_\infty^2\sup_{0\leq t\leq T}\norm{\hat\d(t,\d_0)(\b_0)}_2^2|\nabla\bar\d(s,\omega)|_2^2\notag\\
	&\hspace{5mm}+\sup_{0\leq t\leq T}|Q|_\infty^2|\nabla\hat\d(t,\d_0)(\b_0)|_2^2\norm{\bar\u(t,\omega)}_1^2]\notag\\
	&\times\exp c\bigg\{|Q|_\infty\int_0^T\norm{\u(t,\v_0+h\u_0)+Z(t)}_2^2dt+|Q|_\infty^2\sup_{0\leq t\leq T}\norm{\hat\u(t,\v_0)(\u_0)}_1^2T\notag\\
	&\hspace{15mm}+|Q^{-1}|_\infty^2\sup_{0\leq t\leq T}(\norm{\d(t,\d_0+h\b_0)}_2^2+\norm{\d(t,\d_0)}_2^2)T\notag\\
	&\hspace{15mm}+|Q|_\infty^2\sup_{0\leq t\leq T}(\norm{\d(t,\d_0)}_2^2+\norm{\u(t,\v_0+h\u_0)}_1^2)T\notag\\
	&\hspace{15mm}+\sup_{0\leq t\leq T}\sup_{h}\norm{\d(t,\d_0+h\b_0)}_1^4T+T\bigg\}\notag\\
	=:&g_3(T)g_4(T).
\end{align}
Note that by \eqref{eqn:3.35}, $g_3(T)\to0$ as $h\to0$ and $g_4(T)<\infty$, this gives that
\begin{align}\label{eqn:3.43}
	\lim_{h\to0}\sup_{\norm{\u_0}_1+\norm{\b_0}_2\leq1}\bigg\{&\sup_{0\leq t\leq T}[|X(t)|_2^2+|Y(t)|_2^2+|\nabla Y(t)|_2^2]\notag\\
	&+\int_0^T|\nabla X(t)|_2^2dt+\int_0^T|-\Delta Y(t)+Y(t)|_2^2dt\bigg\}=0.
\end{align}
Now we take inner product between $\Delta Y(t)$ and $\Delta\partial_tY(t)$ in $L^{2}(\mathbf{D})$, one can get
\begin{align}
	\frac{1}{2}\partial_t|\Delta Y(t)|_2^2=&\qv{\Delta Y(t),\Delta\partial_tY(t)}=-\qv{\nabla\Delta Y(t),\nabla\partial_t Y(t)}\notag\\
	=&-|\nabla\Delta Y(t)|_2^2+Q(t)\qv{\nabla\Delta Y(t),B_2(\nabla X(t),\d(t,\d_0))}\notag\\
	&+Q(t)\qv{\nabla\Delta Y(t),B_2(X(t),\nabla \d(t,\d_0))}\notag\\
	&+Q(t)\qv{\nabla\Delta Y(t),B_2(\nabla(\u(t,\v_0+h\u_0)+Z(t)),Y(t))}\notag\\
	&+Q(t)\qv{\nabla\Delta Y(t),B_2(\u(t,\v_0+h\u_0)+Z(t),\nabla Y(t))}\notag\\
	&+Q(t)\qv{\nabla\Delta Y(t),B_2(\nabla\bar\u(t,\omega),\hat\d(t,\d_0)(\b_0))}\notag\\
	&+Q(t)\qv{\nabla\Delta Y(t),B_2(\bar\u(t,\omega),\nabla\hat\d(t,\d_0)(\b_0))}\notag\\
	&+\qv{\nabla\Delta Y(t),\nabla\{\frac{1}{h}[f(\d(t,\d_0+h\b_0))-f(\d(t,\d_0))]-f'(\d(t,\d_0))\hat\d(t,\d_0)(\b_0)\}}
\end{align}
By H\"older's inequality and Young's inequality,
\begin{align}
	&\frac{d|\Delta Y(t)|_2^2}{dt}+|\nabla\Delta Y(t)|_2^2\notag\\
	\leq&c|Q|_\infty^2\norm{\d(t,\d_0)}_2^2|\nabla X(t)|_2^2+c|Q|_\infty^2\norm{\d(t,\d_0)}_3^2|X(t)|_2^2\notag\\
	&+c|Q|_\infty^2\norm{\u(t,\v_0+h\u_0)+Z(t)}_2^2|\nabla Y(t)|_2^2+c|Q|_\infty^2\norm{\u(t,\v_0+h\u_0)+Z(t)}_1^2|\Delta Y(t)|_2^2\notag\\
	&+c|Q|_\infty^2\norm{\hat\d(t,\d_0)(\b_0)}_2^2|\nabla\bar\u(t,\omega)|_2^2+c|Q|_\infty^2|\Delta\hat\d(t,\d_0)(\b_0)|_2^2\norm{\bar\u(t,\omega)}_1^2\notag\\
	&+c\sup_h\norm{\d(t,\d_0+h\b_0)}_2^4|\nabla Y(t)|_2^2.
\end{align}
Applying Gronwall's inequality,
\begin{align}
	&\sup_{0\leq t\leq T}|\Delta Y(t)|_2^2+\int_0^T|\nabla\Delta Y(t)|_2^2dt\notag\\
	\leq&c\bigg[|Q|_\infty^2\sup_{0\leq t\leq T}\norm{\d(t,\d_0)}_2^2\int_0^T|\nabla X(t)|_2^2dt+|Q|_\infty^2\sup_{0\leq t\leq T}|X(t)|_2^2\int_0^T\norm{\d(t,\d_0)}_3^2dt\notag\\
	&\hspace{5mm}+|Q|_\infty^2\sup_{0\leq t\leq T}|\nabla Y(t)|_2^2\int_0^T\norm{\u(t,\v_0+h\u_0)+Z(t)}_2^2dt+|Q|_\infty^2\sup_{0\leq t\leq T}\norm{\hat\d(t,\d_0)(\b_0)}_2^2|\nabla\bar\u(t,\omega)|_2^2T\notag\\
	&\hspace{5mm}+|Q|_\infty^2\sup_{0\leq t\leq T}|\Delta\hat\d(t,\d_0)(\b_0)|_2^2\norm{\bar\u(t,\omega)}_1^2T+\sup_{0\leq t\leq T}\sup_h\norm{\d(t,\d_0+h\b_0)}_2^4|\nabla Y(t)|_2^2T\bigg]\notag\\
	&\times\exp\{c|Q|_\infty^2\norm{\u(t,\v_0+h\u_0)+Z(t)}_1^2T\}.
\end{align}
By \eqref{eqn:3.35} and \eqref{eqn:3.43},
\begin{equation}\label{eqn:3.47}
	\lim_{h\to0}\sup_{\norm{\u_0}_1+\norm{\b_0}_2\leq1}\bigg\{\sup_{0\leq t\leq T}|\Delta Y(t)|_2^2+\int_0^T|\nabla\Delta Y(t)|_2^2dt\bigg\}=0.
\end{equation}
Now we multiply $X(t,\v_0,\u_0,h)$ with $-\Delta X(t,\v_0,\u_0,h)$ and integrate over $\mathbf{D}$,
\begin{align}
	\frac{d|\nabla X(t)|_2^2}{dt}=&-2|\Delta X(t)|_2^2+2Q(t)\qv{B_1(\u(t,\v_0)+Z(t),X(t)),\Delta X(s)}\notag\\
	&+2Q(t)\qv{B_1(X(t),\u(t,\v_0+h\u_0)+Z(t)),\Delta X(s)}\notag\\
	&+2Q(t)\qv{B_1(\hat\u(t,\v_0)(\u_0),\u(t,\v_0+h\u_0)-\u(t,\v_0)),\Delta X(t)}\notag\\
	&+2Q(t)^{-1}\qv{M(\d(t,\d_0),Y(t)),\Delta X(t)}+2Q(t)^{-1}\qv{M(Y(t),\d(t,\d_0+h\b_0)),\Delta X(t)}\notag\\
	&+2Q(t)^{-1}\qv{M(\hat\d(t,\d_0)(\b_0),\d(t,\d_0+h\b_0)-\d(t,\d_0)),\Delta X(t)}\notag\\
	=&m_1+\cdots+m_7.
\end{align}
By Lemma \ref{lemma:2.1}, and Young's inequality, we get
\begin{align*}
	m_2\leq& \sup_{0\leq t\leq T}|Q(t)||\u(t,\v_0)+Z(t)|_\infty|\nabla X(t)|_2|\Delta X(t)|_2\\
	\leq& \eps|\Delta X(t)|_2^2+c|Q|_\infty^2\norm{\u(t,\v_0)+Z(t)}_1^2|\nabla X(t)|_2^2,
\end{align*}
and similarly
\begin{align*}
	m_3\leq&\eps|\Delta X(t)|_2^2+c|Q|_\infty^2\norm{\u(t,\v_0+h\u_0)+Z(t)}_2^2|X(t)|_2^2,
\end{align*}
and
\begin{align*}
	m_4\leq&\eps|\Delta X(t)|_2^2+c|Q|_\infty^2\norm{\hat\u(t,\v_0)(\u_0)}_1^2|\nabla\bar\u(t,\omega)|_2^2.
\end{align*}
By Proposition \ref{prop:2.5}, we get
\begin{align}
	m_5+m_6+m_7=&2Q(t)^{-1}\qv{B_2(\Delta X(t),\d(t,\d_0)),\Delta Y(t)}+2Q(t)^{-1}\qv{B_2(\Delta X(t),Y(t)),\Delta\d(t,\d_0+h b_{0})}\notag\\
	&+2\frac{1}{h}Q(t)^{-1}\qv{B_2(\Delta X(t),\hat\d(t,\omega)),\Delta \bar\d(t,\omega)}\notag\\
	\leq&c|Q^{-1}|_\infty\norm{\d(t,\d_0)}_2|\Delta Y(t)|_2|\Delta X(t)|_2+c|Q^{-1}|_\infty\norm{\d(t,\d_0+h b_{0})}_3|\nabla Y(t)|_2|\Delta X(t)|_2\notag\\
	&+\frac{c}{h}|Q^{-1}|_\infty|\Delta X(t)|_2\norm{\bar\d(t,\omega)}_2|\Delta \bar\d(t,\omega)|_2\notag\\
	\leq&\eps|\Delta X(t)|_2^2+c|Q^{-1}|^2_\infty\norm{\d(t,\d_0)}^2_2|\Delta Y(t)|_2^2+c|Q^{-1}|^2_\infty\norm{\d(t,\d_0+h b_{0})}_3^2|\nabla Y(t)|_2^2\notag\\
	&+c|Q^{-1}|^2_\infty\norm{ \bar\d(t,\omega)}_2^2|\Delta \bar\d(t,\omega)|_2^2.
\end{align}
Applying Gronwall's inequality yields that
\begin{align}
	&\sup_{0\leq t\leq T}|\nabla X(t)|_2^2+\int_0^T|\Delta X(t)|_2^2dt\notag\\
	\leq&c\bigg[|Q|_\infty^2\sup_{0\leq t\leq T}|X(t)|_2^2\int_0^T\norm{\u(t,\v_0+h\u_0)+Z(t)}_2^2dt+|Q|_\infty^2\sup_{0\leq t\leq T}\norm{\hat\u(t,\v_0)(\u_0)}_1^2|\nabla\bar\u(t,\omega)|_2^2T\notag\\
	&\hspace{5mm}+|Q^{-1}|^2_\infty\sup_{0\leq t\leq T}\norm{\d(t,\d_0)}^2_2|\Delta Y(t)|_2^2T+|Q^{-1}|^2_\infty\sup_{0\leq t\leq T}|\nabla Y(t)|_2^2\int_0^T\norm{\d(t,\d_0)}_3^2dt\notag\\
	&\hspace{5mm}+|Q^{-1}|^2_\infty\sup_{0\leq t\leq T}\norm{D(t,\omega)}_2^2|\Delta \bar\d(t,\omega)|_2^2T\bigg]\exp\{c|Q|_\infty^2\sup_{0\leq t\leq T}\norm{\u(t,\v_0)+Z(t)}_1^2T\}
\end{align}
According to \eqref{eqn:3.35}, \eqref{eqn:3.43} and \eqref{eqn:3.47}, we get
\begin{equation}\label{eqn:3.51}
		\lim_{h\to0}\sup_{\norm{\u_0}_1+\norm{\b_0}_2\leq1}\bigg\{\sup_{0\leq t\leq T}|\nabla X(t)|_2^2+\int_0^T|\Delta X(t)|_2^2dt\bigg\}=0.
\end{equation}
Combining \eqref{eqn:3.43}, \eqref{eqn:3.47}, and \eqref{eqn:3.51}, we conclude that
\begin{equation}
		\lim_{h\to0}\sup_{\norm{\u_0}_1+\norm{\b_0}_2\leq1}\bigg\{\sup_{0\leq t\leq T}[\norm{X(t)}_1^2+\norm{Y(t)}_2^2]+\int_0^T\norm{X(t)}_2^2dt+\int_0^T\norm{Y(t)}_3^2dtt\bigg\}=0.
\end{equation}
which proves \eqref{eqn:3.20}.

It remains to show $\V\times\dH^2\ni(\v_0,\d_0)\mapsto(\v(t,\v_0),\d(t,\d_0)\in\V\times\dH^2$ is Fr\'echet $C^{1,1}$, and to see that, it suffices to show $\V\times\dH^2\ni(\v_0,\d_0)\mapsto(\hat\u(t,\v_0),\hat\d(t,\d_0))\in L(\V\times\dH^2)$ is Lipschitz continuous on bounded sets.

Now let $\v_0,\v_0',\u_0\in\V, \d_0,\d_0',\b_0\in\dH^2$ with $\norm{\v_0}_1\leq M,\norm{\v_0'}_1\leq M, \norm{\d_0}_2\leq M,\norm{\d_0'}_2\leq M$ and $\norm{\u_0}_1+\norm{\b_0}_2\leq1$. By \eqref{eqn:3.6} and \eqref{eqn:3.7}, we have for $t\in[0,T]$,
\begin{align}\label{eqn:3.53}
	&\hat\u(t,\v_0)(\u_0)-\hat\u(t,\v_0')(\u_0)\notag\\
	=&-\int_0^tA_1(\hat\u(s,\v_0)(\u_0)-\hat\u(s,\v_0')(\u_0))ds\notag\\
	&-\int_0^tQ(s)B_1(\hat\u(s,\v_0)(\u_0)-\hat\u(s,\v_0')(\u_0),\u(s,\v_0)+Z(s))ds\notag\\
	&-\int_0^tQ(s)B_1(\u(s,\v_0)+Z(s),\hat\u(s,\v_0)(\u_0)-\hat\u(s,\v_0')(\u_0))ds\notag\\
	&-\int_0^tQ(s)B_1(\hat\u(s,\v_0')(\u_0),\u(s,\v_0)-\u(s,\v_0'))ds-\int_0^tQ(s)B_1(\u(s,\v_0)-\u(s,\v_0'),\hat\u(s,\v_0')(\u_0))ds\notag\\
	&-\int_0^tQ(s)^{-1}M(\hat\d(s,\d_0)(\b_0)-\hat\d(s,\d_0')(\b_0),\d(s,\d_0))ds\notag\\
	&-\int_0^tQ(s)^{-1}M(\hat\d(s,\d_0')(\b_0),\d(s,\d_0)-\d(s,\d_0'))ds\notag\\
	&-\int_0^tQ(s)^{-1}M(\d(s,\d_0),\hat\d(s,\d_0)(\b_0)-\hat\d(s,\d_0')(\b_0))ds\notag\\
	&-\int_0^tQ(s)^{-1}M(\d(s,\d_0)-\d(s,\d_0'),\hat\d(s,\d_0')(\b_0))ds;
	\end{align}
	\begin{align}
	\label{eqn:3.54}
	&\hat\d(t,\d_0)(\b_0)-\hat\d(t,\d_0')(\b_0)\notag\\
	=&-\int_0^tA_2(\hat\d(s,\d_0)(\b_0)-\hat\d(s,\d_0')(\b_0))ds-\int_0^tQ(s)B_2(\hat\u(s,\v_0)(\u_0)-\hat\u(s,\v_0')(\u_0),\d(s,\d_0))ds\notag\\
	&-\int_0^tQ(s)B_2(\u(s,\v_0)+Z(s),\hat\d(s,\d_0)(\b_0)-\hat\d(s,\d_0')(\b_0))ds\notag\\
	&-\int_0^tQ(s)B_2(\hat\u(s,\v_0')(\u_0),\d(s,\d_0)-\d(s,\d_0')ds-\int_0^tQ(s)B_2(\u(s,\v_0)-\u(s,\v_0'),\hat\d(s,\d_0')(\b_0))ds\notag\\
	&-\int_0^t f'(\d(s,\d_0))\hat\d(s,\d_0)(\b_0)ds+\int_0^t f'(\d(s,\d_0'))\hat\d(s,\d_0')(\b_0)ds.
\end{align}
For simplicity of notations, we denote by $\hat\u_\triangle(t):=\hat\u(t,\v_0)(\u_0)-\hat\u(t,\v_0')(\u_0)$ and $\hat\d_\triangle(t):=\hat\d(t,\d_0)(\b_0)-\hat\d(t,\d_0')(\b_0)$.

We first take inner product of \eqref{eqn:3.53} with $\hat\u_\triangle(t)$ in $L^{2}(\mathbf{D})$,
\begin{align}
	|\hat\u_\triangle(t)|_2^2=&-2\int_0^t|\nabla\hat\u_\triangle(s)|_2^2ds-2\int_0^tQ(s)\qv{B_1(\hat\u_\triangle(s),\u(s,\v_0)+Z(s)),\hat\u_\triangle(s)}ds\notag\\
	&-2\int_0^tQ(s)\qv{B_1(\u(s,\v_0)+Z(s),\hat\u_\triangle(s)),\hat\u_\triangle(s)}ds\notag\\
	&-2\int_0^tQ(s)\qv{B_1(\hat\u(s,\v_0')(\u_0),\u(s,\v_0)-\u(s,\v_0')),\hat\u_\triangle(s)}ds\notag\\
	&-2\int_0^tQ(s)\qv{B_1(\u(s,\v_0)-\u(s,\v_0'),\hat\u(s,\v_0')(\u_0)),\hat\u_\triangle(s)}ds\notag\\
	&-2\int_0^tQ(s)^{-1}\qv{M(\hat\d_\triangle(s),\d(s,\d_0)),\hat\u_\triangle(s)}ds\notag\\
	&-2\int_0^tQ(s)^{-1}\qv{M(\hat\d(s,\d_0')(\b_0),\d(s,\d_0)-\d(s,\d_0')),\hat\u_\triangle(s)}ds\notag\\
	&+2\int_0^tQ(s)^{-1}\qv{M(\d(s,\d_0),\hat\d_\triangle(s)),\hat\u_\triangle(s)}ds\notag\\
	&+2\int_0^tQ(s)^{-1}\qv{M(\d(s,\d_0)-\d(s,\d_0'),\hat\d(s,\d_0')(\b_0)),\hat\u_\triangle(s)}ds.
\end{align}
Applying Lemma \ref{lemma:2.1}, H\"older's inequality and Young's inequality gives that
\begin{align}\label{eqn:3.56}
	&|\hat\u_\triangle(t)|_2^2+\int_0^t|\nabla\hat\u_\triangle(s)|_2^2ds\notag\\
	\leq&c|Q|_\infty\int_0^t\norm{\u(s,\v_0)+Z(s)}_2|\hat\u_\triangle(s)|_2^2ds+c\int_0^t\norm{\u(s,\v_0)-\u(s,\v_0')}_1^2ds\notag\\
	&+c|Q|_\infty^2\int_0^t\norm{\hat\u(s,\v_0')(\u_0)}_1^2|\hat\u_\triangle(s)|_2^2ds+c|Q^{-1}|_\infty^2\int_0^t\norm{\d(s,\d_0)}_2^2|\nabla\hat\d_\triangle(s)|_2^2ds\notag\\
	&+c|Q^{-1}|_\infty^2\int_0^t\norm{\hat\d(s,\d_0')(\b_0)}_1^2\norm{\d(s,\d_0)-\d(s,\d_0')}_2^2ds.
\end{align}
Similarly to derive (\ref{eqn:3.11}), now we take innner product between $\partial_t\hat\d_\triangle(t)$ and $-\Delta\hat\d_\triangle(t)+\hat\d_\triangle(t)$ in $L^{2}(D)$ obtaining that ,
\begin{align}
		&\qv{\partial_t\hat\d_\triangle(t),-\Delta\hat\d_\triangle(t)+\hat\d_\triangle(t)}=\frac{1}{2}\frac{d|\nabla\hat\d_\triangle(t)|_2^2}{dt}+\frac{1}{2}\frac{d|\hat\d_\triangle(t)|_2^2}{dt}\notag\\
		=&-|-\Delta\hat\d_\triangle(t)+\hat\d_\triangle(t)|_2^2-\qv{-\Delta\hat\d_\triangle(t)+\hat\d_\triangle(t),Q(t)B_2(\hat\u_\triangle(t),\d(t,\d_0))}\notag\\
		&-\qv{-\Delta\hat\d_\triangle(t)+\hat\d_\triangle(t),Q(t)B_2(\u(t,\v_0)+Z(t),\hat\d_\triangle(t))}\notag\\
		&-\qv{-\Delta\hat\d_\triangle(t)+\hat\d_\triangle(t),Q(t)B_2(\hat\u(t,\v_0')(\u_0),\d(t,\d_0)-\d(t,\d_0'))}\notag\\
		&-\qv{-\Delta\hat\d_\triangle(t)+\hat\d_\triangle(t),Q(t)B_2(\u(t,\v_0)-\u(t,\v_0'),\hat\d(t,\d_0')(\b_0))}\notag\\
		&-\qv{-\Delta\hat\d_\triangle(t)+\hat\d_\triangle(t),f'(\d(t,\d_0))\hat\d_\triangle(t)-\hat\d_\triangle(t)}\notag\\
		&-\qv{-\Delta\hat\d_\triangle(t)+\hat\d_\triangle(t),(f'(\d(t,\d_0))-f'(\d(t,\d_0')))\hat\d(t,\d_0')(\b_0)}
\end{align}
With similar discussion as above, we obtain that
\begin{align}\label{eqn:3.58}
	&\frac{d|\nabla\hat\d_\triangle(t)|_2^2}{dt}+\frac{d|\hat\d_\triangle(t)|_2^2}{dt}+|-\Delta\hat\d_\triangle(t)+\hat\d_\triangle(t)|_2^2\notag\\
	\leq&c|Q|_\infty^2\norm{\d(t,\d_0)}_2^2|\hat\u_\triangle(t)|_2^2+c|Q|_\infty^2\norm{\u(t,\v_0)+Z(t)}_1^2|\nabla\hat\d_\triangle(t)|_2^2\notag\\
	&+c|Q|_\infty^2|\hat\u(t,\v_0')(\u_0)|_2^2\norm{\d(t,\d_0)-\d(t,\d_0')}_2^2+c|Q|_\infty^2\norm{\hat\d(t,\d_0')(\b_0)}_1^2\norm{\u(t,\v_0)-\u(t,\v_0')}_1^2\notag\\
	&+c(1+\norm{\d(t,\d_0)}_1^4)|\hat\d_\triangle(t)|_2^2+c\norm{\hat\d(t,\d_0')(\b_0)}_1^2(\norm{\d(t,\d_0)}_1^2+\norm{\d(t,\d_0')}_1^2)|\d(t,\d_0)-\d(t,\d_0')|_2^2.
\end{align}
Combining \eqref{eqn:3.56} and \eqref{eqn:3.58}, and applying Gronwall's inequality yields that
\begin{align*}
&\sup_{0\leq t\leq T}[|\hat\u_\triangle(t)|_2^2+|\nabla\hat\d_\triangle(t)|_2^2+|\hat\d_\triangle(t)|_2^2]+\int_0^T|\nabla\hat\u_\triangle(t)|_2^2dt+\int_0^T|-\Delta\hat\d_\triangle(t)+\hat\d_\triangle(t)|_2^2dt\notag\\
\leq&cT\sup_{0\leq t\leq T}\Big[\norm{\u(t,\v_0)-\u(t,\v_0')}_1^2+|Q^{-1}|_\infty^2\norm{\hat\d(t,\d_0')(\b_0)}_1^2\norm{\d(t,\d_0)-\d(t,\d_0')}_2^2\notag\\
&\hspace{15mm}+|Q|_\infty^2|\hat\u(t,\v_0')(\u_0)|_2^2\norm{\d(t,\d_0)-\d(t,\d_0')}_2^2+|Q|_\infty^2\norm{\hat\d(t,\d_0')(\b_0)}_1^2\norm{\u(t,\v_0)-\u(t,\v_0')}_1^2\notag\\
&\hspace{15mm}+\norm{\hat\d(t,\d_0')(\b_0)}_1^2(\norm{\d(t,\d_0)}_1^2+\norm{\d(t,\d_0')}_1^2)|\d(t,\d_0)-\d(t,\d_0')|_2^2\big]\notag\\
&\times\exp c\bigg\{|Q|_\infty\int_0^T\norm{\u(t,\v_0)+Z(t)}_2dt+|Q|^2_\infty\sup_{0\leq t\leq T}\norm{\hat\u(t,\v_0')(\u_0)}_1^2T\notag\\
&\hspace{15mm}+(|Q|_\infty+|Q^{-1}|_\infty)^2\sup_{0\leq t\leq T}\norm{\d(t,\d_0)}_2^2T+|Q|_\infty^2\sup_{0\leq t\leq T}\norm{\u(t,\v_0)+Z(t)}_1^2T\notag\\
&\hspace{15mm}+T+\sup_{0\leq t\leq T}\norm{\d(t,\d_0)}_1^4T\bigg\}.
\end{align*}
It was shown in \cite{GHZ} that
\begin{align}\label{eqn:3.59}
\sup_{0\leq t\leq T}[\norm{\u(t,\v_0)-\u(t,\v_0')}_1^2+\norm{\d(t,\d_0)-\d(t,\d_0')}_2^2]
\leq c(T)[\norm{\v_0-\v_0'}_1^2+\norm{\d_0-\d_0'}_2^2],
\end{align}
where $c(T)$ is bounded given the initial value norms are bounded by $M$. According to Proposition \ref{prop:3.1}, and \eqref{eqn:3.19}, we arrive at the following estimate:
\begin{align}\label{eqn:3.60}
	&\sup_{0\leq t\leq T}[|\hat\u_\triangle(t)|_2^2+|\nabla\hat\d_\triangle(t)|_2^2+|\hat\d_\triangle(t)|_2^2]+\int_0^T|\nabla\hat\u_\triangle(t)|_2^2dt+\int_0^T|-\Delta\hat\d_\triangle(t)+\hat\d_\triangle(t)|_2^2dt\notag\\
	\leq&g_5(T)[\norm{\v_0-\v_0'}_1^2+\norm{\d_0-\d_0'}_2^2].
\end{align}
Similarly to derive (\ref{eqn:3.14a}), taking inner product between $\Delta\hat\d_\triangle(t)$ and $\Delta\partial_t \hat\d_\triangle(t)$ in $L^{2}(\mathbf{D})$ yields,
\begin{align}
	&\qv{\Delta\hat\d_\triangle(t),\Delta\partial_t\hat\d_\triangle(t)}=\frac{1}{2}\partial_t|\Delta\hat\d_\triangle(t)|_2^2=-\qv{\nabla\Delta\hat\d_\triangle(t),\nabla\partial_t\hat\d(t)}\notag\\
	=&-|\nabla\Delta\hat\d_\triangle(t)|_2^2+\qv{\nabla\Delta\hat\d_\triangle(t),Q(t)B_2(\nabla\hat\u_\triangle(t),\d(t,\d_0))}\notag\\
	&+\qv{\nabla\Delta\hat\d_\triangle(t),Q(t)B_2(\hat\u_\triangle(t),\nabla\d(t,\d_0))}+\qv{\nabla\Delta\hat\d_\triangle(t),Q(t)B_2(\nabla(\u(t,\v_0)+Z(t)),\hat\d_\triangle(t))}\notag\\
	&+\qv{\nabla\Delta\hat\d_\triangle(t),Q(t)B_2(\u(t,\v_0)+Z(t),\nabla\hat\d_\triangle(t))}\notag\\
	&+\qv{\nabla\Delta\hat\d_\triangle(t),Q(t)B_2(\nabla\hat\u(t,\v_0')(\u_0),\d(t,\d_0)-\d(t,\d_0'))}\notag\\
	&+\qv{\nabla\Delta\hat\d_\triangle(t),Q(t)B_2(\hat\u(t,\v_0')(\u_0),\nabla(\d(t,\d_0)-\d(t,\d_0'))}\notag\\
	&+\qv{\nabla\Delta\hat\d_\triangle(t),Q(t)B_2(\nabla(\u(t,\v_0)-\u(t,\v_0')),\hat\d(t,\d_0')(\b_0))}\notag\\
	&+\qv{\nabla\Delta\hat\d_\triangle(t),Q(t)B_2(\u(t,\v_0)-\u(t,\v_0'),\nabla\hat\d(t,\d_0')(\b_0))}\notag\\
	&+\qv{\nabla\Delta\hat\d_\triangle(t),\nabla[f'(\d(t,\d_0))\hat\d_\triangle(t)]}+\qv{\nabla\Delta\hat\d_\triangle(t),\nabla[(f'(\d(t,\d_0))-f'(\d(t,\d_0')))\hat\d(t,\d_0')(\b_0)]}.
\end{align}
Applying H\"older's inequality and Young's inequality gives that
\begin{align}
	&\frac{d|\Delta\hat\d_\triangle(t)|_2^2}{dt}+|\nabla\Delta\hat\d_\triangle(t)|_2^2\notag\\
	\leq&c|Q|_\infty^2\norm{\d(t,\d_0)}_2^2|\nabla\hat\u_\triangle(t)|_2^2+c|Q|_\infty^2\norm{\d(t,\d_0)}_3^2|\hat\u_\triangle(t)|_2^2\notag\\
	&+c|Q|_\infty^2\norm{\u(t,\v_0)+Z(t)}_2^2|\nabla\hat\d_\triangle(t)|_2^2++c|Q|_\infty^2\norm{\u(t,\v_0)+Z(t)}_1^2|\Delta\hat\d_\triangle(t)|_2^2\notag\\
	&+c|Q|_\infty^2\norm{\hat\u(t,\v_0')(\u_0)}_1^2\norm{\d(t,\d_0)-\d(t,\d_0')}_2^2+c\norm{\hat\d(t,\d_0')(\b_0)}_2^2\norm{\u(t,\v_0)-\u(t,\v_0')}_1^2\notag\\
	&+c\norm{\d(t,\d_0)}_1^2\norm{\d(t,\d_0)}_2^2|\hat\d_\triangle(t)|_2^2+c(1+\norm{\d(t,\d_0)}_1^4)|\nabla\hat\d_\triangle(t)|_2^2\notag\\
	&+c(\norm{\d(t,\d_0)}_1^2
	+\norm{\d(t,\d_0')}_1^2)\norm{\hat\d(t,\d_0')(\b_0)}_2^2|\d(t,\d_0)-\d(t,\d_0')|_2^2\notag\\
	&+c(\norm{\d(t,\d_0)}_2^2
	+\norm{\d(t,\d_0')}_2^2)\norm{\hat\d(t,\d_0')(\b_0)}_1^2|\d(t,\d_0)-\d(t,\d_0')|_2^2\notag\\
	&+c(\norm{\d(t,\d_0)}_1^2
	+\norm{\d(t,\d_0')}_1^2)\norm{\hat\d(t,\d_0')(\b_0)}_1^2\norm{\d(t,\d_0)-\d(t,\d_0')}_1^2.
\end{align}
Applying Gronwall's inequality, and with the estimates \eqref{eqn:3.59}, \eqref{eqn:3.60}, there exists $g_6(T)<\infty$ such that
\begin{equation}\label{eqn:3.63}
	\sup_{0\leq t\leq T}|\Delta\hat\d_\triangle(t)|_2^2+\int_0^T|\nabla\Delta\hat\d_\triangle(t)|_2^2dt\leq g_5(T)[\norm{\v_0-\v_0'}_1^2+\norm{\d_0-\d_0'}_2^2].
\end{equation}
We now take inner product between $\hat\u_\triangle(t)$ and $-\Delta\hat\u_\triangle(t)$ in $L^{2}(\mathbf{D})$,
\begin{align}
	\frac{d|\nabla\hat\u_\triangle(t)|_2^2}{dt}=&-2|\Delta\hat\u_\triangle(t)|_2^2\notag\\
	&+2Q(t)\qv{B_1(\hat\u_\triangle(t),\u(t,\v_0)+Z(t)),\Delta\hat\u_\triangle(t)}\notag\\
	&+2Q(t)\qv{B_1(\u(t,\v_0)+Z(t),\hat\u_\triangle(t)),\Delta\hat\u_\triangle(t)}\notag\\
	&+2Q(t)\qv{B_1(\hat\u(t,\v_0')(\u_0),\u(t,\v_0)-\u(t,\v_0')),\Delta\hat\u_\triangle(t)}\notag\\
	&+2Q(t)\qv{B_1(\u(t,\v_0)-\u(t,\v_0'),\hat\u(t,\v_0')(\u_0)),\Delta\hat\u_\triangle(t)}\notag\\
	&+2Q(t)^{-1}\qv{M(\hat\d_\triangle(t),\d(t,\d_0)),\Delta\hat\u_\triangle(t)}\notag\\
	&+2Q(t)^{-1}\qv{M(\hat\d(t,\d_0')(\b_0),\d(t,\d_0)-\d(t,\d_0')),\Delta\hat\u_\triangle(t)}\notag\\
	&+2Q(t)^{-1}\qv{M(\d(t,\d_0),\hat\d_\triangle(t)),\Delta\hat\u_\triangle(t)}\notag\\
	&+2Q(s)^{-1}\qv{M(\d(t,\d_0)-\d(t,\d_0'),\hat\d(t,\d_0')(\b_0)),\Delta\hat\u_\triangle(t)}\notag\\
	=:&n_1+\cdots+n_9.
\end{align}
By Lemma \ref{lemma:2.1}, H\"older's inequality and Young's inequality,
\begin{align*}
	n_2+n_3\leq&c|Q|_\infty|\hat\u_\triangle(t)|_2|\nabla(\u(t,\v_0)+Z(t))|_\infty|\Delta\hat\u_\triangle(t)|_2\\
	&+|Q|_\infty|\u(t,\v_0)+Z(t)|_\infty|\nabla\hat\u_\triangle(t)|_2|\Delta\hat\u_\triangle(t)|_2\\
	\leq&\eps|\Delta\hat\u_\triangle(t)|_2^2+c|Q|_\infty^2\norm{\u(t,\v_0)+Z(t)}_2^2|\hat\u_\triangle(t)|_2^2+c|Q|_\infty^2\norm{\u(t,\v_0)+Z(t)}_1^2|\nabla\hat\u_\triangle(t)|_2^2,
\end{align*}
similarly,
\begin{align*}
	n_4+n_5\leq&c|Q|_\infty|\hat\u(t,\v_0')|_\infty|\nabla(\u(t,\v_0)-\u(t,\v_0'))|_2|\Delta\hat\u_\triangle(t)|_2\\
	&+|Q|_\infty|\u(t,\v_0)-\u(t,\v_0')|_\infty|\nabla\hat\u(t,\v_0')|_2|\Delta\hat\u_\triangle(t)|_2\\
	\leq&\eps|\Delta\hat\u_\triangle(t)|_2^2+c|Q|_\infty^2\norm{\hat\u(t,\v_0')}_1^2\norm{\u(t,\v_0)-\u(t,\v_0')}_1^2.
\end{align*}
By Proposition \ref{prop:2.5}, H\"older's inequality and Young's inequality,
\begin{align*}
	n_6+n_8
	=&2Q(t)^{-1}\qv{B_2(\Delta\hat\u_\triangle(t),\hat\d_\triangle(t)),\Delta\d(t,\d_0)}+2Q(t)^{-1}\qv{B_2(\Delta\hat\u_\triangle(t),\d(t,\d_0)),\Delta\hat\d_\triangle(t)}ds\\
	\leq&c|Q^{-1}|_\infty[|\Delta\hat\u_\triangle(t)|_2|\nabla\hat\d_\triangle(t)|_2|\Delta\d(t,\d_0)|_\infty+|\Delta\hat\u_\triangle(t)|_2|\nabla\d(t,\d_0)|_\infty|\Delta\hat\d_\triangle(t)|_2]\\
	\leq&\eps|\Delta\hat\u_\triangle(t)|_2^2+c|Q^{-1}|_\infty^2\norm{\d(t,\d_0)}_3^2|\nabla\hat\d_\triangle(t)|_2^2+c|Q^{-1}|_\infty^2\norm{\d(t,\d_0)}_2^2|\Delta\hat\d_\triangle(t)|_2^2.
\end{align*}
Similarly,
\begin{align*}
	n_7+n_9=&2Q(t)^{-1}\qv{B_2(\Delta\hat\u_\triangle(t),\hat\d(t,\d_0')(\b_0)),\Delta(\d(t,\d_0)-\d(t,\d_0'))}\\
	&+2Q(t)^{-1}\qv{B_2(\Delta\hat\u_\triangle(t),\d(t,\d_0)-\d(t,\d_0')),\Delta(\hat\d(t,\d_0')(\b_0))}\\
	\leq&c|Q^{-1}|_\infty|\Delta\hat\u_\triangle(t)|_2|\nabla\hat\d(t,\d_0')(\b_0)|_\infty|\Delta(\d(t,\d_0)-\d(t,\d_0'))|_2\\
	&+c|Q^{-1}|_\infty|\Delta\hat\u_\triangle(t)|_2|\nabla(\d(t,\d_0)-\d(t,\d_0'))|_\infty|\Delta\hat\d(t,\d_0')(\b_0)|_2\\
	\leq&\eps|\Delta\hat\u_\triangle(t)|_2^2+c|Q^{-1}|_\infty^2\norm{\hat\d(t,\d_0')(\b_0)}_2^2\norm{\d(t,\d_0)-\d(t,\d_0')}_2^2.
	\end{align*}
	Altogether, applying Gronwall's inequality, one arrives at
\begin{align}
\label{eqn:3.65}
&\sup_{0\leq t\leq T}|\nabla\hat\u_\triangle(t)|_2^2+\int_0^T|\Delta\hat\u_\triangle(t)|_2^2dt
\notag\\
\leq&c\bigg[|Q|_\infty^2\sup_{0\leq t\leq T}|\hat\u_\triangle(t)|_2^2\int_0^T\norm{\u(t,\v_0)+Z(t)}_2^2dt+|Q|_\infty^2\sup_{0\leq t\leq T}\norm{\hat\u(t,\v_0')}_1^2\norm{\u(t,\v_0)-\u(t,\v_0')}_1^2T\notag\\
&\hspace{5mm}+|Q^{-1}|_\infty^2\sup_{0\leq t\leq T}|\nabla\hat\d_\triangle(t)|_2^2\int_0^T\norm{\d(t,\d_0)}_3^2dt+|Q^{-1}|_\infty^2\sup_{0\leq t\leq T}\norm{\d(t,\d_0)}_2^2|\Delta\hat\d_\triangle(t)|_2^2T\notag\\
&\hspace{5mm}+|Q^{-1}|_\infty^2\sup_{0\leq t\leq T}\norm{\hat\d(t,\d_0')(\b_0)}_2^2\norm{\d(t,\d_0)-\d(t,\d_0')}_2^2T\bigg]\notag\\
&\times\exp \{c|Q|_\infty^2\sup_{0\leq t\leq T}\norm{\u(t,\v_0)+Z(t)}_1^2T\}.
\end{align}
Thus, with \eqref{eqn:3.60}, \eqref{eqn:3.63} and \eqref{eqn:3.65}, we conclude that
\begin{align}
&\sup_{0\leq t\leq T}\norm{\hat\u(t,\v_0)(\u_0)-\hat\u(t,\v_0')(\u_0)}_{1}^2+\int_0^T\norm{\hat\u(t,\v_0)(\u_0)-\hat\u(t,\v_0')(\u_0)}_{L^2([0,T];\vH^2)}^2dt\notag\\
&+\sup_{0\leq t\leq T}\norm{\hat\d(t,\d_0)(\b_0)-\hat\d(t,\d_0')(\b_0)}_{2}^2+\int_0^T\norm{\hat\d(t,\d_0)(\b_0)-\hat\d(t,\d_0')(\b_0)}_{L^2([0,T];\dH^3)}^2dt\notag\\
\leq&c[\norm{\v_0-\v_0'}_1^2+\norm{\d_0-\d_0'}_2^2],
\end{align}
where $c$ is a positive constant that is independent of initial date provided that  $\norm{\v_0}_1\leq M,\norm{\v_0'}_1\leq M, \norm{\d_0}_2\leq M,\norm{\d_0'}_2\leq M$. Hence, we prove that the map $\V\times\dH^2\ni(\v_0,\d_0)\mapsto(\hat\u(t,\v_0),\hat\d(t,\d_0)\in L(\V\times\dH^2)$ is Lipschitz continuous on bounded sets.

To see the compactness of the Fr\'echet derivative, we can follow the method in \cite{GHZ} and use the Aubin-Lions Lemma as well as the regularity of solutions. One can also adopt the method in Theorem 3.1 of \cite{MS} to show the compactness of $(D\v(t,\v_0,\omega),D\d(t,\d_0,\omega)):\V\times\dH^2\to\V\times\dH^2$ for $t>0$.
\qed

Now we are ready to discuss the Malliavin regularity for solutions of the stochastic nematic liquid crystal equations.
\begin{thm}
	\label{prop:3.4}
	For $\v_0\in\V,\d_0\in\dH^2$ and $t\geq0$, the solution maps $\omega\mapsto\u(t,\v_0,\omega),\omega\mapsto\d(t,\d_0,\omega)$ are Malliavin differentiable, and for all $t\in[0,T]$, almost surely their Malliavin derivatives $\D_v\u(t,\v_0),\D_v\d(t,\d_0)$ solve the following equations:
	\begin{align}
	\D_v\u(t,\v_0)=&-\int_0^tA_1\D_v\u(s,\v_0)ds\notag\\
	&-\int_0^tQ(s)B_1\left(\D_v\u(s,\v_0)+\D_vZ(s),\u(s,\v_0)+Z(s)\right)ds\notag\\
	&-\int_0^tQ(s)B_1\left(\u(s,\v_0)+Z(s),\D_v\u(s,\v_0)+\D_vZ(s)\right)ds\notag\\
	&-\int_0^t\D_vQ(s)B_1(\u(s,\v_0)+Z(s))ds\notag\\
	&-\int_0^tQ(s)^{-1}M(\D_v\d(s,\d_0),\d(s,\d_0))ds-\int_0^tQ(s)^{-1}M(\d(s,\d_0),\D_v\d(s,\d_0))ds\notag\\
	&-\int_0^t\D_vQ(s)^{-1}M(\d(s,\d_0))ds.\\
	\D_v\d(t,\d_0)=&-\int_0^tA_2\D_v\d(s,\d_0)ds\notag\\
	&-\int_0^tQ(s)B_2(\D_v\u(s,\v_0)+\D_vZ(s), \d(s,\d_0))ds\notag\\
	&-\int_0^tQ(s)B_2\left(\u(s,\v_0)+Z(s),\D_v\d(s,\d_0)\right)ds\notag\\
	&-\int_0^t\D_vQ(s)B_2(\u(s,\v_0)+Z(s),\d(s,\d_0))ds-\int_0^tf'(\d(s,\d_0))\D_v\d(s,\d_0)ds.
	\end{align}
\end{thm}
\proof We will show $\v(t,\v_0)\in\D_{\text{loc}}^{1,2}(\vH), \d(t,\d_0)\in\D_{\text{loc}}^{1,2}(\dH^1)$. First by the  uniqueness of solutions to the model \eqref{eqn:3.3}, we define $\u(t,\v_0)=\u_N(t,\v_0), \d(t,\d_0)=\d_N(t,\d_{0})$ as the solution on $\Omega_N=\{\sup_{0\leq t\leq T}(|W(t)|\vee\norm{Z(t)}_2)\leq N\}$, that is, $\u_N(t,\v_0), \d_N(t,\d_{0,n})$ are solutions to  \eqref{eqn:3.3} with $Q(t), Z(t)$ replaced by $Q_N(t):=\exp\{W(t)1_{\{|W|\leq N\}}\}$ and $Z_N(t):=Z(t)1_{\{\norm{Z}_2\leq N\}}$. For simplicity, we still use $Q, Z$ to represent $Q_N, Z_N$.

Now we use Galerkin approximation and write $\{e_k\}_{k\geq1}$ as an orthonormal basis for $\V$, serving as eigenvectors of $-A_1$ subject to the boundary condition \eqref{BC}, with corresponding eigenvalues $\{r_k\}_{k\geq1}$, that is, $A_1e_k=-r_ke_k$. Let $\V_n$ be $n$-dimensional subspace spanned by $\{e_1,\dots,e_n\}$, and define
$$\v_{0,n}=\sum_{k=1}^{n}\qv{\v_0,e_k}e_k.$$
Similarly, let $\{\rho_k\}_{k\geq1}$ be an orthonormal basis for $\dH^2$, which serves as eigenvectors of $-A_2$ subject to the boundary condition. Let $\dH^2_n$ be $n$-dimensional subspace spanned by $\{\rho_1,\dots,\rho_n\}$ and define
$$\d_{0,n}=\sum_{k=1}^n\qv{\d_0,\rho_k}\rho_k.$$
Now we let $(\u_n(t,\v_{0,n}),\d_n(t,\d_{0,n})\in\V_n\times\dH^2_n$ be the unique solution to the following equations:
\begin{align}
&d\u_n(t,\v_{0,n})=-A_1\u_n(t,\v_{0,n})dt-Q(t)B_1(\u_n(t,\v_{0,n})+Z(t))dt-Q(t)^{-1}M(\d_n(t,\d_{0,n}))dt,\\
&\nabla\cdot(\u_n(t,\v_{0,n})+Z(t))=0,\\
&d\d_{n}(t,\d_{0,n})=-A_2\d_n(t,\d_{0,n})dt-Q(t)B_2(\u_n(t,\v_{0,n})+Z(t),\d_n(t,\d_{0,n}))dt-f(\d_n(t,\d_{0,n})))dt,\\
&\u_n(t,\v_{0,n})|_{\partial \mathbf{D}}=0,\ \d_n(t,\d_{0,n})|_{\partial\mathbf{D}}=\d_{0,n},\\
&\u_n(0,\v_{0,n})=\v_{0,n},\ \d_n(0,\d_{0,n})=\d_{0,n}.
\end{align}
By the proof of global well-posedness of stochastic nematic liquid crystals equations in ( \cite{BMP}), we know there exists a subsequence still denoted by $(\u_n,\d_n)$ such that
\begin{equation*}
\lim_{n\to\infty}\int_{0}^{T}\left(\norm{\u_n(t,\v_{0,n})-\u(t,\v_0)}_1^2+\norm{\d_n(t,\d_{0,n})-\d(t,\d_0)}_2^2\right)dt=0,\ \ a.s..
\end{equation*}
Then by the localization and the dominated convergence theorem we get that
\begin{equation}
\label{eqn:3.74}
\mathbb{E}\lim_{n\to\infty}\int_{0}^{T}\left(\norm{\u_n(t,\v_{0,n})-\u(t,\v_0)}_1^2+\norm{\d_n(t,\d_{0,n})-\d(t,\d_0)}_2^2\right)dt=0.
\end{equation}
For any $t\in[0,T]$, the Malliavin derivatives $\D_v\u_n(t,\v_{0,n})$, $\D_v\d_n(t,\d_{0,n})$ satisfy
\begin{align}
\D_v\u_n(t,\v_{0,n})=&-\int_0^tA_1\D_v\u_n(s,\v_{0,n})ds\notag\\
&-\int_0^tQ(s)B_1\left(\D_v\u_n(s,\v_{0,n})+\D_vZ(s),\u_n(s,\v_{0,n})+Z(s)\right)ds\notag\\
&-\int_0^tQ(s)B_1\left(\u_n(s,\v_{0,n})+Z(s),\D_v\u_n(s,\v_{0,n})+\D_vZ(s)\right)ds\notag\\
&-\int_0^t\D_vQ(s)B_1(\u_n(s,\v_{0,n})+Z(s))ds\notag\\
&-\int_0^tQ(s)^{-1}M(\D_v\d_n(s,\d_{0,n}),\d_n(s,\d_{0,n}))ds-\int_0^tQ(s)^{-1}M(\d_n(s,\d_{0,n}),\D_v\d_n(s,\d_{0,n}))ds\notag\\
&-\int_0^t\D_vQ(s)^{-1}M(\d_n(s,\d_{0,n}))ds.\\
\D_v\d_n(t,\d_{0,n})=&-\int_0^tA_2\D_v\d_n(s,\d_{0,n})ds\notag\\
&-\int_0^tQ(s)B_2(\D_v\u_n(s,\v_{0,n})+\D_vZ(s), \d_n(s,\d_{0,n}))ds\notag\\
&-\int_0^tQ(s)B_2\left(\u_n(s,\v_{0,n})+Z(s),\D_v\d_n(s,\d_{0,n})\right)ds\notag\\
&-\int_0^t\D_vQ(s)B_2(\u_n(s,\v_{0,n})+Z(s),\d_n(s,\d_{0,n}))ds\notag\\
&-\int_0^t f'(\d_n(s,\d_{0,n}))\D_v\d_n(s,\d_{0,n})ds.
\end{align}
Now let $\xi_v,\eta_v$ be the solution to the following random equations as well as the boundary conditions \eqref{BC}
\begin{align}
\xi_v(t,\v_0)=&-\int_0^tA_1\xi_v(s,\v_0)ds\notag\\
&-\int_0^tQ(s)B_1\left(\xi_v(s,\v_0)+\D_vZ(s),\u(s,\v_0)+Z(s)\right)ds\notag\\
&-\int_0^tQ(s)B_1\left(\u(s,\v_0)+Z(s),\xi_v(s,\v_0)+\D_vZ(s)\right)ds\notag\\
&-\int_0^t\D_vQ(s)B_1(\u(s,\v_0)+Z(s))ds\notag\\
&-\int_0^tQ(s)^{-1}M(\eta_v(s,\d_0),\d(s,\d_0))ds-\int_0^tQ(s)^{-1}M(\d(s,\d_0),\eta_v(s,\d_0))ds\notag\\
&-\int_0^t\D_vQ(s)^{-1}M(\d(s,\d_0))ds.\\
\eta_v(t,\d_0)=&-\int_0^tA_2\eta_v(s,\d_0)ds\notag\\
&-\int_0^tQ(s)B_2(\xi_v(s,\v_0)+\D_vZ(s), \d(s,\d_0))ds\notag\\
&-\int_0^tQ(s)B_2\left(\u(s,\v_0)+Z(s),\eta_v(s,\d_0)\right)ds\notag\\
&-\int_0^t\D_vQ(s)B_2(\u(s,\v_0)+Z(s),\d(s,\d_0))ds-\int_0^tf'(\d(s,\d_0))\eta_v(s,\d_0)ds,
\end{align}
for any $t\in[0,T]$. The global well-posedness of the above equations have been studied in \cite{BM}. Since $\D$ is closed, it suffices to show that
\begin{equation}
\label{eqn:3.79}
\lim_{n\to\infty}\Ex\sup_{0\leq v\leq t}\left[|\D_v\u_n(t,\v_{0,n})-\xi_v(t,\v_0)|_2^2+\norm{\D_v\d_n(t,\d_{0,n})-\eta_v\d(t,\d_0)}_1^2\right]=0
\end{equation}
Define the following norm notations:
\begin{align}
\label{eqn:norm}
&C_1^n:=\sup_{0\leq t\leq T}\left(|\u_n(t,\v_{0,n})|_2+|Z(t)|_2\right),&&C_2^n:=\sup_{0\leq t\leq T}\left(\norm{\u_n(t,\v_{0,n})}_1+\norm{Z(t)}_1\right),\notag\\
&C_1:=\sup_{0\leq t\leq T}\left(|\u(t,\v_0)|_2+|Z(t)|_2\right),&&C_2:=\sup_{0\leq t\leq T}\left(\norm{\u(t,\v_0)}_1+\norm{Z(t)}_1\right),\notag\\
&M_1^n:=\sup_{0\leq t\leq T}\left(|\D_v\u_n(t,\v_{0,n})|_2+|\D_vZ(t)|_2\right),&&M_2^n:=\sup_{0\leq t\leq T}\left(\norm{\D_v\u_n(t,\v_{0,n})}_1+\norm{\D_vZ(t)}_1\right),\notag\\
&M_1:=\sup_{0\leq t\leq T}|\xi_v(t,\v_0)|_2,&&M_2:=\sup_{0\leq t\leq T}\norm{\xi_v(t,\v_0)}_1;\notag\\
&D_1^n:=\sup_{0\leq t\leq T}|\d_n(t,\d_{0,n})|_2, \ D_2^n:=\sup_{0\leq t\leq T}\norm{\d_n(t,\d_{0,n})}_1,\ &&D_3^n:=\sup_{0\leq t\leq T}\norm{\d_n(t,\d_{0,n})}_2,\notag\\
&D_1:=\sup_{0\leq t\leq T}|\d(t,\d_0)|_2,\ D_2:=\sup_{0\leq t\leq T}\norm{\d(t,\d_0)}_1,\ &&D_3:=\sup_{0\leq t\leq T}\norm{\d(t,\d_0)}_2,\notag\\
&N_1^n:=\sup_{0\leq t\leq T}|\D_v\d_n(t,\d_{0,n})|_2,\ N_2^n:=\sup_{0\leq t\leq T}\norm{\D_v\d_n(t,\d_{0,n})}_1,\ &&N_3^n:=\sup_{0\leq t\leq T}\norm{\D_v\d_n(t,\d_{0,n})}_2,\notag\\
&N_1:=\sup_{0\leq t\leq T}|\eta_v(t,\d_0)|_2,\ N_2:=\sup_{0\leq t\leq T}\norm{\eta_v(t,\d_0)}_1,\ &&N_3:=\sup_{0\leq t\leq T}\norm{\eta_v(t,\d_0)}_2.
\end{align}
We first estimate the followings:
\begin{align}
\label{eqn:uL2}
&|\D_v\u_n(t,\v_{0,n})-\xi_v(t,\v_0)|_2^2\notag\\
=&-2\int_0^t|\nabla(\D_v\u_n(s,\v_{0,n})-\xi_v(s,\v_0))|_2^2ds\notag\\
&-2\int_0^tQ(s)\qv{B_1\left(\D_v\u_n(s,\v_{0,n})+\D_vZ(s),\u_n(s,\v_{0,n})-\u(s,\v_0)\right), \D_v\u_n(s,\v_{0,n})-\xi_v(s,\v_0)}ds\notag\\
&-2\int_0^tQ(s)\qv{B_1\left(\D_v\u_n(s,\v_{0,n})-\xi_v(s,\v_0),\u(s,\v_0)+Z(s)\right),\D_v\u_n(s,\v_{0,n})-\xi_v(s,\v_0)}ds\notag\\
&-2\int_0^tQ(s)\qv{B_1\left(\u_n(s,\v_{0,n})-\u(s,\v_0),\D_v\u_n(s,\v_{0,n})+\D_vZ(s)\right),\D_v\u_n(s,\v_{0,n})-\xi_v(s,\v_0)}ds\notag\\
&-2\int_0^tQ(s)\qv{B_1\left(\u(s,\v_0)+Z(s),\D_v\u_n(s,\v_{0,n})-\xi_v(s,\v_0)\right),\D_v\u_n(s,\v_{0,n})-\xi_v(s,\v_0)}ds\notag\\
&-2\int_0^t\D_vQ(s)\qv{B_1(\u_n(s,\v_{0,n})+Z(s),\u_n(s,\v_{0,n})-\u(s,\v_0)), \D_v\u_n(s,\v_{0,n})-\xi_v(s,\v_0)}ds\notag\\
&-2\int_0^t\D_vQ(s)\qv{B_1(\u_n(s,\v_{0,n})-\u(s,\v_0),\u(s,\v_0)+Z(s)) ,\D_v\u_n(s,\v_{0,n})-\xi_v(s,\v_0)}\notag\\
&-2\int_0^tQ(s)^{-1}\qv{M(\D_v\d_n(s,\d_{0,n}),\d_n(s,\d_{0,n})-\d(s,\d_0)),\D_v\u_n(s,\v_{0,n})-\xi_v(s,\v_0)}ds\notag\\
&-2\int_0^tQ(s)^{-1}\qv{M(\D_v\d_n(s,\d_{0,n})-\eta_v(s,\d_0),\d(s,\d_0)),\D_v\u_n(s,\v_{0,n})-\xi_v(s,\v_0)}ds\notag\\
&-2\int_0^tQ(s)^{-1}\qv{M(\d_n(s,\d_{0,n})-\d(s,\d_0),\D_v\d_n(s,\d_{0,n})),\D_v\u_n(s,\v_{0,n})-\xi_v(s,\v_0)}ds\notag\\
&-2\int_0^tQ(s)^{-1}\qv{M(\d(s,\d_0),\D_v\d_n(s,\d_{0,n})-\eta_v(s,\d_0)),\D_v\u_n(s,\v_{0,n})-\xi_v(s,\v_0)}ds\notag\\
&-2\int_0^t\D_vQ(s)^{-1}\qv{M(\d_n(s,\d_{0,n}), \d_n(s,\d_{0,n})-\d(s,\d_0)),\D_v\u_n(s,\v_{0,n})-\xi_v(s,\v_0)}ds\notag\\
&-2\int_0^t\D_vQ(s)^{-1}\qv{M(\d_n(s,\d_{0,n})-\d(s,\d_0), \d(s,\d_0)),\D_v\u_n(s,\v_{0,n})-\xi_v(s,\v_0)}ds\notag\\
=:&I_1+\dots+I_{13}.
\end{align}
By H\"older's inequality and Young's inequality, we get that
\begin{align*}
I_2\leq&c|Q|_\infty\int_0^t|\D_v\u_n(s,\v_{0,n})+\D_vZ(s)|_\infty|\nabla(\u_n(s,\v_{0,n})-\u(s,\v_0))|_2|\D_v\u_n(t,\v_{0,n})-\xi_v(t,\v_0)|_2ds\\
\leq&c\int_0^t\norm{\u_n(s,\v_{0,n})-\u(s,\v_0)}_1^2ds+c|Q|_\infty^2[M_2^n]^2\int_0^t|\D_v\u_n(t,\v_{0,n})-\xi_v(t,\v_0)|_2^2ds,
\end{align*}
and
\begin{align*}
I_3\leq&c|Q|_\infty\int_0^t|\u(s,\v_0)+Z(s)|_\infty|\D_v\u_n(s,\v_{0,n})-\xi_v(s,\v_0)|_2^2ds\notag\\
\leq&c|Q|_\infty C_2\int_0^t|\D_v\u_n(s,\v_{0,n})-\xi_v(s,\v_0)|_2^2ds,
\end{align*}
similarly, we have
\begin{align*}
I_4\leq&c|Q|_\infty\int_0^t|\u_n(s,\v_{0,n})-\u(s,\v_0)|_\infty|\nabla(\D_v\u_n(s,\v_{0,n})+\D_vZ(s))|_2|\D_v\u_n(s,\v_{0,n})-\xi_v(s,\v_0)|_2ds\\
\leq&c\int_0^t\norm{\u_n(s,\v_{0,n})-\u(s,\v_0)}_1^2ds+c|Q|_\infty^2[M_2^n]^2\int_0^t|\D_v\u_n(t,\v_{0,n})-\xi_v(t,\v_0)|_2^2ds.
\end{align*}
According to Lemma \ref{lemma:2.1}, $\qv{B_1(\u,\v),\v}=0$, we get $I_5=0$, and
\begin{align*}
I_6\leq&c\sup_{0\leq t\leq T}|\D_vQ|\int_0^t|\u_n(s,\v_{0,n})+Z(s)|_\infty|\nabla(\u_n(s,\v_{0,n})-\u(s,\v_0))|_2|\D_v\u_n(s,\v_{0,n})-\xi_v(s,\v_0)|_2\\
\leq&c\int_0^t\norm{\u_n(s,\v_{0,n})-\u(s,\v_0)}_1^2ds+c|\D_vQ|_\infty^2[C_2^n]^2\int_0^t|\D_v\u_n(t,\v_{0,n})-\xi_v(t,\v_0)|_2^2ds,
\end{align*}
and similarly we obtain that
\begin{align*}
I_7\leq&c\int_0^t\norm{\u_n(s,\v_{0,n})-\u(s,\v_0)}_1^2ds+c|\D_vQ|_\infty^2[C_2]^2\int_0^t|\D_v\u_n(t,\v_{0,n})-\xi_v(t,\v_0)|_2^2ds.
\end{align*}
By Lemma \ref{lemma:2.3}, we obtain that
\begin{align*}
I_8+I_{10}\leq&c|Q^{-1}|_\infty\int_0^t\norm{\D_v\d_n(s,\d_{0,n})}_1\norm{\d_n(s,\d_{0,n})-\d(s,\d_0)}_2|\nabla(\D_v\u_n(s,\v_{0,n})-\xi_v(s,\v_0))|_2ds\\
\leq&\eps\int_0^t|\nabla(\D_v\u_n(s,\v_{0,n})-\xi_v(s,\v_0))|_2^2ds+c|Q^{-1}|_\infty^2[N_2^n]^2\int_0^t\norm{\d_n(s,\d_{0,n})-\d(s,\d_0)}_2^2ds.
\end{align*}
By H\"older's inequality and Young's inequality,
\begin{align*}
I_9+I_{11}\leq&c|Q^{-1}|\int_0^t|\nabla(\D_v\d_n(s,\d_{0,n})-\eta_v(s,\d_0))|_2\norm{\d(s,\d_0)}_2|\nabla(\D_v\u_n(s,\v_{0,n})-\xi_v(s,\v_0))|_2ds\\
\leq&\eps\int_0^t|\nabla(\D_v\u_n(s,\v_{0,n})-\xi_v(s,\v_0))|_2^2ds+c|Q^{-1}|_\infty^2[D_3]^2\int_0^t|\nabla(\D_v\d_n(s,\d_{0,n})-\eta_v(s,\d_0))|_2^2ds.
\end{align*}
Similarly, we get
\begin{align*}
I_{12}+I_{13}\leq&c|\D_vQ(s)^{-1}|_\infty\int_0^t[\norm{\d_n(s,\d_{0,n})}_1+\norm{\d(s,\d_0)}_1]\\
&\hspace{3cm}\times\norm{\d_n(s,\d_{0,n})-\d(s,\d_0)}_2|\nabla(\D_v\u_n(s,\v_{0,n})-\xi_v(s,\v_0))|_2ds\\
\leq&\eps\int_0^t|\nabla(\D_v\u_n(s,\v_{0,n})-\xi_v(s,\v_0))|_2^2ds+c|\D_vQ^{-1}|_\infty^2[D_2^n+D_2]^2\int_0^t\norm{\d_n(s,\d_{0,n})-\d(s,\d_0)}_2^2ds.
\end{align*}
For simplicity of notations, we use $\bar D(t)$ to represent $\D_v\d_n(t,\d_{0,n})-\eta_v(t,\d_0)$. Now taking inner product of $\partial_t\bar D(t)$ with $-\Delta\bar D(t)+\bar D(t)$,
\begin{align}
\label{eqn:dH1}
&\qv{\partial_t\bar D(t),-\Delta\bar D(t)+\bar D(t)}=\frac{1}{2}\frac{d|\nabla\bar D(t)|_2^2}{dt}+\frac{1}{2}\frac{d|\bar D(t)|_2^2}{dt}\notag\\
=&-|-\Delta\bar D(t)+\bar D(t)|_2^2-Q(t)\qv{-\Delta\bar D(t)+\bar D(t),B_2(\D_v\u_n(t,\v_{0,n})+\D_vZ(t),\d_n(t,\d_{0,n})-\d(t,\d_0))}\notag\\
&-Q(t)\qv{-\Delta\bar D(t)+\bar D(t),B_2(\D_v\u_n(t,\v_{0,n})-\xi_v(t,\v_0),\d(t,\d_0))}\notag\\
&-Q(t)\qv{-\Delta\bar D(t)+\bar D(t),B_2(\u_n(t,\v_{0,n})-\u(t,\v_0),\D_v\d_n(t,\d_{0,n}))}\notag\\
&-Q(t)\qv{-\Delta\bar D(t)+\bar D(t),B_2(\u(t,\v_{0})+Z(t),\bar D(t))}\notag\\
&-\D_vQ(t)\qv{-\Delta\bar D(t)+\bar D(t),B_2(\u_n(t,\v_{0,n})-\u(t,\v_0),\d(t,\d_0))}\notag\\
&-\D_vQ(t)\qv{-\Delta\bar D(t)+\bar D(t),B_2(\u_n(t,\v_{0,n})+Z(t),\d_n(t,\d_{0,n})-\d(t,\d_0))}\notag\\
&-\qv{-\Delta\bar D(t)+\bar D(t), (f'(\d_n(t,\d_{0,n}))-1)\bar D(t)}\notag\\
&-\qv{-\Delta\bar D(t)+\bar D(t), (f'(\d_n(t,\d_{0,n}))-f'(\d(t,\d_{0})))\eta_v(t,\d_0)}\notag\\
=:&J_1+\cdots +J_9
\end{align}
By H\"older's inequality and Young's inequality,
\begin{align*}
J_2\leq&|Q|_\infty|-\Delta\bar D(t)+\bar D(t)|_2|\D_v\u_n(t,\v_{0,n})+\D_vZ(t)|_2|\nabla(\d_n(t,\d_{0,n})-\d(t,\d_0))|_\infty\\
\leq&\eps|-\Delta\bar D(t)+\bar D(t)|_2^2+c|Q|_\infty^2[M_1^n]^2\norm{\d_n(t,\d_{0,n})-\d(t,\d_0)}_2^2,
\end{align*}
and
\begin{align*}
J_3\leq&|Q|_\infty|-\Delta\bar D(t)+\bar D(t)|_2|\D_v\u_n(t,\v_{0,n})-\xi_v(t,\v_0)|_2|\nabla\d(t,\d_0)|_\infty\\
\leq&\eps|-\Delta\bar D(t)+\bar D(t)|_2^2+c|Q|_\infty^2[D_3]^2|\D_v\u_n(t,\v_{0,n})-\xi_v(t,\v_0)|_2^2.
\end{align*}
Similarly, we get
\begin{align*}
J_4\leq&|Q|_\infty|-\Delta\bar D(t)+\bar D(t)|_2|\u_n(t,\v_{0,n})-\u(t,\v_0)|_\infty|\nabla\D_v\d_n(t,\d_{0,n})|_2\\
\leq&\eps|-\Delta\bar D(t)+\bar D(t)|_2^2+c|Q|_\infty^2[N_2^n]^2\norm{\u_n(t,\v_{0,n})-\u(t,\v_0)}_1^2,
\end{align*}
and
\begin{align*}
J_5\leq&|Q|_\infty|-\Delta\bar D(t)+\bar D(t)|_2|\u(t,\v_{0})+Z(t)|_\infty|\nabla\bar D(t)|_2\\
\leq&\eps|-\Delta\bar D(t)+\bar D(t)|_2^2+c|Q|_\infty^2[C_2]^2|\nabla\bar D(t)|_2^2.
\end{align*}
Still by H\"older's inequality and Young's inequality,
\begin{align*}
J_6+J_7\leq&|\D_vQ|_\infty|-\Delta\bar D(t)+\bar D(t)|_2|\u_n(t,\v_{0,n})-\u(t,\v_0)|_\infty|\nabla\d(t,\d_0)|_2\\
&+|\D_vQ|_\infty|-\Delta\bar D(t)+\bar D(t)|_2|\u_n(t,\v_{0,n})+Z(t)|_2|\nabla(\d_n(t,\d_{0,n})-\d(t,\d_0))|_\infty\\
\leq&\eps|-\Delta\bar D(t)+\bar D(t)|_2^2+c|\D_vQ|_\infty^2[D_2]^2\norm{\u_n(t,\v_{0,n})-\u(t,\v_0)}_1^2\\
&+c|\D_vQ|_\infty^2[C_1^n]^2\norm{\d_n(t,\d_{0,n})-\d(t,\d_0)}_2^2,
\end{align*}
and
\begin{align*}
J_8+J_9\leq&c|-\Delta\bar D(t)+\bar D(t)|_2|f'(\d_n(t,\d_{0,n}))-1|_\infty|\bar D(t)|_2\\
&+c|-\Delta\bar D(t)+\bar D(t)|_2|(f'(\d_n(t,\d_{0,n}))-f'(\d(t,\d_{0})))|_2|\eta_v(t,\d_0)|_\infty\\
\leq&\eps|-\Delta\bar D(t)+\bar D(t)|_2^2+c(1+[D_2^n]^4)|\bar D(t)|_2^2\\
&+c[N_2(D_1^n+D_1)]^2|\d_n(t,\d_{0,n})-\d(t,\d_0)|_2^2.
\end{align*}
Altogether with the estimates in \eqref{eqn:uL2} and \eqref{eqn:dH1}, we have
\begin{align}
&\frac{1}{2}\frac{d|\D_v\u_n(t,\v_{0,n})-\xi_v(t,\v_0)|_2^2}{dt}+|\nabla(\D_v\u_n(t,\v_{0,n})-\xi_v(t,\v_0))|_2^2\notag\\
&+\frac{1}{2}\frac{d|\D_v\d_n(t,\d_{0,n})-\eta_v(t,\d_0)|_2^2}{dt}+\frac{1}{2}\frac{d|\nabla(\D_v\d_n(t,\d_{0,n})-\eta_v(t,\d_0))|_2^2}{dt}\notag\\
&+|-\Delta(\D_v\d_n(t,\d_{0,n})-\eta_v(t,\d_0))+\D_v\d_n(t,\d_{0,n})-\eta_v(t,\d_0)|_2^2\notag\\
\leq&c\norm{\u_n(t,\v_{0,n})-\u(t,\v_0)}_1^2+c|Q^{-1}|_2^2(N_2^n+C_2)^2\norm{\d_n(t,\d_{0,n})-\d(t,\d_0)}_2^2\notag\\
&+c|\D_vQ^{-1}|_\infty^2(D_2^n+D_2)^2\norm{\d_n(t,\d_{0,n})-\d(t,\d_0)}_2^2+c|Q|_\infty^2[M_1^n]^2\norm{\d_n(t,\d_{0,n})-\d(t,\d_0)}_2^2\notag\\
&+c|Q|_\infty^2[N_2^n]^2\norm{\u_n(t,\v_{0,n})-\u(t,\v_0)}_1^2+c|\D_vQ|_\infty^2[D_2]^2\norm{\u_n(t,\v_{0,n})-\u(t,\v_0)}_1^2\notag\\
&+c|\D_vQ|_\infty^2[C_1^n]^2\norm{\d_n(t,\d_{0,n})-\d(t,\d_0)}_2^2+c[N_2(D_1^n+D_1)]^2|\d_n(t,\d_{0,n})-\d(t,\d_0)|_2^2\notag\\
&+c|Q|_\infty C_2|\D_v\u_n(t,\v_{0,n})-\xi_v(t,\v_0)|_2^2+c|Q|_\infty^2[M_2^n+D_3]^2|\D_v\u_n(t,\v_{0,n})-\xi_v(t,\v_0)|_2^2\notag\\
&+c|\D_vQ|_\infty^2[C_2^n+C_2]^2|\D_v\u_n(t,\v_{0,n})-\xi_v(t,\v_0)|_2^2+c|Q^{-1}|_\infty^2[D_3]^2|\nabla(\D_v\d_n(t,\d_{0,n})-\eta_v(t,\d_0))|_2^2\notag\\
&+c|Q|_\infty^2[C_2]^2|\nabla(\D_v\d_n(t,\d_{0,n})-\eta_v(t,\d_0))|_2^2+c(1+[D_2^n]^4)|\D_v\d_n(t,\d_{0,n})-\eta_v(t,\d_0)|_2^2.
\end{align}
Applying Gronwall inequality, we get that
\begin{align}
\label{eqn:3.85}
&\sup_{v\leq t\leq T}\left\{|\D_v\u_n(t,\v_{0,n})-\xi_v(t,\v_0)|_2^2+\norm{\D_v\d_n(t,\d_{0,n})-\eta_v(t,\d_0)}_1^2\right\}\notag\\
&+\int_0^T\norm{\D_v\u_n(t,\v_{0,n})-\xi_v(t,\v_0)}_1^2ds+\int_0^T\norm{\D_v\d_n(t,\d_{0,n})-\eta_v(t,\d_0)}_2^2ds\notag\\
\leq & \int_{0}^{T}L_1(\omega)\norm{\u_n(t,\v_{0,n})-\u(t,\v_0)}_1dt\times\exp \left\{cT(|Q|_\infty C_2+|Q|_\infty^2[M_2^n+D_3]^2+|\D_vQ|_\infty^2[C_2^n+C_2]^2)\right\}\notag\\
&+ \int_{0}^{T}L_2(\omega)\norm{\d_n(s,\d_{0,n})-\d(s,\d_0)}_2^{2}ds\times\exp \left\{cT(|Q^{-1}|_\infty^2[D_3]^2+|Q|_\infty^2[C_2]^2+1+[D_2^n]^4)\right\},
\end{align}
where
\begin{align*}
L_1(\omega)=&cT(1+|Q|_\infty^2[N_2^n]^2+|\D_vQ|_\infty^2[D_2]^2)\\
L_2(\omega)=&cT(|Q^{-1}|_2^2(N_2^n+C_2)^2+|\D_vQ^{-1}|_\infty^2(D_2^n+D_2)^2+|Q|_\infty^2[M_1^n]^2+|\D_vQ|_\infty^2[C_1^n]^2+[N_2(D_1^n+D_1)]^2).
\end{align*}
As we localize $Q,Z$ at the beginning of the proof, they are bounded by $N$. Moreover, since the initial conditions are deterministic, by Theorem \ref{thm:3.2} and Proposition \ref{prop:3.5}, all the norms defined in \eqref{eqn:norm} are uniformly bounded with respect to $\omega,n$.  Hence by (\ref{eqn:3.85}) and dominated convergence theorem,
$$\lim\limits_{n\to\infty}\Ex\sup_{v\leq t\leq T}\{|\D_v\u_n(t,\v_{0,n})-\xi_v(t,\v_0)|_2^2+\norm{\D_v\d_n(t,\d_{0,n})-\eta_v(t,\d_0)}_1^2\}=0.$$
Thus, \eqref{eqn:3.79} gets proved. \qed

\begin{prop}
	\label{prop:3.5}
	For $\v\in\V,\d_0\in\dH^2$, the Malliavin derivative $\D_v\u(t,\v_0),\D_v\d(t,\d_0)$ satisfy the following estimates:
	\begin{align}
	&\sup_{0\leq t\leq T}[|\D_v\u(t,\v_0)|_2^2+\norm{\D_v\d(t,\d_0)}_1^2]+\int_0^T\norm{\D_v\u(t,\v_0)}_1^2dt+\int_0^T\norm{\D_v\d(t,\d_0)}_2^2dt\notag\\\leq& c(|\v_0|_2,\norm{\d_0}_1,|Q|_\infty,\sup\limits_{0\leq t\leq T}\norm{Z}_2,T),
	\end{align}
	and
	\begin{align}
	&\sup_{0\leq t\leq T}[\norm{\D_v\u(t,\v_0)}_1^2+\norm{\D_v\d(t,\d_0)}_2^2]+\int_0^T\norm{\D_v\u(t,\v_0)}_2^2dt+\int_0^T\norm{\D_v\d(t,\d_0)}_3^2dt\notag\\
	\leq&c(\norm{\v_0}_1,\norm{\d_0}_2,|Q|_\infty,\sup\limits_{0\leq t\leq T}\norm{Z}_2,\int_0^T\norm{Z(t)}_3^2dt,T).
	\end{align}
\end{prop}
\proof We first estimate the following:
\begin{align}
|\D_v\u(t,\v_0)|_2^2=&-2\int_0^t|\nabla\D_v\u(s,\v_0)|_2^2ds\notag\\
&-2\int_0^tQ(s)\qv{B_1\left(\D_v\u(s,\v_0)+\D_vZ(s),\u(s,\v_0)+Z(s)\right),\D_v\u(s,\v_0)}ds\notag\\
&-2\int_0^tQ(s)\qv{B_1\left(\u(s,\v_0)+Z(s),\D_v\u(s,\v_0)+\D_vZ(s)\right),\D_v\u(s,\v_0)}ds\notag\\
&-2\int_0^t\D_vQ(s)\qv{B_1(\u(s,\v_0)+Z(s)),\D_v\u(s,\v_0)}ds\notag\\
&-2\int_0^tQ(s)^{-1}\qv{M(\D_v\d(s,\d_0),\d(s,\d_0)),\D_v\u(s,\v_0)}ds\notag\\
&-2\int_0^tQ(s)^{-1}\qv{M(\d(s,\d_0),\D_v\d(s,\d_0)),\D_v\u(s,\v_0)}ds\notag\\
&-2\int_0^t\D_vQ(s)^{-1}\qv{M(\d(s,\d_0)),\D_v\u(s,\v_0)}ds\notag\\
=:&P_1+\cdots+P_7.
\end{align}
By H\"older's inequality and Young's inequality,
\begin{align*}
P_2\leq&c|Q|_\infty\int_0^t|\D_v\u(s,\v_0)+\D_vZ(s)|_2|\nabla(\u(s,\v_0)+Z(s))|_\infty|\D_v\u(s,\v_0)|_2ds\\
\leq& c\int_0^t|\D_v\u(s,\v_0)|_2^2ds+c|Q|_\infty\int_0^t\norm{\u(s,\v_0)+Z(s)}_2^2|\D_v\u(s,\v_0)|_2^2ds\\
&+c|Q|_\infty^2\int_0^t|\D_vZ(s)|_2^2\norm{\u(s,\v_0)+Z(s)}_2^2ds.
\end{align*}
By Lemma \ref{lemma:2.1}, we have
\begin{align*}
P_3=&-2\int_0^tQ(s)\qv{B_1\left(\u(s,\v_0)+Z(s),\D_vZ(s)\right),\D_v\u(s,\v_0)}ds\notag\\
\leq&c|Q|_\infty\int_0^t|\u(s,\v_0)+Z(s)|_\infty|\nabla\D_vZ(s)|_2|\D_v\u(s,\v_0)|_2ds\\
\leq&c\int_0^t|\D_v\u(s,\v_0)|_2^2ds+c|Q|_\infty^2\int_0^t\norm{\u(s,\v_0)+Z(s)}_1^2\norm{\D_vZ(s)}_1^2ds,
\end{align*}
and
\begin{align*}
P_4\leq&c|\D_vQ|_\infty\int_0^t|\u(s,\v_0)+Z(s)|_\infty|\nabla(\u(s,\v_0)+Z(s))|_2|\D_v\u(s,\v_0)|_2ds\\
\leq&c\int_0^t|\D_v\u(s,\v_0)|_2^2ds+c|Q|_\infty^2\int_0^t\norm{\u(s,\v_0)+Z(s)}_1^4ds.
\end{align*}
Similarly,
\begin{align*}
P_5+P_6+P_7\leq&c|Q^{-1}|\int_0^t|\nabla\D_v\d(s,\d_0)|_2|\nabla\d(s,\d_0)|_\infty|\nabla\D_v\u(s,\v_0)|_2\\
&+c|\D_vQ^{-1}|_\infty\int_0^t|\nabla\d(s,\d_0)|_2|\nabla\d(s,\d_0)|_\infty|\nabla\D_v\u(s,\v_0)|_2ds\\
\leq&\eps\int_0^t|\nabla\D_v\u(s,\v_0)|_2ds+c|Q^{-1}|_\infty^2\int_0^t\norm{\d(s,\d_0)}_2^2|\nabla\D_v\d(s,\d_0)|_2^2ds\\
&+c|\D_vQ^{-1}|_\infty^2\int_0^t\norm{\d(s,\d_0)}_1^2\norm{\d(s,\d_0)}_2^2ds.
\end{align*}
We now take inner product of $\partial_t\D_v\d(t,\d_0)$ with $-\Delta\D_v\d(t,\d_0)+\D_v\d(t,\d_0)$ in $L^2(\mathbf{D})$,
\begin{align}
&\qv{-\Delta\D_v\d(t,\d_0)+\D_v\d(t,\d_0),\partial_t\D_v\d(t,\d_0)}=\frac{1}{2}\frac{d|\D_v\d(t,\d_0)|_2^2}{dt}+\frac{1}{2}\frac{d|\nabla\D_v\d(t,\d_0)|_2^2}{dt}\notag\\
=&-|-\Delta\D_v\d(t,\d_0)+\D_v\d(t,\d_0)|_2^2\notag\\
&-Q(t)\qv{-\Delta\D_v\d(t,\d_0)+\D_v\d(t,\d_0),B_2(\D_v\u(t,\v_0)+\D_vZ(t), \d(t,\d_0))}\notag\\
&-Q(t)\qv{-\Delta\D_v\d(t,\d_0)+\D_v\d(t,\d_0),B_2\left(\u(t,\v_0)+Z(t),\D_v\d(t,\d_0)\right)}\notag\\
&-\D_vQ(t)\qv{-\Delta\D_v\d(t,\d_0)+\D_v\d(t,\d_0),B_2(\u(t,\v_0)+Z(t),\d(t,\d_0))}\notag\\
&-\qv{-\Delta\D_v\d(t,\d_0)+\D_v\d(t,\d_0),f'(\d(t,\d_0))\D_v\d(t,\d_0)-\D_v\d(t,\d_0)}ds\notag\\
=:&Q_1+\cdots+Q_5.
\end{align}
By H\"older's inequality and Young's inequality, we get
\begin{align*}
Q_2+Q_3\leq&|Q|_\infty|-\Delta\D_v\d(t,\d_0)+\D_v\d(t,\d_0)|_2|\D_v\u(t,\v_0)+\D_vZ(t)|_2|\nabla\d(t,\d_0)|_\infty\\
&+|Q|_\infty|-\Delta\D_v\d(t,\d_0)+\D_v\d(t,\d_0)|_2|\u(t,\v_0)+Z(t)|_\infty|\nabla\D_v\d(t,\d_0)|_2\\
\leq&\eps|-\Delta\D_v\d(t,\d_0)+\D_v\d(t,\d_0)|_2^2+c|Q|_\infty^2\norm{\d(t,\d_0)}_2^2|\D_v\u(t,\v_0)|_2^2\\
&+c|Q|_\infty^2\norm{\d(t,\d_0)}_2^2|\D_vZ(t)|_2^2+c|Q|_\infty^2\norm{\u(t,\v_0)+Z(t)}_1^2|\nabla\D_v\d(t,\d_0)|_2^2.
\end{align*}
Similarly,
\begin{align*}
Q_4\leq& c|\D_vQ|_\infty|-\Delta\D_v\d(t,\d_0)+\D_v\d(t,\d_0)|_2|\u(t,\v_0)+Z(t)|_\infty|\nabla\d(t,\d_0)|_2\\
\leq&\eps|-\Delta\D_v\d(t,\d_0)+\D_v\d(t,\d_0)|_2^2+c|\D_vQ|_\infty^2\norm{\u(t,\v_0)+Z(t)}_1^2\norm{\d(t,\d_0)}_1^2.
\end{align*}
Finally,  we have
\begin{align*}
Q_5\leq& c|-\Delta\D_v\d(t,\d_0)+\D_v\d(t,\d_0)|_2|f'(\d(t,\d_0))-1|_\infty|\D_v\d(t,\d_0)|_2\\
\leq&\eps|-\Delta\D_v\d(t,\d_0)+\D_v\d(t,\d_0)|_2^2+c(1+\norm{\d(t,\d_0)}_1^4)|\D_v\d(t,\d_0)|_2^2.
\end{align*}
Altogether, we get
\begin{align}
&\frac{1}{2}\frac{d|\D_v\u(t,\v_0)|_2^2}{dt}+\frac{1}{2}\frac{d|\D_v\d(t,\d_0)|_2^2}{dt}+\frac{1}{2}\frac{d|\nabla\D_v\d(t,\d_0)|_2^2}{dt}\notag\\
&+|\nabla\D_v\u(t,\v_0)|_2^2+|-\Delta\D_v\d(t,\d_0)+\D_v\d(t,\d_0)|_2^2\notag\\
\leq&c|Q|_\infty^2|\D_vZ(t)|_2^2\norm{\u(t,\v_0)+Z(t)}_2^2+c|Q|_\infty\norm{\u(t,\v_0)+Z(t)}_1^2\norm{\D_vZ(t)}_1^2\notag\\
&+c|Q|_\infty^2\norm{\u(t,\v_0)+Z(t)}_1^4+c|\D_vQ^{-1}|^2_\infty\norm{\d(t,\d_0)}_1^2\norm{\d(t,\d_0)}_2^2\notag\\
&+c|Q|_\infty^2\norm{\d(t,\d_0)}_2^2|\D_vZ(t)|_2^2+c|\D_vQ|_\infty^2\norm{\u(t,\v_0)+Z(t)}_1^2\norm{\d(t,\d_0)}_1^2\notag\\
&+c|Q|_\infty\norm{\u(t,\v_0)+Z(t)}_2^2|\D_v\u(t,\v_0)|_2^2+c|\D_v\u(t,\v_0)|_2^2+c|Q|_\infty^2\norm{\d(t,\d_0)}_2^2|\D_v\u(t,\v_0)|_2^2\notag\\
&+c|Q^{-1}|_\infty^2\norm{\d(t,\d_0)}_2^2|\nabla\D_v\d(t,\d_0)|_2^2+c|Q|_\infty^2\norm{\u(t,\v_0)+Z(t)}_1^2|\nabla\D_v\d(t,\d_0)|_2^2\notag\\
&+c(1+\norm{\d(t,\d_0)}_1^4)|\D_v\d(t,\d_0)|_2^2.
\end{align}
Applying Gronwall's inequality yields that
\begin{align}\label{eqn:3.90}
&\sup_{0\leq t\leq T}[|\D_v\u(t,\v_0)|_2^2+\norm{\D_v\d(t,\d_0)}_1^2]+\int_0^T\norm{\D_v\u(t,\v_0)}_1^2dt+\int_0^T|-\Delta\D_v\d(t,\d_0)+\D_v\d(t,\d_0)|_2^2dt\notag\\
\leq& ch_1(T)h_2(T),
\end{align}
where
\begin{align*}
h_1(T):=&|Q|_\infty^2\sup_{0\leq t\leq T}|\D_vZ(t)|_2^2\int_0^T\norm{\u(t,\v_0)+Z(t)}_2^2dt+|Q|_\infty\sup_{0\leq t\leq T}\norm{\u(t,\v_0)+Z(t)}_1^2\norm{\D_vZ(t)}_1^2T\\
&+|Q|_\infty^2\sup_{0\leq t\leq T}\norm{\u(t,\v_0)+Z(t)}_1^4T+|\D_vQ^{-1}|_\infty\sup_{0\leq t\leq T}\norm{\d(t,\d_0)}_1^2\norm{\d(t,\d_0)}_2^2T\\
&+|Q|_\infty^2\sup_{0\leq t\leq T}\norm{\d(t,\d_0)}_2^2|\D_vZ(t)|_2^2T+|\D_vQ|_\infty^2\sup_{0\leq t\leq T}\norm{\u(t,\v_0)+Z(t)}_1^2\norm{\d(t,\d_0)}_1^2T,
\end{align*}
and
\begin{align*}
h_2(T):=\exp c\bigg\{&T+|Q|_\infty\int_0^T\norm{\u(t,\v_0)+Z(t)}_2^2dt+|Q|_\infty^2\sup_{0\leq t\leq T}\norm{\d(t,\d_0)}_2^2T\\
&+|Q^{-1}|_\infty^2\sup_{0\leq t\leq T}\norm{\d(t,\d_0)}_2^2T+|Q|_\infty^2\sup_{0\leq t\leq T}\norm{\u(t,\v_0)+Z(t)}_1^2T+\sup_{0\leq t\leq T}\norm{\d(t,\d_0)}_1^4T\bigg\}.
\end{align*}
In view of Proposition \ref{prop:3.1}, $h_1(T), h_2(T)$ are constants depending on $\norm{\v_0}_1,\norm{\d_0}_2,|Q|_\infty,\sup\limits_{0\leq t\leq T}\norm{Z}_2,T$.
Now taking inner product between $\Delta\D_v\d(t,\d_0)$ and $\Delta\partial_t\D_v\d(t,\d_0)$ in $L^2(\mathbf{D})$ gives that
\begin{align}
&\qv{\Delta\D_v\d(t,\d_0),\Delta\partial_t\D_v\d(t,\d_0)}=\frac{1}{2}\partial_t|\Delta\D_v\d(t,\d_0)|_2^2=-\qv{\nabla\Delta\D_v\d(t,\d_0),\nabla\partial_t\D_v\d(t,\d_0)}\notag\\
=&-|\nabla\Delta\D_v\d(t,\d_0)|_2^2+Q(t)\qv{\nabla\Delta\D_v\d(t,\d_0),B_2(\nabla(\D_v\u(t,\v_0)+\D_vZ(t)),\d(t,\d_0))}\notag\\
&+Q(t)\qv{\nabla\Delta\D_v\d(t,\d_0),B_2(\D_v\u(t,\v_0)+\D_vZ(t),\nabla\d(t,\d_0)}\notag\\
&+Q(t)\qv{\nabla\Delta\D_v\d(t,\d_0),B_2(\nabla(\u(t,\v_0)+Z(t)),\D_v\d(t,\d_0))}\notag\\
&+Q(t)\qv{\nabla\Delta\D_v\d(t,\d_0),B_2(\u(t,\v_0)+Z(t),\nabla\D_v\d(t,\d_0))}\notag\\
&+\D_vQ(t)\qv{\nabla\Delta\D_v\d(t,\d_0),B_2(\nabla(\u(t,\v_0)+Z(t)),\d(t,\d_0))}\notag\\
&+\D_vQ(t)\qv{\nabla\Delta\D_v\d(t,\d_0),B_2(\u(t,\v_0)+Z(t),\nabla\d(t,\d_0))}\notag\\
&+\qv{\nabla\Delta\D_v\d(t,\d_0),\nabla(f'(\d(t,\d_0))\D_v\d(t,\d_0))}.
\end{align}
With a similar discussion, we get that
\begin{align}
&\partial_t|\Delta\D_v\d(t,\d_0)|_2^2+|\nabla\Delta\D_v\d(t,\d_0)|_2^2\notag\\
\leq&c|Q|_\infty^2\norm{\d(t,\d_0)}_2^2|\nabla\D_v\u(t,\v_0)|_2^2+c|Q|_\infty^2|\nabla\D_vZ(t)|_2^2\norm{\d(t,\d_0)}_2^2\notag\\
&+c|Q|_\infty^2\norm{\d(t,\d_0)}_3^2|\D_v\u(t,\v_0)|_2^2+c|Q|_\infty^2\norm{\d(t,\d_0)}_3^2|\D_vZ(t)|_2^2\notag\\
&+c|Q|_\infty^2\norm{\u(t,\v_0)+Z(t)}_2^2|\nabla\D_v\d(t,\d_0)|_2^2+c|Q|_\infty^2\norm{\u(t,\v_0)+Z(t)}_1^2|\Delta\D_v\d(t,\d_0)|_2^2\notag\\
&+c|\D_vQ|_\infty^2\norm{\u(t,\v_0)+Z(t)}_1^2\norm{\d(t,\d_0)}_2^2+c\norm{\d(t,\d_0)}_1^2\norm{\d(t,\d_0)}_2^2|\D_v\d(t,\d_0)|_2^2\notag\\
&+c\norm{\d(t,\d_0)}_1^4|\nabla\D_v\d(t,\d_0)|_2^2.
\end{align}
Applying Gronwall's inequality, and combining the estimate in \eqref{eqn:3.90}, we conclude that
\begin{align}\label{eqn:3.93}
&\sup_{0\leq t\leq T}[|\D_v\u(t,\v_0)|_2^2+\norm{\D_v\d(t,\d_0)}_2^2]+\int_0^T\norm{\D_v\u(t,\v_0)}_1^2dt+\int_0^T\norm{\D_v\d(t,\d_0)}_3^2dt\notag\\
\leq&ch_1(T)h_2(T)h_3(T),	
\end{align}
where
\begin{align*}
h_3(T):=&\bigg[\sup_{0\leq t\leq T}|Q|_\infty^2\norm{\d(t,\d_0)}_2^2+|Q|_\infty^2\sup_{0\leq t\leq T}|\nabla\D_vZ(t)|_2^2\norm{\d(t,\d_0)}_2^2T\notag\\
&\hspace{5mm}+|Q|_\infty^2\int_0^T\norm{\d(t,\d_0)}_3^2dt+|Q|_\infty^2\sup_{0\leq t\leq T}|\D_vZ(t)|_2^2\int_0^T\norm{\d(t,\d_0)}_3^2dt\\
&\hspace{5mm}+|Q|_\infty^2\int_0^T\norm{\u(t,\v_0)+Z(t)}_2^2dt+|\D_vQ|_\infty^2\sup_{0\leq t\leq T}\norm{\u(t,\v_0)+Z(t)}_1^2\norm{\d(t,\d_0)}_2^2T\\
&\hspace{5mm}+\sup_{0\leq t\leq T}\norm{\d(t,\d_0)}_1^2\norm{\d(t,\d_0)}_2^2T+\sup_{0\leq t\leq T}\norm{\d(t,\d_0)}_1^4T\bigg]\\
&\times\exp\{c|Q|_\infty^2\sup_{0\leq t\leq T}\norm{\u(t,\v_0)+Z(t)}_1^2T\}.
\end{align*}
In view of Proposition \ref{prop:3.1}, $h_3(T)$ is a constant depending on $\norm{\v_0}_1,\norm{\d_0}_2,|Q|_\infty,\sup\limits_{0\leq t\leq T}\norm{Z}_2,T$. Finally, we take inner product between $\D_v\u(t,\v_0)$ and $-\Delta\D_v\u(t,\v_0)$ in $L^2(\mathbf{D})$,
\begin{align}
\frac{d|\nabla\D_v\u(t,\v_0)|_2^2}{dt}=&-2|\Delta\D_v\u(t,\v_0)|_2^2\notag\\
&+2Q(t)\qv{B_1\left(\D_v\u(t,\v_0)+\D_vZ(t),\u(t,\v_0)+Z(t)\right),\Delta\D_v\u(t,\v_0)}\notag\\
&+2Q(t)\qv{B_1\left(\u(t,\v_0)+Z(t),\D_v\u(t,\v_0)+\D_vZ(t)\right),\Delta\D_v\u(t,\v_0)}\notag\\
&+2 \D_vQ(t)\qv{B_1(\u(t,\v_0)+Z(t)),\Delta\D_v\u(t,\v_0)} \notag\\
&+2 Q(t)^{-1}\qv{M(\D_v\d(t,\d_0),\d(t,\d_0)),\Delta\D_v\u(t,\v_0)} \notag\\
&+2Q(t)^{-1}\qv{M(\d(t,\d_0),\D_v\d(t,\d_0)),\Delta\D_v\u(t,\v_0)} \notag\\
&+2 \D_vQ(t)^{-1}\qv{M(\d(t,\d_0)),\Delta\D_v\u(t,\v_0)}\notag\\
=:&Q_1+\cdots+Q_7.
\end{align}
By H\"older's inequality and Young's inequality,
\begin{align*}
Q_2+Q_3\leq&c|Q|_\infty|\D_v\u(t,\v_0)+\D_vZ(t)|_2|\nabla(\u(t,\v_0)+Z(t))|_\infty|\Delta\D_v\u(t,\v_0)|_2\notag\\
&+c|Q|_\infty|\u(t,\v_0)+Z(t)|_\infty|\nabla(\D_v\u(t,\v_0)+\D_vZ(t))|_2|\Delta\D_v\u(t,\v_0)|_2\notag\\
\leq&\eps|\Delta\D_v\u(t,\v_0)|_2^2+c|Q|_\infty^2\norm{\u(t,\v_0)+Z(t)}_2^2(|\D_v\u(t,\v_0)|_2^2+|\D_vZ(t)|_2^2)\notag\\
&+c|Q|_\infty^2\norm{\u(t,\v_0)+Z(t)}_1^2(|\nabla\D_v\u(t,\v_0)|_2^2+|\nabla\D_vZ(t)|_2^2),
\end{align*}
and
\begin{align*}
Q_4\leq&c|\D_vQ|_\infty|\u(t,\v_0)+Z(t)|_\infty|\nabla(\u(t,\v_0)+Z(t))|_2|\Delta\D_v\u(t,\v_0)|_2\\
\leq&\eps|\Delta\D_v\u(t,\v_0)|_2^2+c|\D_vQ|_\infty^2\norm{\u(t,\v_0)+Z(t)}_1^4.
\end{align*}
By Proposition \ref{prop:2.5}, H\"older's inequality and Young's inequality, we get
\begin{align*}
Q_5+Q_6+Q_7=&2Q(t)^{-1}\qv{B_2(\Delta\D_v\u(t,\v_0),\d(t,\d_0)),\Delta\D_v\d(t,\d_0)}\\
&+2Q(t)^{-1}\qv{B_2(\Delta\D_v\u(t,\v_0),\D_v\d(t,\d_0)),\Delta\d(t,\d_0)}\\
&+2\D_vQ(t)^{-1}\qv{B_2(\Delta\D_v\u(t,\v_0),\d(t,\d_0)),\Delta\d(t,\d_0)}\\
\leq&c|Q^{-1}|_\infty|\Delta\D_v\u(t,\v_0)|_2|\nabla\d(t,\d_0)|_\infty|\Delta\D_v\d(t,\d_0)|_2\\
&+c|Q^{-1}|_\infty|\Delta\D_v\u(t,\v_0)|_2|\Delta\d(t,\d_0)|_\infty|\nabla\D_v\d(t,\d_0)|_2\\
&+c|\D_vQ^{-1}|_\infty|\Delta\D_v\u(t,\v_0)|_2|\nabla\d(t,\d_0)|_\infty|\Delta\d(t,\d_0)|_2\\
\leq&\eps|\Delta\D_v\u(t,\v_0)|_2^2+c|Q^{-1}|_\infty^2\norm{\d(t,\d_0)}_2^2|\Delta\D_v\d(t,\d_0)|_2^2\\
&+c|Q^{-1}|_\infty^2\norm{\d(t,\d_0)}_3^2|\nabla\D_v\d(t,\d_0)|_2^2+c|\D_vQ^{-1}|_\infty^2\norm{\d(t,\d_0)}_2^4.
\end{align*}
Altogether, we get
\begin{align}
&\frac{d|\nabla\D_v\u(t,\v_0)|_2^2}{dt}+|\Delta\D_v\u(t,\v_0)|_2^2\notag\\
\leq&c|Q|_\infty^2\norm{\u(t,\v_0)+Z(t)}_2^2(|\D_v\u(t,\v_0)|_2^2+|\D_vZ(t)|_2^2)\notag\\
&+c|Q|_\infty^2\norm{\u(t,\v_0)+Z(t)}_1^2(|\nabla\D_v\u(t,\v_0)|_2^2+|\nabla\D_vZ(t)|_2^2)+c|\D_vQ|_\infty^2\norm{\u(t,\v_0)+Z(t)}_1^4\notag\\
&+c|Q^{-1}|_\infty^2\norm{\d(t,\d_0)}_2^2|\Delta\D_v\d(t,\d_0)|_2^2+c|Q^{-1}|_\infty^2\norm{\d(t,\d_0)}_3^2|\nabla\D_v\d(t,\d_0)|_2^2\notag\\
&+c|\D_vQ^{-1}|_\infty^2\norm{\d(t,\d_0)}_2^4.
\end{align}
Applying Gronwall's inequality, together with the estimates \eqref{eqn:3.90} and \eqref{eqn:3.93}, we arrive at the required result:
\begin{align}
&\sup_{0\leq t\leq T}[\norm{\D_v\u(t,\v_0)}_1^2+\norm{\D_v\d(t,\d_0)}_2^2]+\int_0^T\norm{\D_v\u(t,\v_0)}_2^2dt+\int_0^T\norm{\D_v\d(t,\d_0)}_3^2dt\notag\\
\leq& ch_1(T)h_2(T)h_3(T)h_4(T),
\end{align}
where
\begin{align*}
h_4(T):=&\bigg[|Q|_\infty^2(1+\sup_{0\leq t\leq T}|\D_vZ(t)|_2^2)\int_0^T\norm{\u(t,\v_0)+Z(t)}_2^2dt\\
&\hspace{5mm}+|Q|_\infty^2\sup_{0\leq t\leq T}\norm{\u(t,\v_0)+Z(t)}_1^2|\nabla\D_vZ(t)|_2^2T+|\D_vQ|_\infty^2\sup_{0\leq t\leq T}\norm{\u(t,\v_0)+Z(t)}_1^4T\\
&\hspace{5mm}+|Q^{-1}|_\infty^2\sup_{0\leq t\leq T}\norm{\d(t,\d_0)}_2^2T+|Q^{-1}|_\infty^2\int_0^T\norm{\d(t,\d_0)}_3^2dt+|\D_vQ^{-1}|_\infty^2\sup_{0\leq t\leq T}\norm{\d(t,\d_0)}_2^4T\bigg]\\
&\times\exp\{c|Q|_\infty^2\sup_{0\leq t\leq T}\norm{\u(t,\v_0)+Z(t)}_1^2T\}.
\end{align*}
The proof is done, since in view of Proposition \ref{prop:3.1}, $h_4(T)$ is a constant depending on $\norm{\v_0}_1,\norm{\d_0}_2,|Q|_\infty, \sup\limits_{0\leq t\leq T}\norm{Z}_2$, $\int_0^T\norm{Z(t)}_3^2dt,T$.

\qed
%%%%%%%%%%%%%%%%%%%%%%%
\section{The global well-posedness of stochastic nematic liquid crystals flows with random initial and random boundary conditions}
%%%%%%%%%%%%%%%%%%%%%%%
In this section, we will prove our main result Theorem \ref{thm:main}, the global well-posedness for the stochastic nematic liquid crystals flows with random initial conditions.

\proof To show the Theorem \ref{thm:main}, we first reformulate the equations (\ref{eqn:1.1})-(\ref{eqn:1.3}) in the following equivalent integral form

	\begin{align}
	\label{eqn:4.1}
		\v(t,R_{\nu})=&\v_0-\int_0^tA_1\v(s,R_{\nu})ds-\int_0^tB_1(\v(s,R_{\nu}))ds-\int_0^tM(\d(s,R_d))ds\notag\\
		&+\int_0^t\v(s,R_{\nu})\circ dW(s)+\sigma_0W_0(t),\\
		\d(t,R_d)=&\d_0-\int_0^tA_2\d(s,R_d)ds-\int_0^tB_2(\v(s,R_{\nu}),\d(s,R_d))ds-\int_0^tf(\d(s,R_d))ds.
	\end{align}

Inspired by \cite{MZ}, we try to utilize Malliavin calculus techniques to show the global well-posedness of random initial conditions and random boundary conditions problem (\ref{eqn:1.1})-(\ref{eqn:1.3}). For any fixed $t>0$, denote by $\{0=t_0<t_1<\cdots t_n=n\}$ an
arbitrary partition such that $\tau_n:=\max\limits_{1\leq k\leq n}(t_{k}-t_{k-1})\to0$ as $n\to\infty$. Then we have
\begin{align}
\label{eqn:4.3}
&\v(t,R_{\nu})-R_{\nu}-Q(t)Z(t)=Q(t)\u(t,R_{\nu})-R_{\nu}\notag\\
=&\sum_{k=1}^{n}\left(Q(t_k)\u(t_k,R_{\nu})-Q(t_{k-1})\u(t_{k-1},R_{\nu})\right)\notag\\
=&\sum_{k=1}^{n}Q(t_k)(\u(t_k,R_{\nu})-\u(t_{k-1},R_{\nu}))+\sum_{k=1}^{n}(Q(t_k)-Q(t_{k-1}))\u(t_{k-1},R_{\nu})\notag\\
=&-\sum_{k=1}^{n}\int_{t_{k-1}}^{t_k}Q(t_k)A_1\u(s,R_{\nu})ds-\sum_{k=1}^{n}\int_{t_{k-1}}^{t_k}Q(t_k)Q(s)B_1(\u(s,R_{\nu})+Z(s))ds\notag\\
&-\sum_{k=1}^{n}\int_{t_{k-1}}^{t_k}Q(t_k)Q(s)^{-1}M(\d(s,R_d))ds+\sum_{k=1}^{n}\u(t_{k-1},R_{\nu})\int_{t_{k-1}}^{t_k}Q(s)dW(s)+\frac{1}{2}\sum_{k=1}^{n}\u(t_{k-1},R_{\nu})\int_{t_{k-1}}^{t_k}Q(s)ds\notag\\
=:&I_1^n+\cdots+I_5^n;
\end{align}
\begin{align}
	&\d(t,R_d)-\d_0=\sum_{k=1}^{n}\d(t_k,R_d)-\d(t_{k-1},R_d)\notag\\
	=&-\sum_{k=1}^{n}\int_{t_{k-1}}^{t_k}A_2\d(s,R_d)ds-\sum_{k=1}^{n}\int_{t_{k-1}}^{t_k}Q(s)B_2(\u(s,R_{\nu})+Z(s),\d(s,R_d))ds-\sum_{k=1}^{n}\int_{t_{k-1}}^{t_k}f(\d(s,R_d))ds\notag\\
	=:&J_1^n+J_2^n+J_3^n.
\end{align}
Since $\u(s,R_{\nu}),\d(s,R_d)$ and $Q(s)$ are contiuous with respect to time $s$, we have
\begin{align*}
	&\lim_{n\to\infty}I_1^n=-\int_0^tA_1\v(s,R_{\nu})ds+\int_0^tQ(s)A_1Z(s)ds,\\
	&\lim_{n\to\infty}I_2^n=-\int_0^tQ(s)^2B_1(\u(s,R_{\nu})+Z(s))ds=-\int_0^tB_1(\v(s,R_{\nu}))ds,\\
	&\lim_{n\to\infty}I_3^n=-\int_0^tM(\d(s,R_d))ds,\hspace{1cm}\lim_{n\to\infty}I_5^n=\frac{1}{2}\int_0^tQ(s)\u(s,R_{\nu})ds;\\
	&\lim_{n\to\infty}J_1^n=-\int_0^tA_2\d(s,R_d)ds,\\
	&\lim_{n\to\infty}J_2^n=-\int_0^tQ(s)B_2(\u(s,R_{\nu})+Z(s),\d(s,R_d))ds=-\int_0^tB_2(\v(s,R_{\nu}),\d(s,R_d))ds,\\
	&\lim_{n\to\infty}J_3^n=-\int_0^tf(\d(s,R_d))ds.
\end{align*}
It remains to deal with $I_4^n$. By the property of the Skorohod integral (See Proposition 1.3.5 in \cite{N}) we get
\begin{align*}
	I_4^n=&-\sum_{k=1}^{n}\int_{t_{k-1}}^{t_k}\u(t_{k-1},R_{\nu})Q(s)dW(s)+\sum_{k=1}^{n}\int_{t_{k-1}}^{t_k}\D_s(\u(t_{k-1},R_{\nu}))Q(s)ds\\
	=&\int_0^t\sum_{k=1}^{n}\u(t_{k-1},R_{\nu})Q(s)1_{(t_{k-1},t_k]}(s)dW(s)+\int_0^t\sum_{k=1}^{n}\D_s(\u(t_{k-1},R_{\nu}))Q(s)1_{(t_{k-1},t_k]}(s)ds\\
	=:&\int_0^tK^n(s)dW(s)+\int_0^tL^n(s)ds.
\end{align*}
Denote by $\mathbb L^{1,2}(\vH)$ the class of $\vH$-valued process $\v(t)\in\D^{1,2}(\vH)$ for almost all $t$, and there always exists a measurable version of $\D_s\v(t)$ satsifying $\Ex\left(\int_0^T\int_0^T|\D_s\v(t)|_2^2dsdt\right)<\infty$. We say $\v(t)\in\mathbb L^{1,2}_{\text{loc}}(\vH)$ if there exists a sequence $\{\Omega_n\}\subset\F$ such that $\Omega_n$ increases to $\Omega$ and $\v1_{\Omega_n}\in \mathbb L^{1,2}(\vH)$. Without loss of generality, we can assume $\norm{R_{\nu}}_1\leq M, \|R_{d}\|_{2}\leq M$, $\norm{Z(s)}_2\leq M$ and $Q=Q_N$, or we can always do truncation otherwise. As $\u(s,R_{\nu})$ is continuous in $s$, we have for all $s>0$, $K^n(s)\to\u(s,R_{\nu})Q(s)$, and by Proposition \ref{prop:3.1} and the localization, we have
\begin{align}
	\sup_{0\leq s\leq t}|K^n(s)|_2\leq &\sup_{0\leq s\leq t}|Q(s)|\sup_{0\leq s\leq t}|\u(s,R_{\nu})|_2\notag\\
	\leq&\sup_{0\leq s\leq t}|Q(s)|c(|R_{\nu}|_2,\norm{R_d}_1,\norm{Z}_2).
\end{align}
Applying domniated convergence theorem yields that
\begin{equation}
\label{eqn:4.6}
	\lim_{n\to\infty}\Ex\left[\int_0^t|K^n(s)-\u(s,R_{\nu})Q(s)|_2^2\right]=0.
\end{equation}
Moreover, the Malliavin derivative of $K^n$ is given by
\begin{align}
	\D_vK^n(s)=&\sum_{k=1}^{n}[Q(s)\D_v\u(t_{k-1},R_{\nu})+\D_vQ(s)\u(t_{k-1},R_{\nu})]1_{(t_{k-1},t_k]}(s)\notag\\
	=&\sum_{k=1}^{n}\D_vQ(s)\u(t_{k-1},R_{\nu})1_{(t_{k-1},t_k]}(s)\notag\\
	&+\sum_{k=1}^{n}Q(s)[\D_v\u(t_{k-1},R_{\nu})+D\u(t_{k-1},R_{\nu})\D_vR_{\nu}]1_{(t_{k-1},t_k]}(s),
\end{align}
where $D\u(s,\gamma)$ represents the Fr\'echet derivative at $\gamma\in\V$, and $\D_v\u(t_{k-1},R_{\nu}):=\D_v\u(t_{k-1},\gamma)|_{\gamma=R_{\nu}}$. As $\u(s,R_{\nu}),D\u(s,R_{\nu})$ and $\D_v\u(s,R_{\nu})$ are continuous in $s$, we have $\lim\limits_{n\to\infty}\D_v K^n(s)=\D_v[\u(s,R_{\nu})Q(s)]$ for any $s\geq0$, and by Proposition \ref{prop:3.1}, Theorem \ref{thm:3.2} and Theorem \ref{prop:3.4}, we get
\begin{align}
	|\D_vK^n(s)|_2\leq&
	c(\norm{R_{\nu}}_1,\|R_{d}\|_{2},|Q|_\infty, |Q^{-1}|_\infty, \sup_{0\leq t\leq T}\norm{Z}_2,\int_0^T\norm{Z}_3^2ds,T).
\end{align}
Hence, following the dominated convergence theorem yields that
\begin{equation}
\label{eqn:4.9}
	\lim_{n\to\infty}\Ex\left[\int_0^t\int_0^\infty|\D_v K^n(s)-\D_v[\u(s,R_{\nu})Q(s)]|_2^2dvds\right]
\end{equation}
From \eqref{eqn:4.6} and \eqref{eqn:4.9}, we conclude that $K^n(\cdot)\to\u(\cdot,Y_0)Q(\cdot)$ in $\mathbb L_{\text{loc}}^{1,2}(\vH)$. Hence,
\begin{equation}
	\lim_{n\to\infty}\int_0^tK^n(s)dW(s)=\int_0^t\u(s,R_{\nu})Q(s)dW(s).
\end{equation}
Now we estimate $L^n(s)$,
\begin{align}
	L^n(s)=&Q(s)\sum_{k=1}^{n}[\D_s\u(t_{k-1},R_{\nu})+D\u(t_{k-1},R_{\nu})\D_sR_{\nu}]1_{(t_{k-1},t_k]}(s)\notag\\
	=&Q(s)\sum_{k=1}^{n}D\u(t_{k-1},R_{\nu})\D_sR_{\nu}1_{(t_{k-1},t_k]}(s).
\end{align}
Since $D\u(s,R_{\nu})$ is continuous in $s$, we get that
\begin{equation}
	\lim_{n\to\infty}\int_0^tL^n(s)ds=\int_0^tQ(s)D\u(s,R_{\nu})\D_sR_{\nu}ds.
\end{equation}
Back to \eqref{eqn:4.3}, sending $n\to\infty$ yields that
\begin{align}
	\v(t,R_{\nu})-R_{\nu}-Q(t)Z(t)=&Q(t)\u(t,R_{\nu})-R_{\nu}=\sum_{k=1}^{n}\left(Q(t_k)\u(t_k,R_{\nu})-Q(t_{k-1})\u(t_{k-1},R_{\nu})\right)\notag\\
	=&-\int_0^tA_1\v(s,R_{\nu})ds-\int_0^tB_1(\v(s,R_{\nu}))ds-\int_0^tM(\d(s,R_d))ds\notag\\
	&+\int_0^tQ(s)A_1Z(s)ds+\int_0^tQ(s)\u(s,R_{\nu})dW(s)\notag\\
	&+\frac{1}{2}\int_0^tQ(s)\u(s,R_{\nu})ds+\int_0^tQ(s)D\u(s,R_{\nu})\D_sR_{\nu}ds,
\end{align}
for all $t\geq0$. By It\^o's formula, we first have
\begin{align}
	Q(t)Z(t)=&\int_0^tZ(s)\circ dQ(s)+\int_0^tQ(s)\circ dZ(s)\notag\\
	=&-\int_0^tQ(s)A_1Z(s)ds+\sigma_0W_0(t)+\int_0^tZ(s)Q(s)\circ dW(s).
\end{align}
To obtain the form \eqref{eqn:4.1}, it remains to show for $t\geq0$,
\begin{align}
\label{eqn:4.15}
	\int_0^t\u(s,R_{\nu})Q(s)\circ dW(s)=&\int_0^tQ(s)\u(s,R_{\nu})dW(s)+\frac{1}{2}\int_0^tQ(s)\u(s,R_{\nu})ds\notag\\
	&+\int_0^tQ(s)D\u(s,R_{\nu})\D_sR_{\nu}ds
\end{align}
In view of Theorem 3.1.1 in \cite{N}, the left hand side can be written as
\begin{equation}
		\int_0^t\u(s,R_{\nu})Q(s)\circ dW(s)=\int_0^tQ(s)\u(s,R_{\nu})dW(s)+\frac{1}{2}\int_0^t(\nabla[\u(\cdot,R_{\nu})Q(\cdot)])_sds,
\end{equation}
where
\begin{equation}
(\nabla[\u(\cdot,R_{\nu})Q(\cdot)])(s)=\frac{1}{2}\left(\lim_{\eps\to0+}\D_s[\u(s+\eps,R_{\nu})Q(s+\eps)]+\lim_{\eps\to0+}\D_s[\u(s-\eps,R_{\nu})Q(s-\eps)]\right)
\end{equation}
By chain rule, we know that
\begin{equation}
	\D_s[\u(t,R_{\nu})Q(t)]=\D_s\u(t,R_{\nu})Q(t)+D\u(t,R_{\nu})(\D_sR_{\nu})Q(t)+\u(t,R_{\nu})\D_sQ(t).
\end{equation}
Now replacing $t$ in the above identity by $s+\eps$, $s-\eps$, respectively, and using the face that
$\D_s\u(s-\eps,R_{\nu})=0, \D_sQ(s-\eps)=0$, one can get
\begin{align}
	\D_s[\u(s-\eps,R_{\nu})Q(s-\eps)]=&D\u(s-\eps,R_{\nu})(\D_sR_{\nu})Q(s-\eps);\\
	\D_s[\u(s+\eps,R_{\nu})Q(s+\eps)]=&\D_s\u(s+\eps,R_{\nu})Q(s+\eps)+D\u(s+\eps,R_{\nu})(\D_sR_{\nu})Q(s+\eps)\notag\\
	&+\u(s+\eps,R_{\nu})\D_sQ(s+\eps).
\end{align}
Sending $\eps\to0+$, by the continuity of $\D_s\u(t,R_{\nu})$ and $Q(t)$ in $t$, we get
\begin{equation}
	(\nabla[\u(\cdot,R_{\nu})Q(\cdot)])(s)=D\u(s,R_{\nu})(\D_sR_{\nu})Q(s)+\frac{1}{2}\u(s,R_{\nu})Q(s).
\end{equation}
This proves \eqref{eqn:4.15} and is the end of proving the existence result.

For the uniqueness result, with the arguments in Section 3, we note that the model \eqref{eqn:2.13}-\eqref{eqn:2.14} is equivalent to \eqref{eqn:3.3} when the initial random fields $R_{\nu}\in\D^{1,2}_{\text{loc}}(\vH)\cap\V,R_d\in\D^{1,2}_{\text{loc}}(\dH^1)\cap\dH^2$. The proof of uniqueness is then very close to that in \cite{GHZ}, \cite{Zgl}, so we omit here.
\qed

\end{document}